\documentclass[10pt]{iopart}

\usepackage{lineno}

\usepackage{a4wide}
\usepackage{parskip}
\usepackage{amsfonts}
\usepackage{url}
\usepackage{graphicx}
\usepackage{epstopdf}
\usepackage{bm}
\usepackage{fancyvrb}
\usepackage{todonotes}

\newcommand{\be}{\begin{equation}}
\newcommand{\ee}{\end{equation}}
\newcommand{\ds}{\dot{\sigma}}

\begin{document}
\title[Complex-plane singularity dynamics for blow up in a NLH]{Complex-plane singularity dynamics for blow up in a nonlinear heat equation: analysis and computation}

\author{M.~Fasondini$^1$, J.R. King$^2$, J.A.C. Weideman$^3$}

\address{$^1$ School of Computing and Mathematical Sciences, University of Leicester, 
Leicester LE1~7RH, UK }
\address{$^2$ School of Mathematical Sciences, University of Nottingham, 
Nottingham NG7~2RD, UK}
\address{$^3$ Department of Mathematical Sciences, Stellenbosch University, Stellenbosch 7600, South Africa}

\eads{\mailto{m.fasondini@leicester.ac.uk},  \mailto{John.King@nottingham.ac.uk}, \mailto{weideman@sun.ac.za}}

\begin{abstract}
Blow-up solutions to a heat equation with spatial periodicity and a quadratic nonlinearity are studied through asymptotic analyses and a variety of numerical methods. The focus is on the dynamics of the singularities in the complexified space domain.  Blow up  in finite time is caused by these singularities eventually
reaching the real axis.   The analysis provides a distinction between small and large nonlinear effects, as well as insight into the various time scales on which blow up is approached.  It is shown that an ordinary differential equation with quadratic nonlinearity plays a central role in the asymptotic analysis.  This  
equation is studied in detail, including its numerical
computation on multiple Riemann sheets, and  the far-field solutions are shown to be given at leading order  by a   Weierstrass elliptic function.
\end{abstract}


\section{Introduction}    \label{sec:intro}
\begin{figure}
	\centering
 	\includegraphics[width = 0.45\textwidth]{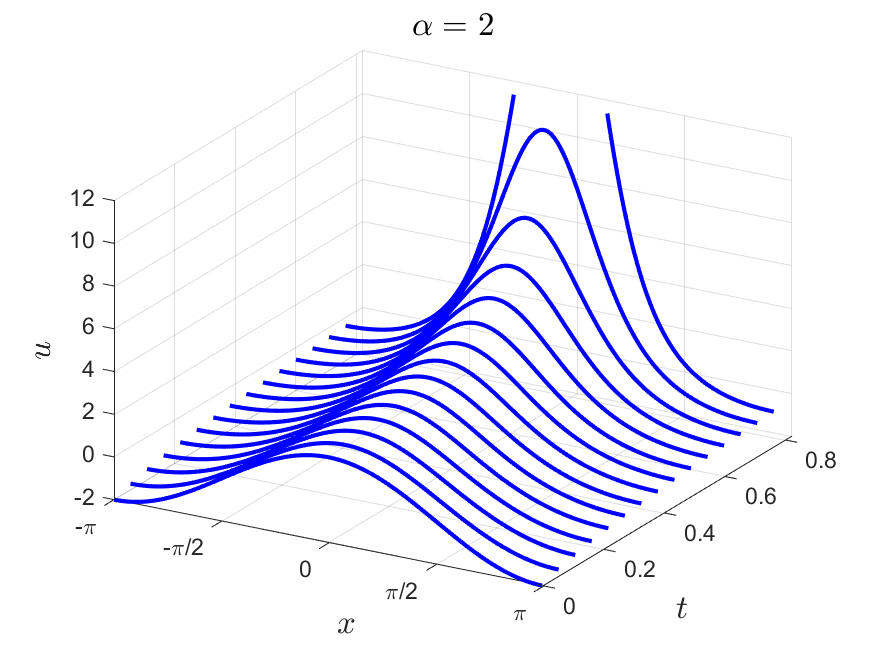} \quad \includegraphics[width = 0.45\textwidth]{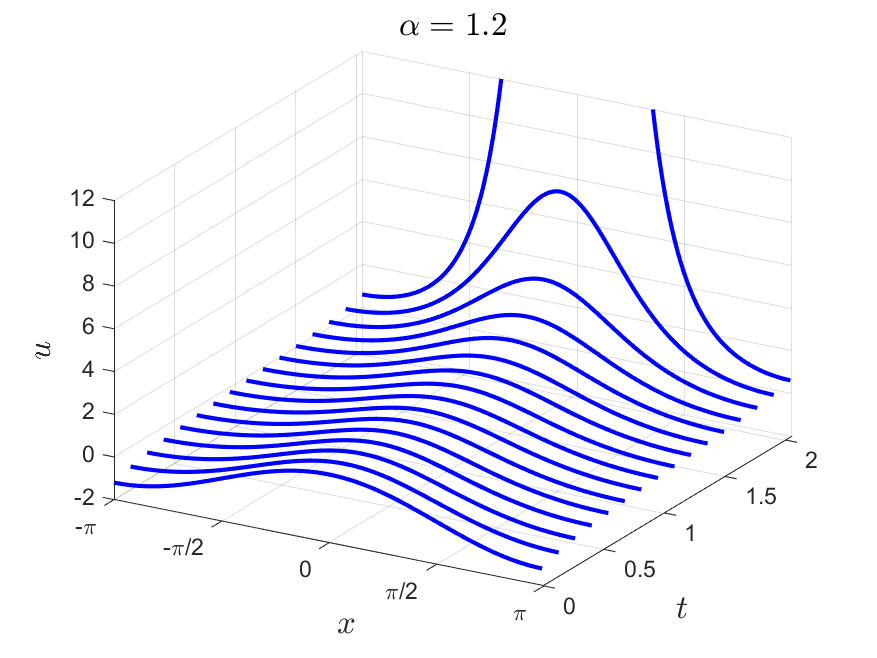}  \\
	 	\includegraphics[width = 0.45\textwidth]{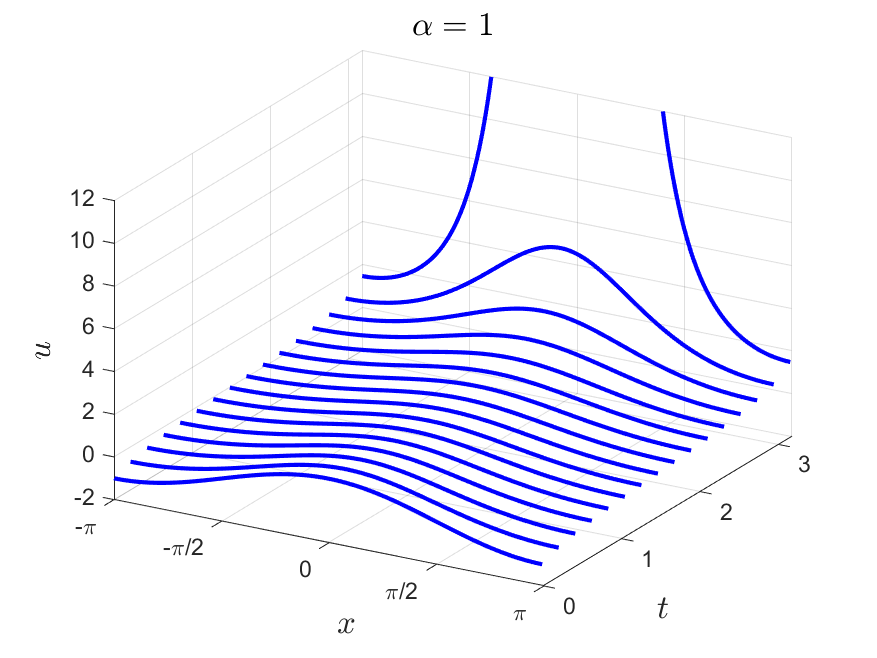} \quad \includegraphics[width = 0.45\textwidth]{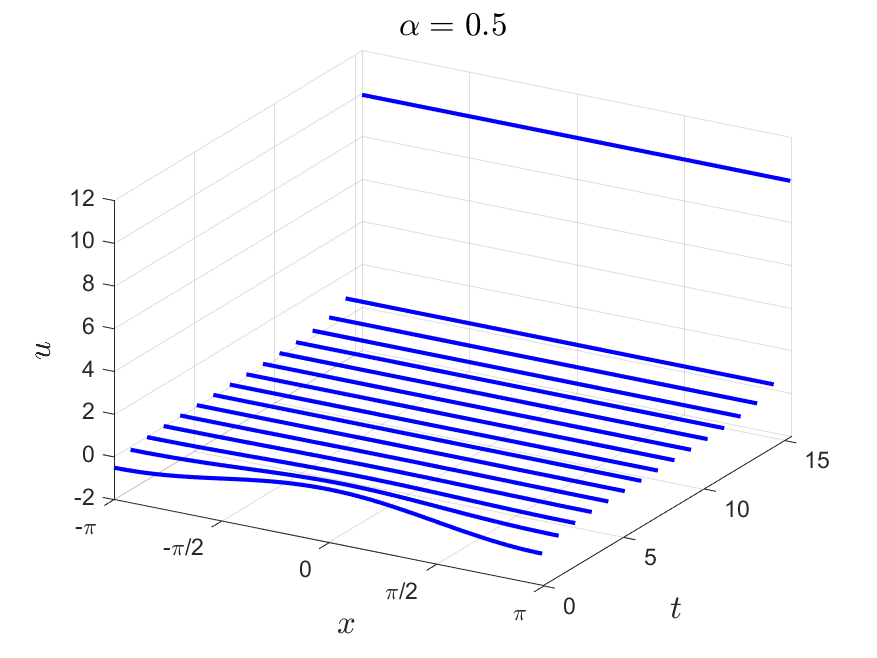}\\
	\caption{Solution profiles for blow up in the PDE~(\ref{eq:pde}) for 
		decreasing values of $\alpha$ in the initial condition~(\ref{eq:ic}).  (Note that
		the scales on the time axes are different in each case.)   In the $\alpha = 0.5$ case blow up appears to be
		uniform, however, in section~\ref{sec:smallamp} we shall show that in fact point blow up occurs at $x = 0$.
  }
	\label{fig:blowup}
\end{figure}

The nonlinear heat equation (NLH)
\be
\frac{\partial u}{\partial t} =  \frac{\partial^2 u}{\partial x^2} + u^2, 
\label{eq:pde}
\ee
is known to exhibit finite-time blow up; see~\cite{bandle98,Galaktionov2002,galakt97,quittner,samarskii}  
for reviews 
and Figure~\ref{fig:blowup} for typical
solution profiles.  The NLH is a model for reaction-diffusion processes and has appeared in numerous applications including fluid dynamics~\cite{braun2,braun1,hocking}, chemical kinetics~\cite{dold1991,Herrero93,Lacey83} and biology~\cite{jabbari2013discrete}.   The blow-up behaviour of solutions to the NLH and more general nonlinear parabolic PDEs has been studied extensively on the real line; see the previously mentioned review papers as well as~\cite{Berger88,budd98,Keller1993}.   The novel approach of the present paper
is to investigate in detail how the blow up relates to singularity dynamics of the
solution when it is  viewed as an analytic function in the complex $z$ plane, with $z = x + i y$ (only the spatial variable is complexified here). 
We adopt a combination of asymptotic and numerical methods with the
goal of obtaining a comprehensive description of the 
complex-plane behaviour leading to blow up. 

The complex-analytic viewpoint for the NLH was introduced in \cite{W03}, 
which was inspired by similar approaches to the Burgers equation
\cite{BessisFournier90,bessis1984,Senouf97a,Senouf97b} and the 
Korteweg--de~Vries equation~\cite{kruskal1974}.   The particular
example of \cite{W03} involved a real-valued, $2\pi$-periodic solution in $x$, associated with the initial condition
\be u(x,0) = \alpha \cos x,  
\label{eq:ic} \ee
with only $\alpha = 1$ considered in that paper.  Here, we 
consider all $\alpha > 0$ and will thus be able to distinguish between initially small and large nonlinear effects.   

\begin{figure}
	\centering
		\includegraphics[width = 0.35\textwidth]{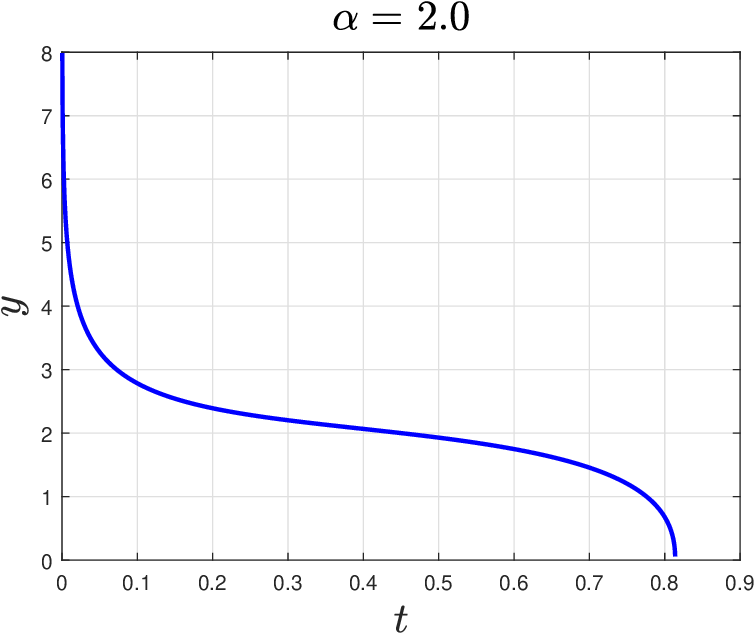} \qquad \includegraphics[width = 0.35\textwidth]{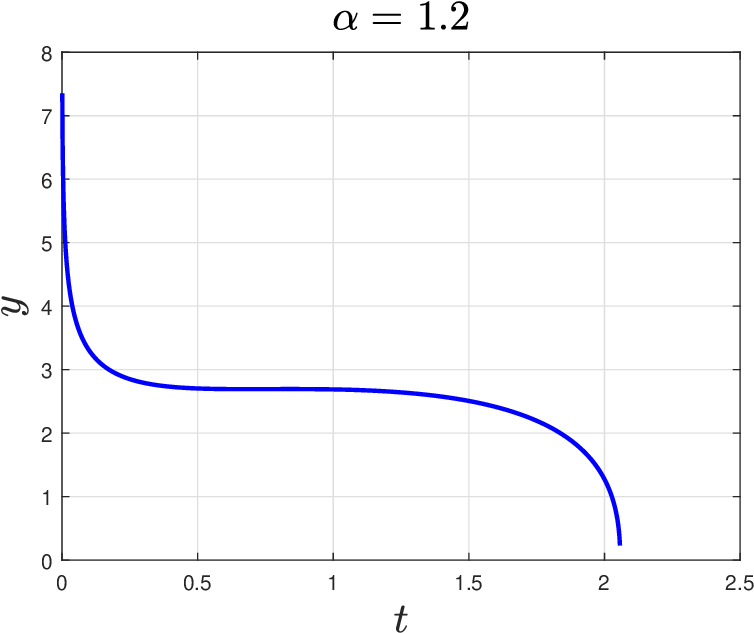}  \\
	\includegraphics[width = 0.35\textwidth]{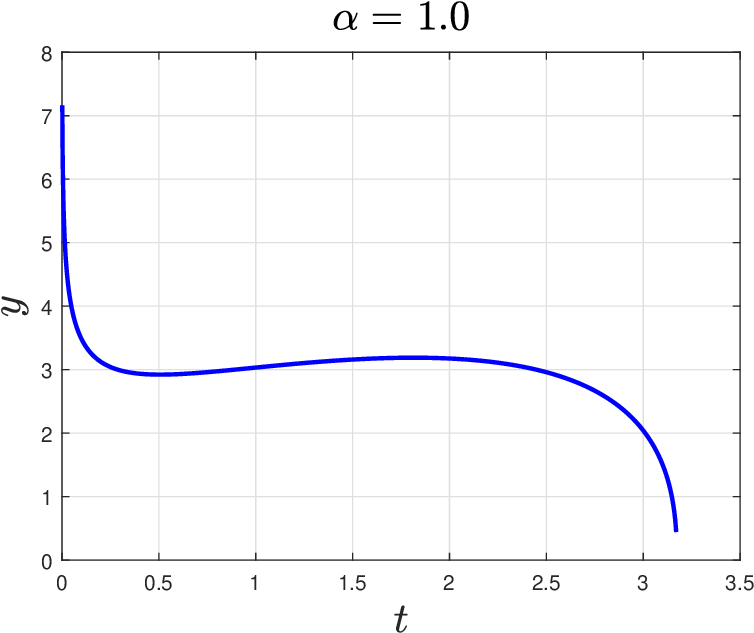} \qquad \includegraphics[width = 0.35\textwidth]{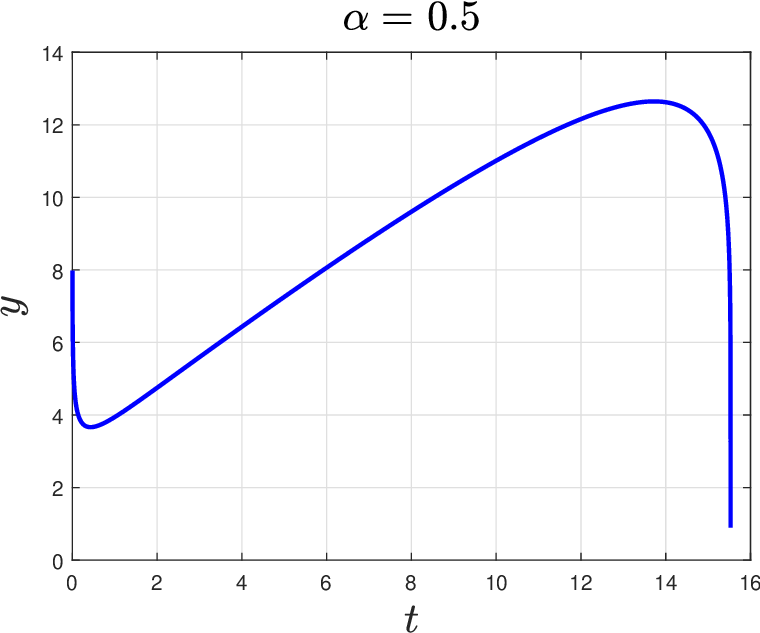} 
	\caption{Singularity locations at $\pm i y$ as functions of time,
		for the solutions shown in Figure~\ref{fig:blowup}.}   
	\label{fig:sing}
\end{figure}

Viewed in the complex $z$ plane, the initial condition (\ref{eq:ic}) is an entire function.  For small $t > 0$, singularities are born at infinity and the ones closest to the real axis, which we locate at $y = \mathrm{Im}\: z = \pm \sigma(t)$, rapidly move along the imaginary axis towards the real axis, with $u$ real on the imaginary axis for $\vert y \vert < \sigma(t)$.  Figure~\ref{fig:sing} shows the position of the closest singularity on the positive imaginary axis for the solutions in Figure~\ref{fig:blowup}. Since the solution in the upper and lower half planes are equal up to complex conjugation ($u(\overline{z},t) = \overline{u(z,t)}$, $z = x + iy$), we shall consider the solution only in the upper half plane. 

 
 

 The singularity dynamics in Figure~\ref{fig:sing} exemplify the competing effects of diffusion and nonlinearity in the NLH. For large enough $\alpha$, as in the top left frames of Figures~\ref{fig:blowup} and \ref{fig:sing}, the focussing nonlinear term dominates diffusion in the sense that the singularity approaches the real axis monotonically, albeit not at constant speed.  The solution profile on the real axis steepens as the singularity approaches the real axis, with point blow up occurring when the two closest singularities in the upper and lower half planes collide on the real axis.  If $\alpha$ is small enough, as in the bottom row of Figures~\ref{fig:blowup} and~\ref{fig:sing}, diffusion dominates initially and we see that the solution profile flattens while the singularities reverse direction after zooming in from infinity to move away from the real axis. 
 However, even as the solution flattens, its mean increases\footnote{The mean of the solution, $\langle u \rangle$, satisfies $\frac{d}{dt}\langle u \rangle = \langle u^2 \rangle$, as can be derived from~(\ref{eq:pde}) using integration by parts.} and thus eventually nonlinearity reasserts itself, the solution profile becomes steeper, the singularity changes direction again and rapidly moves towards the real axis and point blow up ensues.  For the case $\alpha \approx 1.2$ in the top-right frame of Figure~\ref{fig:sing} there is a balance between nonlinearity and diffusion in which the singularity is near stationary for a while after zooming in from infinity and before zooming  in again as point blow-up occurs.
 



\begin{figure}
    \centering
  \mbox{\includegraphics[width = 0.49\textwidth]{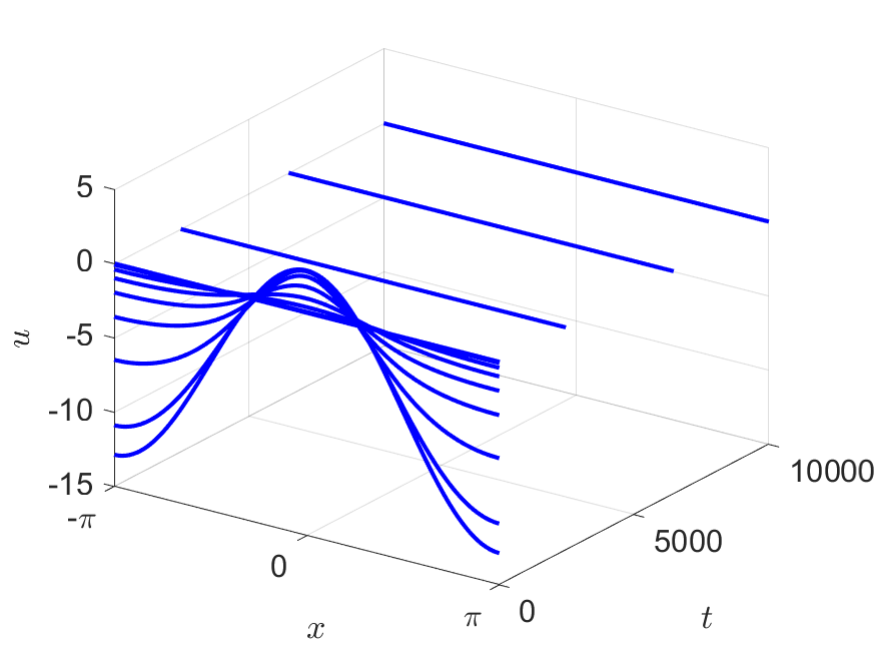}\includegraphics[width = 0.49\textwidth]{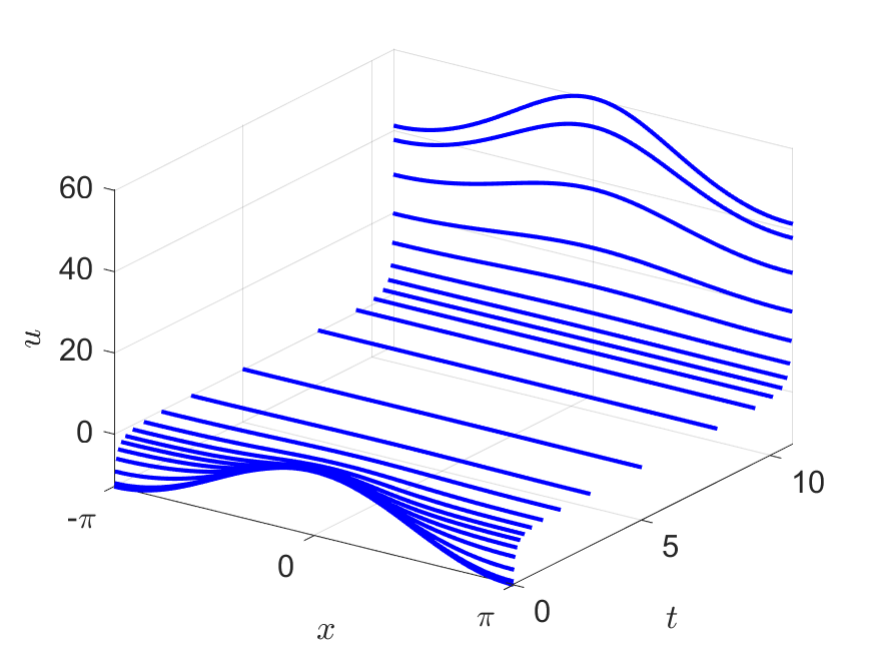}} 
  \caption{Solutions to the NLH for the initial data (\ref{eq:icab}) with $\beta = -5$, $\alpha = 7.856$ (left) and $\beta = -5$,
$\alpha = 7.892$ (right).
On the left, heat death occurs, but with a 
small perturbation in the initial condition the solution blows up, same as in Figure~\ref{fig:blowup}.
}\label{fig:reddotfigs}
   \end{figure}

\begin{figure}[h!]
    \centering
    \includegraphics[width = 0.49\textwidth]{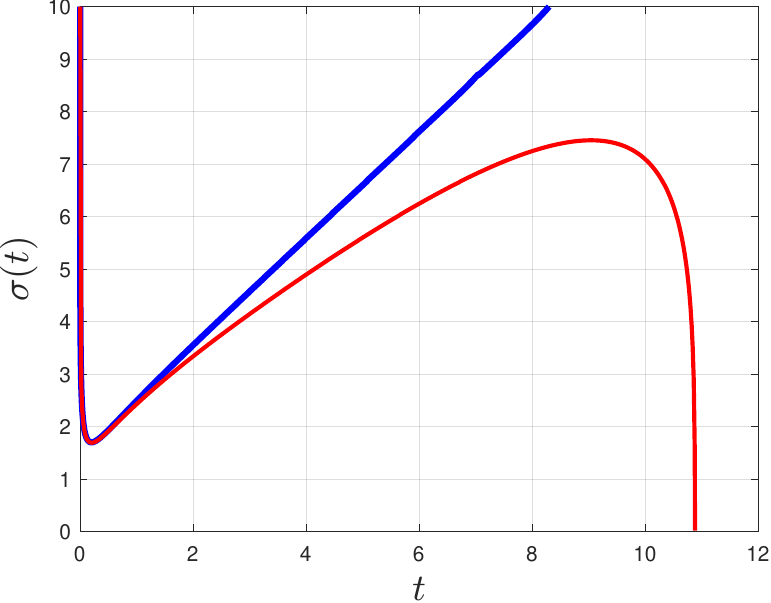}
  \caption{Singularity location of the solution in the left frame (blue curve) and right frame (red curve) of 
  Figure~\ref{fig:reddotfigs}. 
}  \label{fig:heatdeathsing}
   \end{figure}

With the initial data (\ref{eq:ic}), the solution has a point blow up at $x = 0$ for any $\alpha > 0$.  The blow up may occur at other points, however, or not at all. For example, note that $\alpha < 0$ corresponds to a translation of the solution $x \mapsto x \pm \pi$ and thus shifts the blow-up location from the maximum of the initial data (\ref{eq:ic}) at $x =0$ to $x = \pm \pi$.  On the other hand, we shall also 
consider 
\begin{equation}
    u(x,0) =  \alpha \cos x  + \beta,  \label{eq:icab}
\end{equation}
briefly in an appendix. If $\beta$ is sufficiently large negative for a fixed $\alpha$, then the solution does not blow up but extinguishes (`heat death' occurs) according to $u \sim -1/t$, $t \to \infty$; see the left frame of Figure~\ref{fig:reddotfigs}. 
 Otherwise, point blow up occurs at (i) $x = 0$ if $\alpha > 0$ and at (ii) $x = \pm \pi$ if $\alpha < 0$\footnote{On the (unstable) borderline between blow up and extinction one has $u \sim K e^{-t}\cos x$ as $t \to \infty$.}.   For an example of (i), we perturb the initial data of the heat-death solution in Figure~\ref{fig:reddotfigs} and obtain the solution in the right frame of Figure~\ref{fig:reddotfigs}.  For the heat-death solution, the singularities zoom in from infinity, change direction and then move back (linearly in time, as we shall show) to infinity, see Figure~\ref{fig:heatdeathsing}.  This figure also shows that when perturbing the heat-death solution, the singularities can switch direction a second time and move towards the real axis, which leads to blow up.

\begin{figure}
    \centering
  \mbox{\includegraphics[width = 0.49\textwidth]{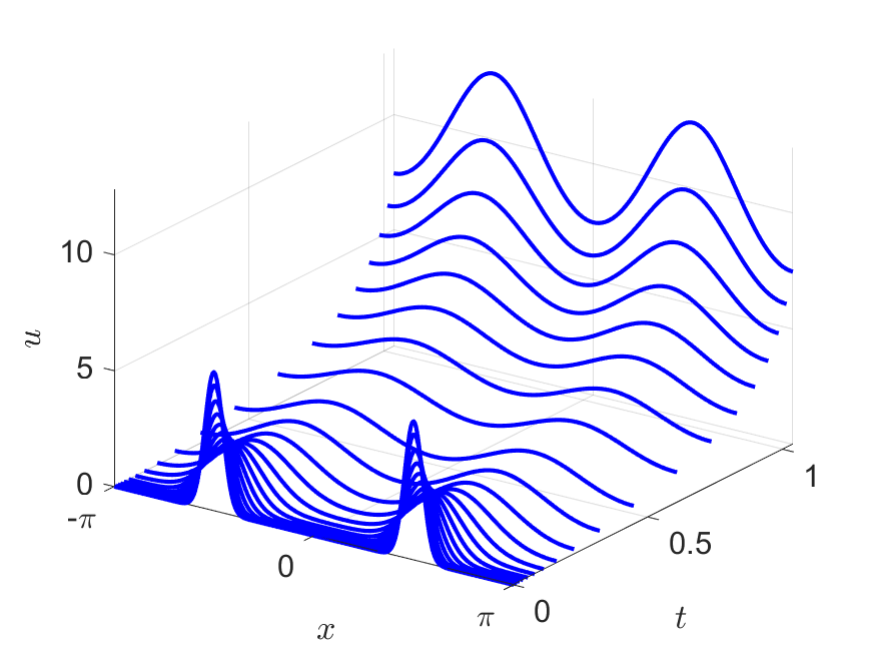}\includegraphics[width = 0.49\textwidth]{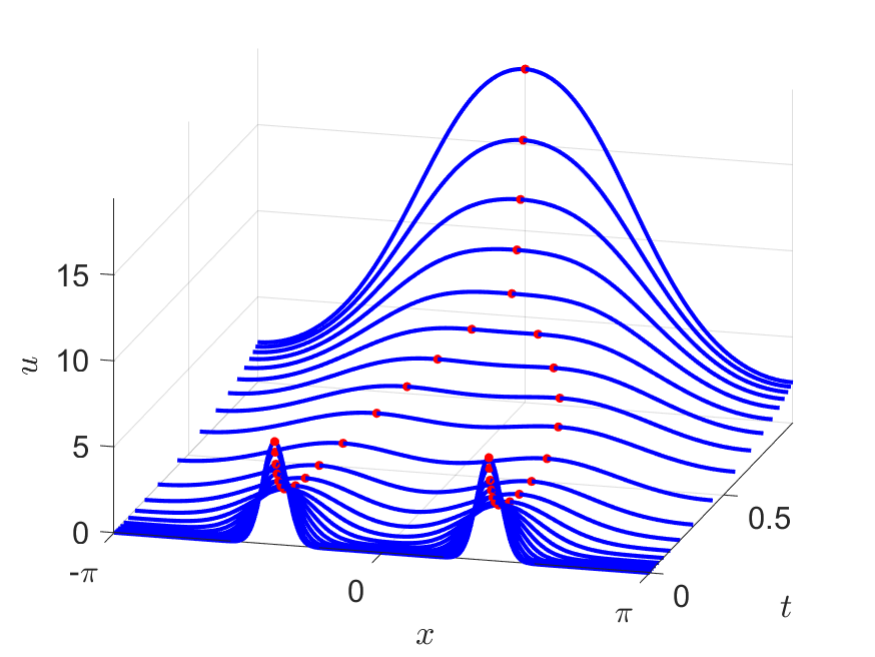}} 
  \caption{NLH solutions for two-peaked  initial data defined by $u(x,0) = \alpha \exp(\mu \cos(x + \delta )-\mu) + \alpha \exp(\mu \cos(x + \delta )-\mu)$.  For both solutions, $\alpha = 6$, $\mu = 50$, while on the left, $\delta = \pi/2$ and on the right, $\delta = 0.4\pi$.  The red dots in the right frame indicate the maxima of the solution. 
}\label{fig:2peaksreal}
   \end{figure}

In our analysis of the solution in the blow-up limit, we shall also briefly consider the blow-up scenarios, and associated singularity dynamics, for initial data with two local maxima.  Figure~\ref{fig:2peaksreal} shows two possibilities: blow up occurs at two distinct points for even initial data with two sufficiently separated and concentrated peaks (left frame) and (right frame) two maxima are sufficiently close to diffuse and combine into a single maximum, and then blow up occurs at a single point.  We shall also consider the non-generic blow up that occurs when the two maxima combine precisely at the blow-up time, which represents the borderline between the cases shown in Figure~\ref{fig:2peaksreal}.  
   
Other complex-plane studies of the NLH
have been reported in~\cite{Cho16,Takayasu22}. In both papers
the equation was 
studied numerically in the complex $t$-plane. 
Related studies reported in~\cite{fasondini2023blow} focused on the case of  nearly flat initial data
\begin{equation}
    u(x,0) = \frac{1}{\alpha - \epsilon \cos x}, \qquad 0 < \epsilon \ll \alpha.  \label{eq:icflat}
\end{equation}
The properties of blow up were investigated  but  
singularity dynamics in the  complex $x$-plane were only mentioned briefly.


The paper consists of seven sections and six appendices, 
the latter containing the bulk of the more technical analyses. 
In section~\ref{sec:fourier} we
review two numerical approaches based on Fourier analysis for locating
and classifying the singularities alluded to above.  
In section~\ref{sec:local} we 
 complement the 
numerical investigation by a local analysis that characterises the type of singularities admitted by the equation. Sections~\ref{sec:smalltime}--\ref{sec:smallamp} are devoted to the analysis and numerical verification of the dynamics of the singularities and the relation of the NLH solution in the neighbourhood of the singularities to certain nonlinear ODE solutions. The ODE solutions are studied in detail in \ref{sect:appode}. Section~\ref{sec:smalltime} and \ref{sec:appsmalltime} concern the small-time limit $t \to 0^+$ for all amplitudes $\alpha > 0$. 
Section~\ref{sec:largeamp} and \ref{sec:applargeamp} consider the large-amplitude limit, $\alpha \to \infty$, and section~\ref{sec:smallamp} and \ref{sec:appsmallamp} treat the small amplitude limit, $\alpha \to 0$.
As Figures~\ref{fig:blowup} 
 and~\ref{fig:sing} indicate, the small-amplitude singularity dynamics are the most complicated. This will be borne out by the asymptotic analysis, in which a complete description requires five distinct time scales. To complete our study, in section~\ref{sec:blowupmain} and \ref{sect:appsectblowup} we consider the solution in the blow-up limit.  
Sections~\ref{sec:smalltime}--\ref{sec:blowupmain} focus on NLH solutions subject to the initial data (\ref{eq:ic}) but the final appendix, \ref{sect:appicab}, considers NLH solutions corresponding to the more general initial condition (\ref{eq:icab}).   
Throughout we shall reuse symbols with different meanings in different sections. 


\section{Numerical method and singularity tracking} \label{sec:fourier}
The numerical solutions reported in this paper were computed by a Fourier spectral method.   Considering solutions $2\pi$-periodic in space, the approximation is based on the Fourier series
\begin{equation}
u(x,t) = \sum_{k=-\infty}^{\infty} c_k(t) e^{i k x}.  \label{eq:uFourier}
\end{equation}
When this is substituted into~(\ref{eq:pde})
an infinite dynamical system for the evolution of the coefficients $c_k(t)$ is obtained.  Upon truncation to modes $|k|\leq N$, a 
finite-dimensional system is obtained, which we integrated in time with the adaptive time-step functions
\texttt{ode45} and \texttt{ode15s} available in MATLAB.   The main advantage of these integrators is the error control that they provide. By experimentation, the number of modes, $2N+1$, was chosen sufficiently large so that all results presented here have fully converged.

For computing solutions close to blow up, a powerful strategy has been suggested in~\cite{Keller1993} and used in~\cite{fasondini2023blow}. Namely, the substitution $u = 1/v$ converts~(\ref{eq:pde}) into an equivalent PDE whose solution approaches zero rather than infinity at the blow-up point. This method can maintain high accuracy for blow-up solutions for values on the order of $u = 1/v = \mathcal{O}(1/\varepsilon)$, where $\varepsilon$ here denotes the machine precision, which can be made arbitrarily small with variable precision arithmetic but which is approximately $10^{-16}$ for IEEE double precision.   In the present paper, we  use the Fourier spectral method in the variable $v$ whenever $u$ is strictly positive, otherwise, if $u$ changes sign on $[-\pi, \pi]$, we solve the PDE in the $u$ variable.  The latter method is accurate for solution values only up to approximately $10^8$ in double precision.  If $u$ changes sign and arbitrarily large solution values are required then more complicated numerical methods such as domain decomposition, dynamical rescaling and/or adaptive mesh refinement are needed; see~\cite{Berger88,budd98}.

As for singularity tracking, 
there are two main numerical techniques.   The first is based on the examination of the rate of decay of the Fourier coefficients in~(\ref{eq:uFourier}). The 
other is based
on Fourier-Pad\'e methods for numerical analytic continuation.

The first of these methods is described in the well-known paper by Sulem~et~al.~\cite{sulem1983}.   
Suppose at a fixed time $t$ the coefficients
$c_k$ of the series~(\ref{eq:uFourier}) 
are available. 
If the singularity closest to the real axis is at $z_* = x_* + i y_*$ and 
\begin{equation}
u(z) \sim C(z - z_*)^{-\mu}, \qquad z \to z_*,  \label{eq:sing}
\end{equation}
for some constant $C$, then~\cite{sulem1983}
\be
|c_k| \sim |C| k^{\mu -1} e^{-k y_*}, \qquad k \to +\infty. 
\label{eq:coeff}
\ee
Given values of $c_k$ for a range $k \gg 1$, the values of 
$\mu$ and $y_*$ can be estimated  by a linear least squares fit, after taking logarithms in (\ref{eq:coeff}).  The singularity locations shown
in Figure~\ref{fig:sing} were computed by this method. 







 

\begin{figure}[htb]
	\centering
	\includegraphics[width = 0.4\textwidth]{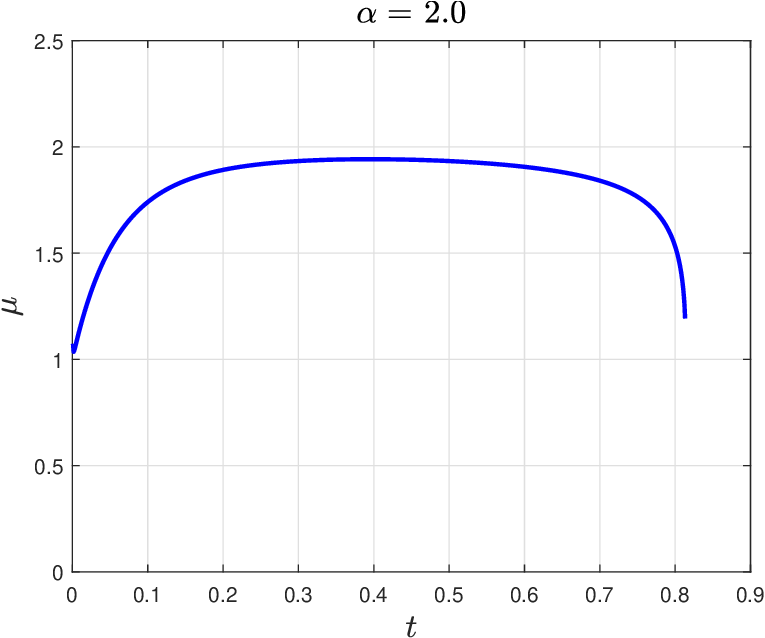} \qquad
	\includegraphics[width = 0.4\textwidth]{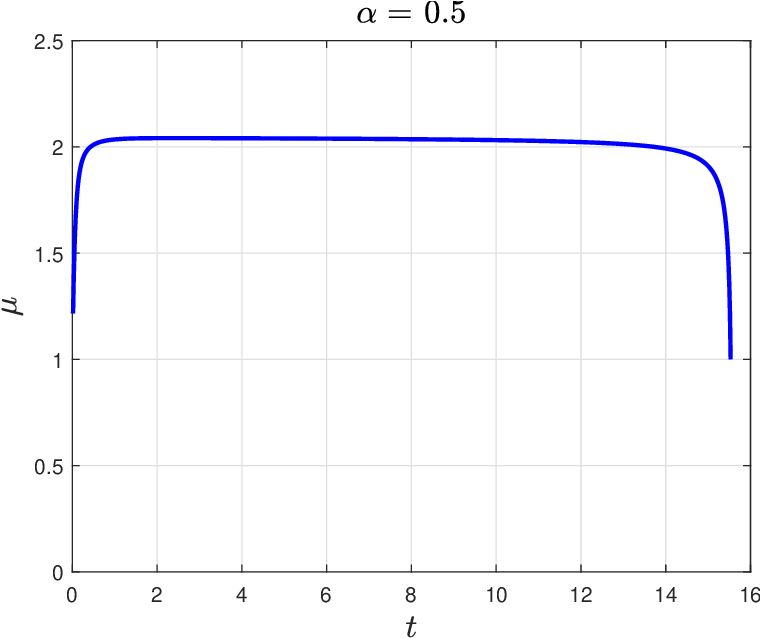}
	\caption{Estimated values of $\mu$ as defined in (\ref{eq:coeff}).  Only two of the cases of Figures~\ref{fig:blowup} and~\ref{fig:sing} are displayed here, 
	 the others being qualitatively similar.
  The numerical results single out the value $\mu=2$, which suggests poles of order $2$ as the nearest singularities to the real axis. The true nature of the singularities is more complicated, involving a logarithmic branch point,   as the analysis of section~\ref{sec:local} will show.}  
	\label{fig:prefactor}
\end{figure}

In the next section we proceed with a theoretical analysis of the
singularities of (\ref{eq:pde}), but we assume for now an expression
of the form (\ref{eq:sing}). The least-squares procedure then produces
the estimates for the exponent $\mu$ shown in Figure~\ref{fig:prefactor}. These results suggest that the singularities could be poles of order two, except in an intermediate asymptotic sense for initial
times and as blow up is approached. The asymptotic analyses in the next sections will clarify these numerical estimates.

The other method for singularity tracking is based on Fourier--Pad\'e methods as considered for the NLH (and other PDEs) in~\cite{W03}.
   The Fourier series~(\ref{eq:uFourier}) is converted
   to  rational trigonometric form,
   which can be continued into the complex plane to some strip around the real axis.
   The advantage of this method 
   over the method~(\ref{eq:sing})--(\ref{eq:coeff}) is that it often gives information further into the complex plane,  beyond the singularities closest to the real axis.   In a further improvement this method was  recently extended to 
quadratic Fourier--Pad\'e, 
which incorporates
a square-root singularity into the approximant in an attempt to capture branch point
singularities more accurately~\cite{Fasondini2019}.

\begin{figure}
	\centering
	\includegraphics[width = 0.2\textwidth]{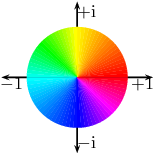}\\
	\includegraphics[width = 0.48\textwidth]{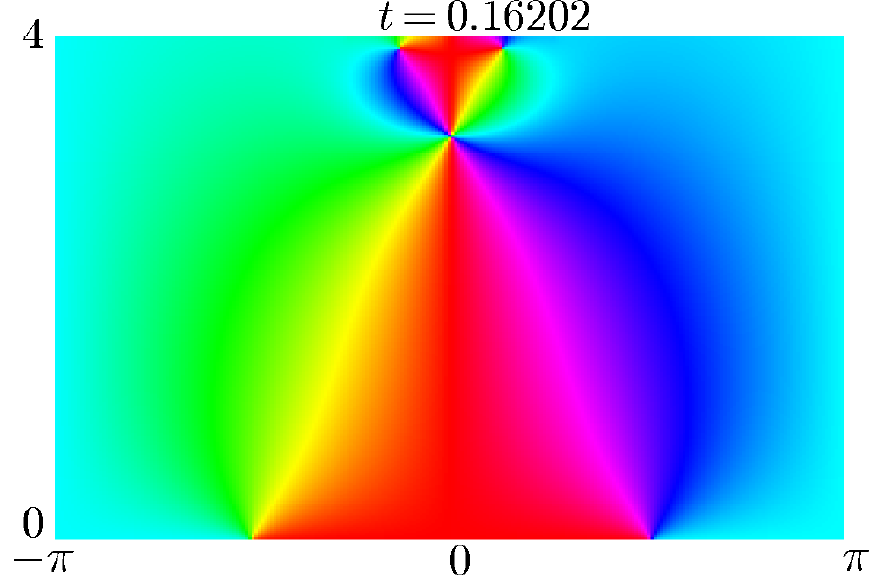} \quad
	\includegraphics[width = 0.48\textwidth]{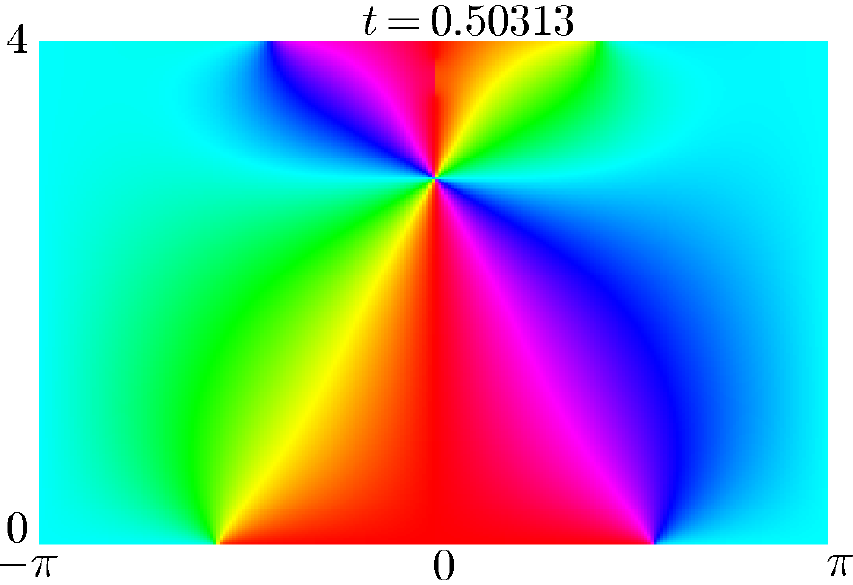}  \\
	\includegraphics[width = 0.48\textwidth]{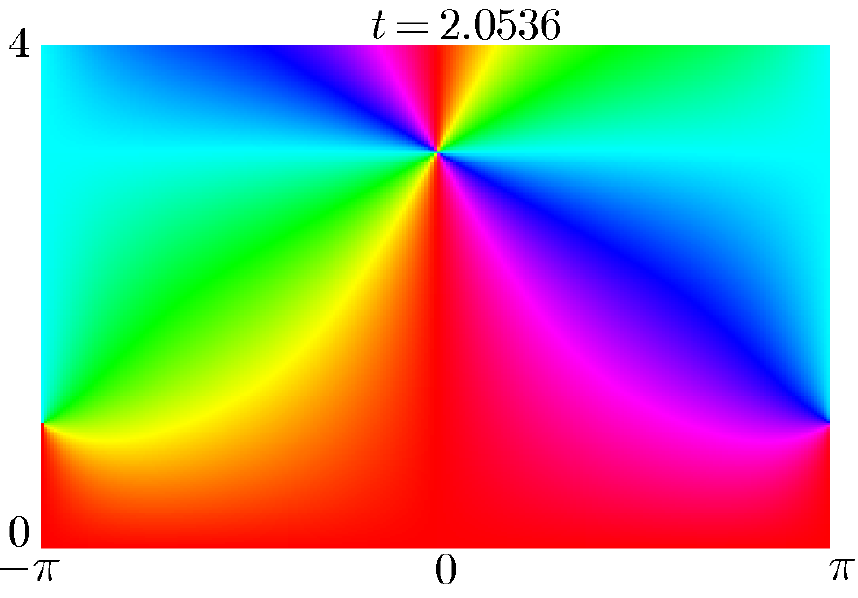} \quad
	\includegraphics[width = 0.48\textwidth]{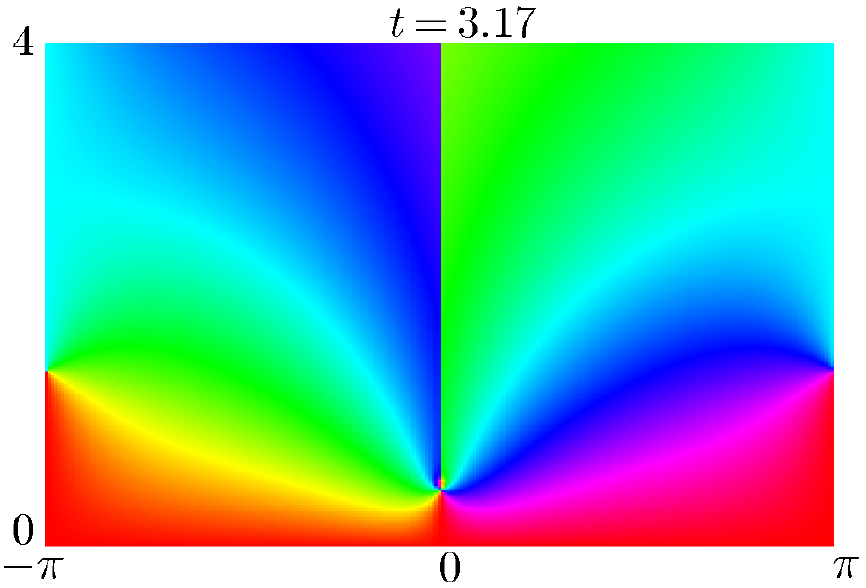}
	\caption{Approximate complex phase plots of the $\alpha = 1$ solution  
		when continued from the real line into the upper half of the complex $z$-plane. The colours indicate the phase $\phi(z,t) \in [-\pi, \pi)$ of the solution (where $u(z,t) = \vert u(z,t) \vert \exp i \phi(z,t)$) according to the colour wheel at the top (taken from~\cite {NIST:DLMF}).  
		In the first frame the singularity travels down the imaginary axis, until it 
		momentarily stops and turns around at the location in the second frame.   
		It then moves away from the
		real axis until it reaches the position of the third frame, at which time
		it turns around once more and
		rushes onto the real axis and blow up occurs.     
		The fact that the colours go twice around the colour wheel when encircling
		the singularity suggests that the dominant contribution is a pole of order two.    In the last frame
	the discontinuity in the phase on the imaginary axis is suggestive of a branch point, 
 however, 
as will be discussed in the next section.
  (The interpretation of complex phase plots is discussed in~\cite{wegert2012visual} and plotting software can be found at~\cite{Wegert}.)
   } 
	\label{fig:phase}
\end{figure}

In Figure~\ref{fig:phase} we show 
phase plots for the case $\alpha = 1$,
computed by this quadratic Fourier--Pad\'e method. (It provides a different view of the singularity dynamics shown in the third frame of Figure~\ref{fig:sing}.)   This improves on the figure   given in~\cite[Fig.~5.2]{W03}, which was computed with the standard Fourier--Pad\'e approach.
By examining the phase plots, it again appears that the dominant singularity is a pole of order two, as was evident in Figure~\ref{fig:prefactor} as well.   As the singularity approaches the real axis, however, the characteristics of a branch point become evident.  The analysis of the 
next section will clarify this. 

\section{Local analysis of the singularities}  \label{sec:local} 
In this section we examine more closely the nature of the singularities 
described in sections~\ref{sec:intro} and~\ref{sec:fourier}. 
Locating a singularity in the complex plane at 
$x = i\sigma(t)$
and
setting 
$x = i\sigma(t)+\zeta$ 
gives
\be
\frac{\partial u}{\partial t} 
- i\dot{\sigma} \frac{\partial u}{\partial \zeta}
= \frac{\partial^2 u}{\partial \zeta^2} + u^2, 
\label{eq:local}
\ee
where the dot represents differentiation with respect to $t$.
The right-hand side dominates the local behaviour, so that as
$\zeta \to 0$,
\begin{equation*}
u \sim -\frac{6}{\zeta^2}.
\end{equation*}
Setting
\begin{equation*}
u = -\frac{6}{\zeta^2} + V(\zeta,t),
\end{equation*}
and linearising, at leading order the `complementary function' $V$ satisfies the Euler equation
\begin{equation*}
\frac{\partial^2 V}{\partial \zeta^2} - \frac{12}{\zeta^2}V = 0,
\end{equation*}
so that
\begin{equation*}
V \sim A(t)\zeta^{-3} + B(t) \zeta^4.
\end{equation*}
Self-consistency requires $A = 0$ (this contribution being associated with the $\zeta$-translation-invariance of (\ref{eq:local})),
but $\sigma(t)$ and $B(t)$ are arbitrary in terms of the local analysis, representing the two degrees of freedom expected of a generic singularity in this second-order problem. 
Reinstating the intervening terms in the local expansion about the singularity in (\ref{eq:local}), one finds that
\be
\hspace{-2cm}u \sim -\frac{6}{\zeta^2} + \frac{6 \, i\ds}{5 \, \zeta} - \frac{1}{50} \ds^2 + a(t) \zeta + b(t)\zeta^2 + c(t) \zeta^3 + d(t) \zeta^4 \log \zeta + B(t) \zeta^4, \qquad \zeta \to 0,
\label{eq:expansion}
\ee
where
\begin{eqnarray*}
& \hspace{-2cm} a(t) = -{\frac { i}{250}}\dot{\sigma}^{3}-\frac{i}{10}\ddot{\sigma},  \qquad \qquad \qquad \:
b(t) = {\frac { \dot{\sigma}  \left( 7\, \dot{\sigma}^{3}+190\,\ddot{\sigma}  \right) }{5000}},\\
& \hspace{-2cm} c(t)= {\frac {79i }{75000}}\dot{\sigma}^{5}+{\frac {229i }{7500}}\dot{\sigma}^{2}\ddot{\sigma} +{\frac {i }{60}}\ddot{\sigma}, \qquad
d(t)= {\frac {18\, }{21875}}\dot{\sigma}^{6}+{\frac {108\,  }{4375}}\dot{\sigma}^{3}\ddot{\sigma}+{\frac {16\, }{875}}\dot{\sigma}\ddot{\sigma}+{\frac {6\, }{875}}\ddot{\sigma}^2.
\end{eqnarray*}

Note that the presence of the $d(t) \zeta^4 \log \zeta$ term implies that the singularity is actually a branch point. However, the double-pole nature of the dominant term in (\ref{eq:expansion}) significantly precedes this in the local expansion. This explains why both singularity tracking methods
in section~\ref{sec:fourier} singled out this
type of singularity.  

The expansion might also explain why the method
~(\ref{eq:sing})--(\ref{eq:coeff}) suggests 
simple pole singularities near $t=0$ and near blow up; see Figure~\ref{fig:prefactor}. This is because the second term on the right in~(\ref{eq:expansion}) has a residue proportional to $\dot{\sigma}$,
and this value is large initially and again near blow up, as can be confirmed in 
Figure~\ref{fig:sing}.  A more detailed asymptotic analysis is required fully to clarify the matter; see below.

The method based on the quadratic Fourier-Pad\'e method is able, at least ultimately, to classify the singularity as a branch point; see  the last frame of Figure~\ref{fig:phase}. The approximant is a two-valued function, so it is unable to capture the logarithmic 
singularity perfectly, but at least the presence of a branch point is indicated strongly.




\section{Small-time limit: singularity dynamics and  singularity structure}~\label{sec:smalltime}  
 In the first of our appendices, \ref{sec:appsmalltime},   it is shown that for  the initial condition (\ref{eq:ic}) and for small values of $t$, the NLH solution on the real axis  is
\begin{equation*}
u(x,t) = \alpha \cos x + t \Big(-\alpha \cos x + \frac12 \alpha^2 (1+\cos 2x)\Big) + \mathcal{O}(t^2), \qquad t \to 0. 
\end{equation*}
Furthermore, it is shown that the singularity closest to the real axis on the positive imaginary axis is located  at $z = i\sigma(t) \sim i\log(2/(\alpha t))$, $t \to 0$. In the vicinity of $i\sigma(t)$, the solution to the NLH in the complex plane is given by 
\begin{equation}
\hspace{-2.5 cm}u(z,t) \sim \frac{\phi(\zeta)}{t^2}, \qquad z = 
i\Big(\log(2/(\alpha t)) + 2t\log (1/t) - (\zeta+1)t \Big), \qquad t \to 0,\qquad \zeta t = o(1),  \label{eq:NLHcomplexsmallt}
\end{equation}
where $\phi$ is the solution to the nonlinear ODE 
\begin{equation}
 {\frac {d^{2}\phi}{d\zeta^{2}}}      -  {\frac{d \phi }
	{d\zeta}}   = \phi^{2},  
 \label{eq:V0m}
\end{equation}
subject to the asymptotic  condition
\begin{equation} 
 \phi = \frac{1}{\zeta} + \frac{2\log(\zeta)}{\zeta^2} + o\left(\zeta^{-2}\right),  
 \qquad  \zeta \to \infty.
 \label{eq:V0mbc}
\end{equation}
Let $\zeta_*$ be the first singularity that is encountered on the real axis as the ODE problem (\ref{eq:V0m})--(\ref{eq:V0mbc}) is integrated from $\infty$, then it follows from (\ref{eq:NLHcomplexsmallt}) that a small-time approximation to the singularity location is given by $z = i\sigma(t)$, where 
\begin{equation}
  \sigma(t) \sim \log(2/(\alpha t)) + 2t\log (1/t) - (\zeta_*+1)t. \label{eq:sigmasmalltime}
\end{equation}

In subsequent sections, we shall find that the ODE (\ref{eq:V0m}) subject to (\ref{eq:V0mbc}) or 
variations thereof
arises repeatedly in our analysis of NLH solutions in the complex plane. Consequently, a detailed analysis of these ODE solutions is given 
in~\ref{sect:appode}.  

\begin{figure}[h]
\centering
    \hspace*{-0cm}\includegraphics[width = 0.85\textwidth]{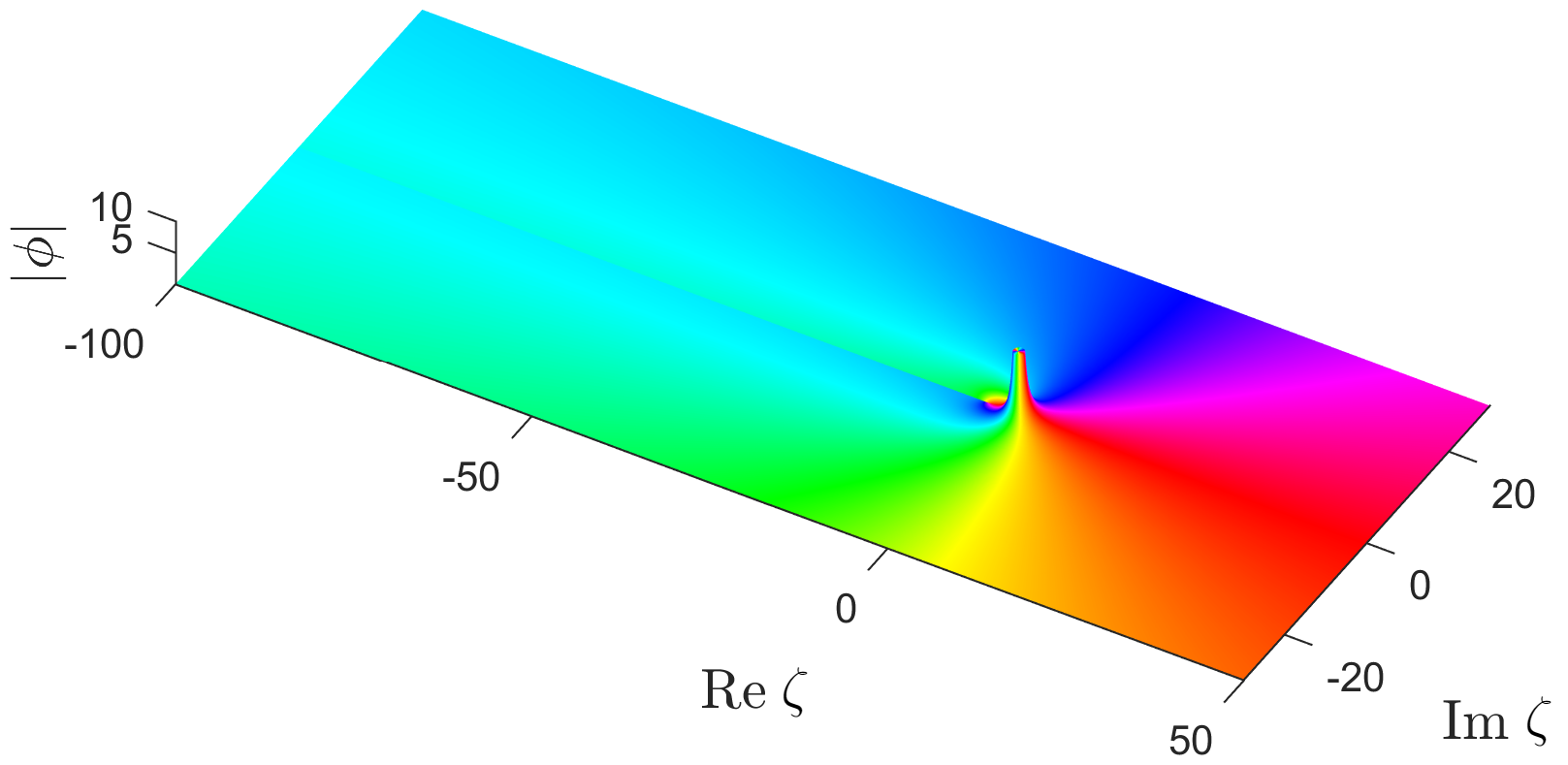}
    \caption{The solution to (\ref{eq:V0m}) satisfying (\ref{eq:V0mbc}). The modulus is shown as the height, and the phase is represented according to the color wheel at the top of Figure~\ref{fig:phase}.   The plot reveals a branch cut along the negative real axis originating from a singularity at $\zeta_* = 0.05695$.  The initial conditions used for this solution are given in (\ref{eq:ODE2ics}) and were computed with more than 800 terms of the asymptotic expansion (\ref{eq:V0mbc}), which are shown in Figure~\ref{fig:asexpterms}.  }\label{fig:ODE_complex_phase_modulus}
\end{figure}

Figure~\ref{fig:ODE_complex_phase_modulus} shows the solution to (\ref{eq:V0m})--(\ref{eq:V0mbc}), revealing that $\zeta_* \approx 0.05695$.  In \ref{sect:appode}, higher-order terms in the asymptotic expansion (\ref{eq:V0mbc}) are derived in order to compute the initial condition to an accuracy on the order of machine precision. 
The solution to the ODE (\ref{eq:V0m}) is then computed in the complex $\zeta$-plane using the `pole field solver,' a method based on an adaptive Pad\'e one-step method~\cite{fasondini2017,fornberg2011}.

\begin{figure}[h!]
	\centering
	\begin{minipage}[c]{0.52\textwidth}
\includegraphics[width = 0.95\textwidth]{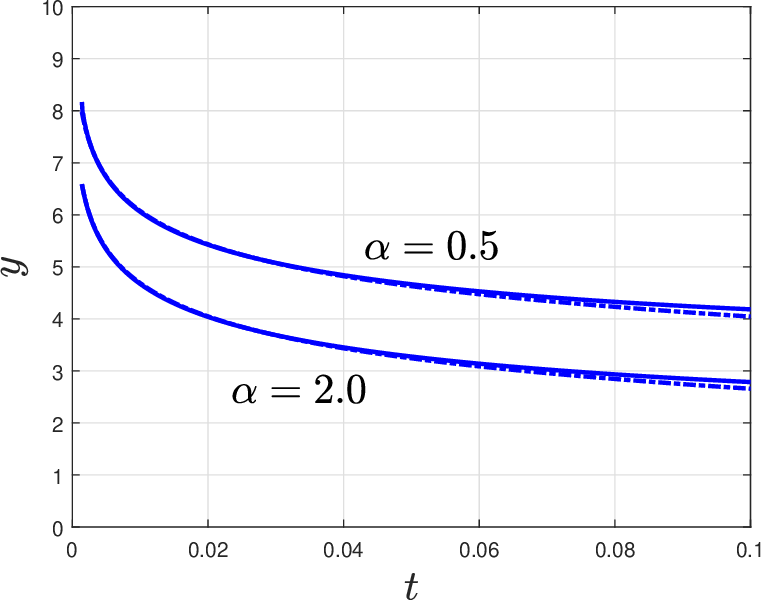}  
\end{minipage}
	\begin{minipage}[c]{0.45\textwidth}
	\begin{tabular}{c|cc|cc}
		& \multicolumn{2}{c|}{$\alpha = 2.0$}  & \multicolumn{2}{c}{$\alpha = 0.5$} \\
	$t$ &   (\ref{eq:coeff})   &   (\ref{eq:sigmasmalltime}) &  
	 (\ref{eq:coeff})     &   (\ref{eq:sigmasmalltime}) \\ \hline
	0.02 & 4.04   &  4.05  &  5.43 &  5.43  \\
	0.04 & 3.45   &  3.43  &  4.84 &  4.82 \\
	0.06 & 3.13   &  3.09  &  4.53 &  4.47 \\
   0.08 & 2.93   &  2.85  &  4.33 &  4.23 \\
   0.10 & 2.79   &  2.66  &  4.18 &  4.04  \\	
	\end{tabular}
\end{minipage}
\caption{
The singularity locations $\pm i y$ as functions of time as estimated
by the procedure of section~\ref{sec:fourier}, which is based on eq.~(\ref{eq:coeff}) 
(solid line) and the estimate 
of eq.~(\ref{eq:sigmasmalltime}) (dashed line).
}  \label{fig:firstphase}
\end{figure}

Figure~\ref{fig:firstphase} shows the estimate (\ref{eq:sigmasmalltime}) compared to the numerical estimate of the singularity location described 
by~(\ref{eq:sing})--(\ref{eq:coeff}).
Both the graphical and the numerical comparisons confirm that the asymptotic and numerical estimates match well for small $t$. As to be expected, the accuracy of the asymptotic estimate deteriorates as~$t$ increases.

The asymptotic analysis   in~\ref{sec:appsmalltime} provides another perspective
on why the singularity exponent in Figure~\ref{fig:prefactor} increases from roughly $1$ to $2$ as $t$ increases from small to intermediate values.  It is shown in (\ref{eq:riccati}) that in the small-time limit, the leading-order behaviour of the singularity (after a single rescaling) is that of a simple pole. It is this leading-order simple pole nature of the singularity that is detected by the numerical method.
It is only after a second rescaling that one arrives at an equation, viz.~(\ref{eq:V0}), whose singularity type matches that of the NLH. 


The type of singularities admitted by the ODE (\ref{eq:V0m}) are (unsurprisingly) of the same type as those of the NLH since if one assumes a singularity of the ODE is located at $\zeta_*$, then the local expansion takes the form
\begin{equation}
\hspace{-2 cm}\phi(\zeta+\zeta_*) \sim \frac{6}{\zeta^2} + \frac{6}{5\zeta} - \frac{1}{50} + \frac{\zeta}{250} - \frac{7\zeta^2}{5000} + \frac{79\zeta^3}{75000} + \frac{18}{21875}\zeta^4\log \zeta + b \zeta^4, \qquad \zeta \to 0, \label{eq:philocexpm}
\end{equation}
where $b$ is an arbitrary constant.  This shows that the ODE (\ref{eq:V0m}) does not possess the Painlev\'e property since it has branch point singularities whose locations are dependent on the initial conditions (i.e., movable branch points). This also confirms that the essential singularity that is present in the initial condition at $\infty$, and which (\ref{eq:expansion}) indicates  to be incompatible with the NLH, is transformed for any $t>0$ into singularities of the compatible form (\ref{eq:expansion}).  







\begin{figure}[h]
\centering
    \includegraphics[width = 1\textwidth]{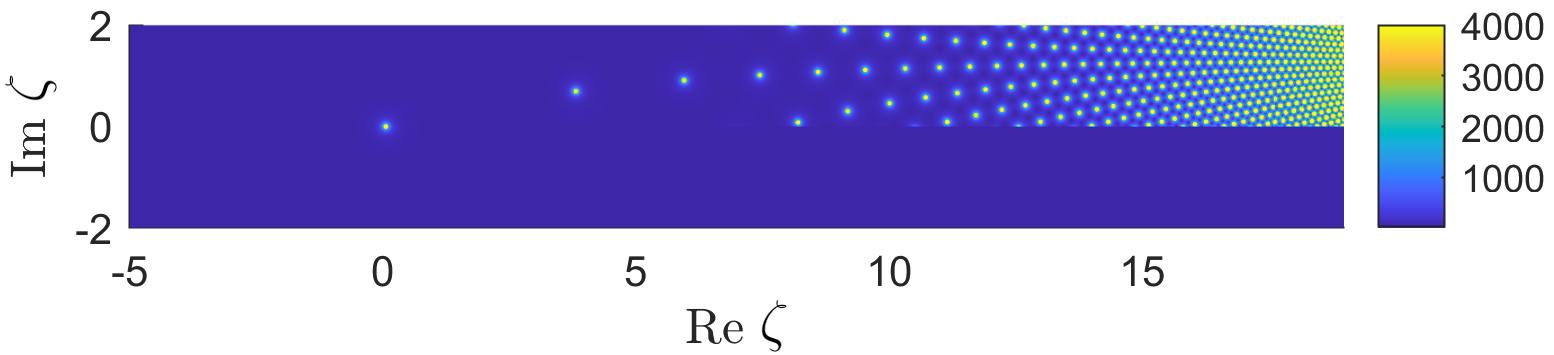}\\
   \includegraphics[width = 0.25\textwidth]{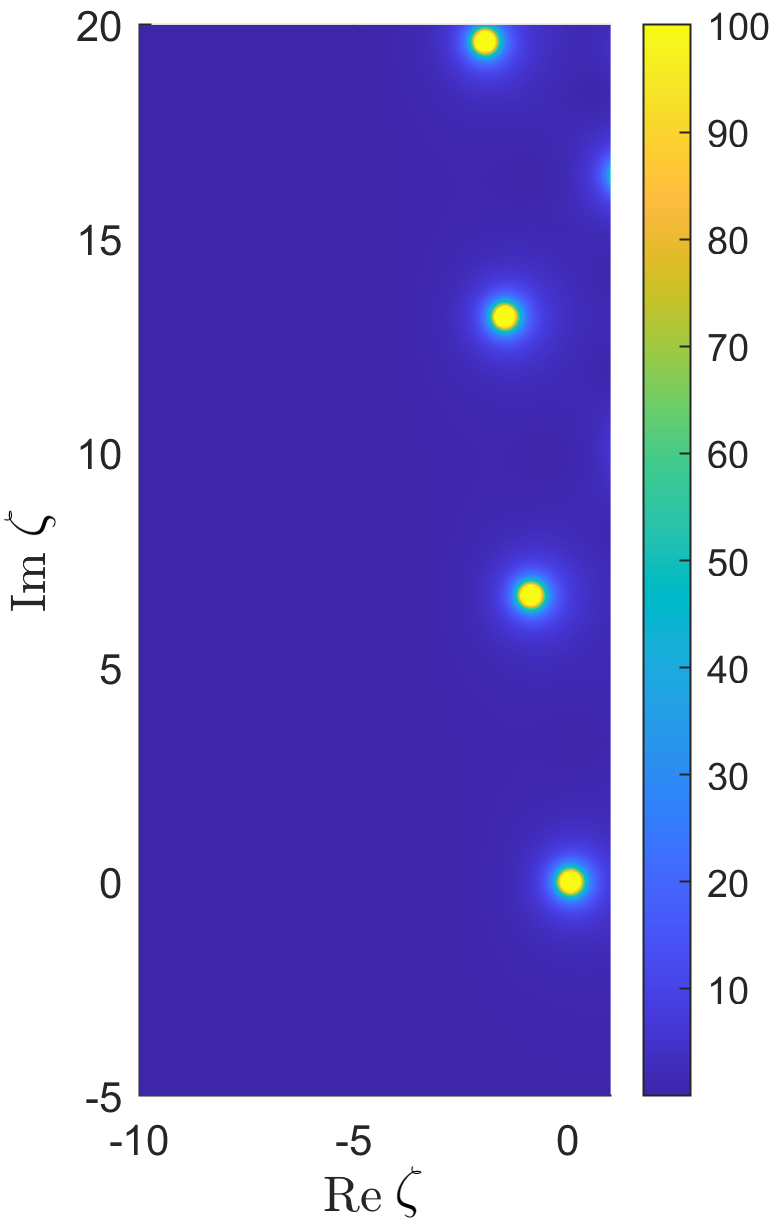} \includegraphics[width = 0.7\textwidth]{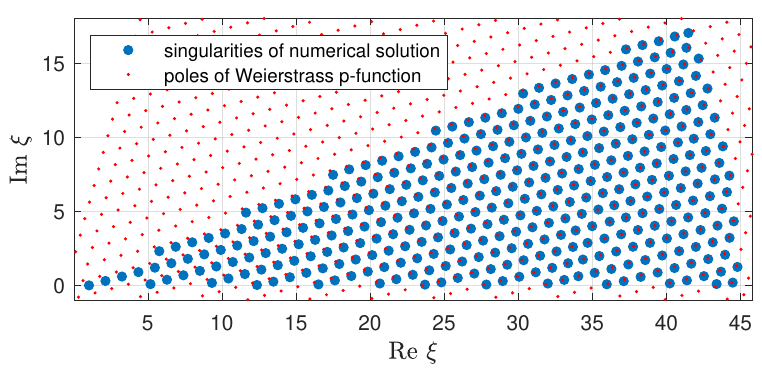}
    \caption{ Top and bottom-left: Modulus plots of  the solution in Figure~\ref{fig:ODE_complex_phase_modulus} integrated clockwise around the branch point at $\zeta_* =0.05695$, hence the solution on the lower half-planes in these figures is on the principal Riemann sheet (the sheet shown in Figure~\ref{fig:ODE_complex_phase_modulus})  and the upper half-plane lies on the next Riemann sheet.
    Bottom: 
    Comparison of the pole locations (\ref{eq:weierpoles}) of the Weierstrass function (with $\alpha = r\exp(i\theta)$, $r \approx 0.2087$, $\theta \approx 0.2524$ and $\zeta_0 \approx -0.5113 + 0.03149i$) and the singularity locations on the upper half-plane in the top frame of the figure after
    mapping to the $\xi$-plane ($\xi = e^{\zeta/5}$).  The bottom-left
    figure illustrates how the singularities first arise in the neighbourhood of the anti-Stokes line identified in~\ref{sec:appivp2}.
    }
\label{fig:ODE_complex_mod_sheet2_rotate}
\end{figure}

Figure~\ref{fig:ODE_complex_phase_modulus} suggests that the branch point singularity at $\zeta_*$ is the only singularity of $\phi(\zeta)$ on the principal Riemann sheet of the solution to (\ref{eq:V0m})--(\ref{eq:V0mbc}). 
It is also of interest to examine if, and where, other singularities might occur.
Integrating clockwise around $\zeta_*$ onto the next Riemann sheet, we obtain the solution shown in Figure~\ref{fig:ODE_complex_mod_sheet2_rotate}, which reveals the presence of a multitude of singularities.  
Such singularities could, at least in principle, subsequently move onto the principal Riemann sheet and influence the real-time behaviour; however, such effects seem not to be of importance in the current context.

\ref{sec:appivp2} shows that, to leading order, the solution for $\zeta \to \infty$ on the second Riemann sheet is expressible in terms of the (equianharmonic) Weierstrass elliptic function in the variable $\xi = e^{\zeta/5}$ (see (\ref{eq:apptrans1}) and (\ref{eq:appvsol})).  Hence, in the bottom frame of Figure~\ref{fig:ODE_complex_phase_modulus}, we find that in the $\xi$-plane the far-field singularities on the second Riemann sheet lie approximately on the same  lattice as the singularities of the Weierstrass elliptic function.

The asymptotic result (\ref{eq:NLHcomplexsmallt}) and Figure~\ref{fig:ODE_complex_mod_sheet2_rotate} suggest that the essential singularity at $\infty$ that is present in the initial condition is instantaneously transformed into infinitely many singularities of the form (\ref{eq:expansion}) for $t > 0$ that lie on a non-compact (infinitely sheeted) Riemann surface.

\section{Large-amplitude initial conditions}\label{sec:largeamp}
\begin{figure}
	\centering
\mbox{	
}\\
\mbox{	
\includegraphics[width = 0.495\textwidth]{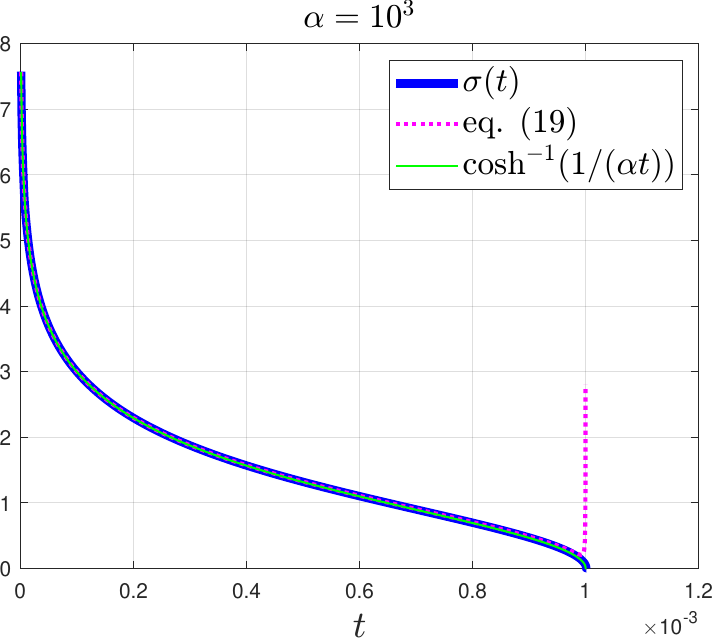}
\includegraphics[width = 0.495\textwidth]{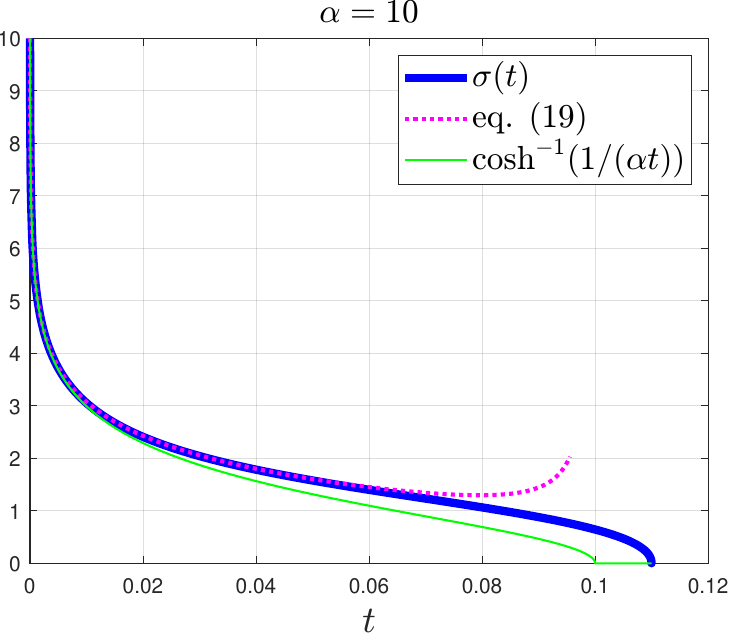}
}
	\caption{
 The singularity locations $\pm i\sigma$ as a function of time, as estimated by the numerical procedure of section~\ref{sec:fourier} (least squares) and the asymptotic estimates (\ref{eq:singloclargeamp}) and (\ref{eq:sigmalargeamp}). As can be expected, the higher-order estimate (\ref{eq:sigmalargeamp}) loses accuracy close to the blow-up time, as it should do; however, in the right frame it is clear that this estimate is more accurate for intermediate times away from blow up than the leading-order approximation (\ref{eq:singloclargeamp}).
	} 
	\label{fig:largeamp}
\end{figure}

In \ref{sec:applargeamp} it is shown that for the initial condition (\ref{eq:ic}) with large amplitude, a leading-order approximation to the solution on the real line is, for $t = \mathcal{O}(1/\alpha)$
\begin{equation}
u \sim \frac{\alpha \cos x}{1 - \alpha t \cos x}, \qquad \alpha \to \infty.  \label{eq:ulargeampreal}
\end{equation}
Hence, to leading order, the singularity locations are $x = \pm i y = \pm i \sigma(t)$, where
\be
\sigma(t) \sim \cosh^{-1}(1/(\alpha t)) = \log(1/(\alpha t)) + \log(1 + (1 - \alpha^2 t^2)^{1/2}), \qquad \alpha \to \infty. \label{eq:singloclargeamp}
\ee
The formula on the right makes it clear that as $t \to 0$,  this estimate is  consistent with the result of the small-time analysis (\ref{eq:sigmasmalltime}).
The singularities move along the imaginary axis and collide with the real axis at $x=0$ for $t = t_c \sim 1/\alpha$.  The motion of the singularities is monotonically towards the real axis, albeit not at constant speed.    Graphs of the singularity locations as functions of time are shown in Figure~\ref{fig:largeamp}.   

The singularities of (\ref{eq:ulargeampreal}) are simple poles and therefore (\ref{eq:ulargeampreal}) ceases to be valid close to the singularities. It is shown in \ref{sec:applargeamp} that in a neighbourhood of the closest singularities,
\begin{equation}
   u(z,t) \sim \frac{\phi(\zeta)}{t (1 - \alpha^2t^2)}, \qquad \alpha \to \infty,  \label{eq:uzlargeamp}
\end{equation}
where $\phi$ is the ODE solution defined by (\ref{eq:V0m})--(\ref{eq:V0mbc})
and the formula
\begin{equation}
    \hspace{-2.5cm} z = i\left[ \cosh^{-1}(1/(\alpha t)) + t \sqrt{1-\alpha^2t^2}\left( 2\log(\alpha) - \frac{1 - 2\alpha^2t^2}{1 - \alpha^2t^2} - 2\log (\alpha t) - 2\log(1 - \alpha^2 t^2) - \zeta  \right)  \right], \label{eq:largeampzzeta}
\end{equation}
defines $\zeta$.
Therefore
\begin{equation}
\hspace{-2.5cm}\sigma(t) \sim \cosh^{-1}(1/(\alpha t)) + t \sqrt{1-\alpha^2t^2}\left( 2\log(\alpha) - \frac{1 - 2\alpha^2t^2}{1 - \alpha^2t^2} - 2\log (\alpha t) - 2\log(1 - \alpha^2 t^2) - \zeta_*  \right),\label{eq:sigmalargeamp}
\end{equation}
where $\zeta_*$ is the location of the first singularity of $\phi$ on the real axis, which, as stated in section~\ref{sec:smalltime}, is $\zeta_* \approx 0.05695$.

\subsection{Approach to blow up} As the singularities approach the real axis, i.e.~in the double limit $x \to 0$, $t \to 1/\alpha$, it follows from (\ref{eq:ulargeampreal}) that
\begin{equation*}
u \sim \frac{\alpha}{1- \alpha t+ x^2/2},
\end{equation*}
holds. 
This might suggest that as blow up is approached, the solution takes the self-similar form,
\begin{equation}
u = \frac{1}{t_c - t}f\left( \frac{x}{\sqrt{t_c - t}}  \right),  \label{eq:largeampblowup}
\end{equation}
this representing a classical similarity reduction to the NLH. As is well known, there were early conjectures in related problems that blow-up solutions indeed take such a self-similar form (see e.g., \cite{kapila1980,kassoy1980}), although in \cite{hocking} the appropriate logarithmic corrections had already been established for a specific form of nonlinearity; see also \cite{dold1991}. Indeed, it is known that no suitable solution to the resulting ODE in fact exists, see \cite{bebernes1987,eberly1988}.

In \ref{sect:appsectblowup} we apply appropriate modifications to existing approaches (pioneered in~\cite{hocking} for a closely related PDE) to analyse the approach to blow up.  The analysis is valid for all the types of initial data mentioned in this paper (namely, (\ref{eq:ic})--(\ref{eq:icflat})).
 

\section{Small-amplitude initial conditions}~\label{sec:smallamp}  
\begin{figure}[h!]
	\centering
	\includegraphics[width = 0.495\textwidth]{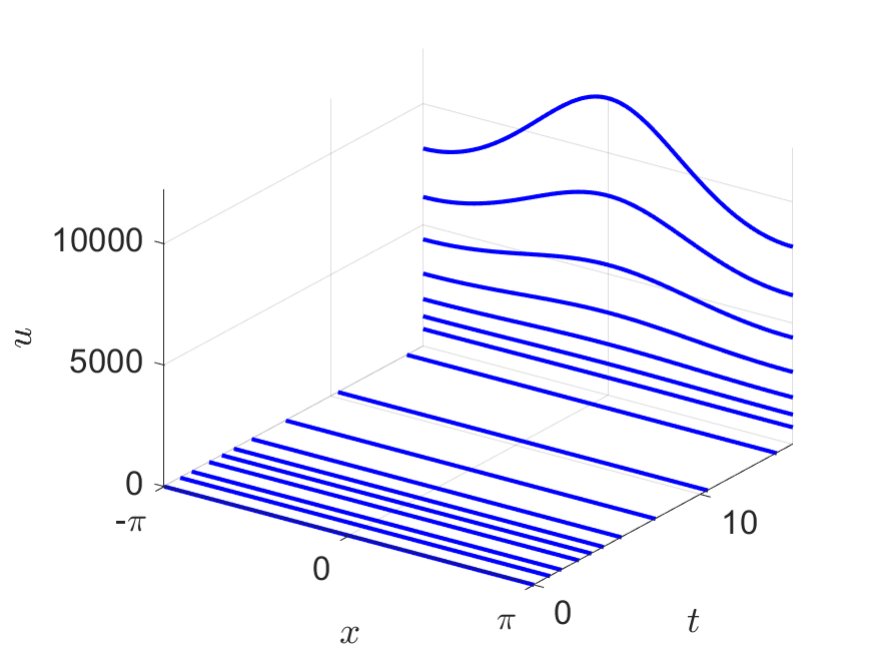}
 \caption{NLH solution for the initial data $u(x,0) = \epsilon\cos x$, with $\epsilon = \alpha = 0.5$.  This is the same solution displayed in the bottom-right frame of Figure~\ref{fig:blowup} but shown here on a different vertical scale.
 }\label{fig:small_amp_point_bu}
\end{figure}

In \ref{sec:appsmallamp}, the NLH solution is analysed for the initial data (\ref{eq:ic}) with $0 < \alpha = \epsilon \ll 1$ and it is found that five different time scales are relevant for the asymptotic analysis. Here we summarise the results on the different time scales.  The bottom right frames of Figures~\ref{fig:blowup} and~\ref{fig:sing} show, respectively, the qualitative behaviour of the solution profile and the singularity dynamics.   Figure~\ref{fig:small_amp_point_bu} illustrates that the  $\epsilon=0.5$ solution in Figure~\ref{fig:blowup}  appears to blow up uniformly  but in fact it blows up at the point $x = 0$.

\subsection{$t = \mathcal{O}(1)$}\label{sec:smallamptscale1main} On the first timescale, the solution on the real axis is given by
\begin{eqnarray}
  &\hspace{-2cm}u(x,t) \sim \epsilon \mathrm{e}^{-t}\cos x + \frac{\epsilon^2}{4}\left[ 1-\,{{\rm e}^{-2\,t}}+\, \left( {{\rm e}^{-2\,t}}-{{\rm e}^{-4
\,t}} \right) \cos  2\,x  \right] \nonumber \\
 & \hspace{-1cm} + \frac{\epsilon^3}{48}\left[ \left( 24 t+6\,{{\rm e}^{-2\,t}}+3\,{{\rm e}^{-4\,t}}-9
 \right) {{\rm e}^{-t}}\cos  x  + \left( 2 -3\,{{\rm e}^{-2\,t}}+{{\rm e}^{-6\,t}}\right) {{\rm e}^{-3\,t}}\cos
  3\,x   \right].
 \label{eq:utimescale1real}
\end{eqnarray}
This approximation, as well as other asymptotic estimates in the upcoming sections~\ref{sec:timescale2}--\ref{sec:timescale4}, will be compared to  numerical results in section~\ref{sec:smallampnumcomp}.

The closest singularities are at $z= \pm i \sigma(t)$, where
\begin{equation}
    \sigma(t) \sim  2t - \log(\sinh t) - \log(\epsilon/2) + \zeta_*(t), \qquad \epsilon \to 0,  \label{eq:sigmatimescale1m}
\end{equation}
and $\zeta_*(t)$ is the location of a singularity of a nonlinear backward diffusion PDE (namely, (\ref{eq:backdiff1})--(\ref{eq:backdiffbcs})) whose solution is not known explicitly. However, the limiting behaviour of $\zeta_*(t)$ is found to be (i) $\zeta_*(t) \to 0$ as $t \to 0$ (hence (\ref{eq:sigmatimescale1m}) is consistent with (\ref{eq:sigmasmalltime}) as $t \to 0$) and (ii) as $t$ becomes large according to $t = \mathcal{O}(1/\epsilon)$, $\zeta_* \to \widetilde{\zeta}_*-\log 2$.\footnote{The $\log 2$ shift arises because~(\ref{eq:backdiff1}) implies $U \sim 2 e^{\zeta}$
as $\zeta \to -\infty$, $t \to \infty$.} 
Here $\widetilde{\zeta}_*$ is the location of the first singularity of the following nonlinear ODE problem:
\begin{equation}
    \frac{d^2 \phi}{d\zeta^2} - \frac{d\phi}{d\zeta} = \phi^2, \qquad \phi \sim e^{\zeta}, \qquad \zeta \to -\infty.   \label{eq:nonlinODEprob2}
\end{equation}
The limiting behaviour (ii) follows because for $t = \mathcal{O}(1/\epsilon)$, $\epsilon \to 0$ the NLH solution in the neighbourhood of the singularity, i.e. for $\zeta = \mathcal{O}(1)$, is
\begin{equation}
u(z,t) \sim \phi(\zeta), \qquad     z = i\left( 2t - \log(\sinh t) + \log(1/\epsilon) + \zeta   \right),  \label{eq:NLHcomplextimescale1}
\end{equation}
where $\phi(\zeta)$ is the solution to (\ref{eq:nonlinODEprob2}).

\begin{figure}
	\centering
	\includegraphics[width = 0.495\textwidth]{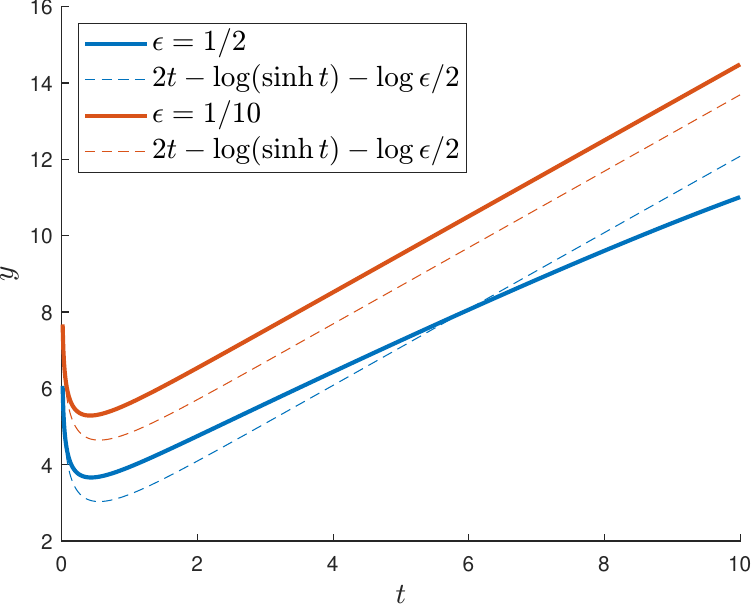}
\includegraphics[width = 0.495\textwidth]{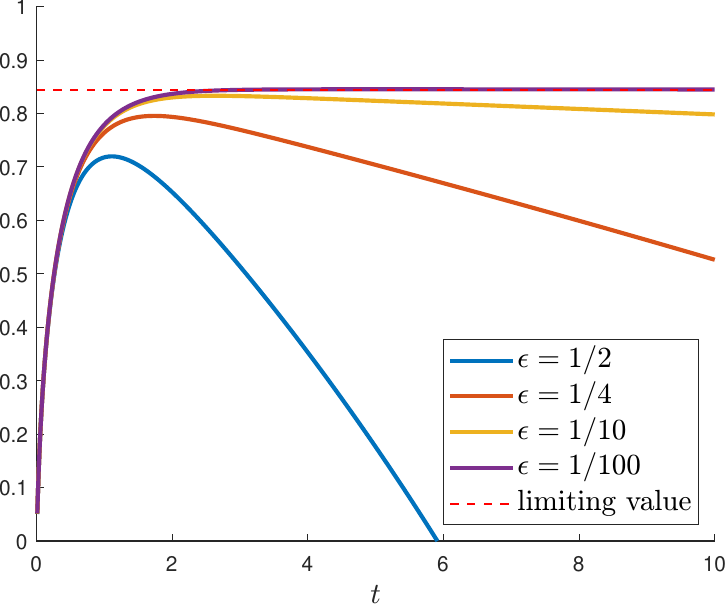}
	\caption{The left frame shows the numerical (solid lines) and asymptotic (dotted lines) approximations (with $\zeta_*(t)=0$ in (\ref{eq:sigmatimescale1m})) for the singularity locations $\pm i\sigma(t)$. The difference between these approximations, shown in the right frame, approaches the predicted limiting value, $\widetilde{\zeta}_*-\log 2\approx 0.84452$, (indicated by the dotted line) for large $t$ as $\epsilon \to 0$.  
 } 
	\label{fig:1stturnaround}
\end{figure}

\begin{figure}[h]
    \hspace*{-0cm}\includegraphics[width = 0.9\textwidth]{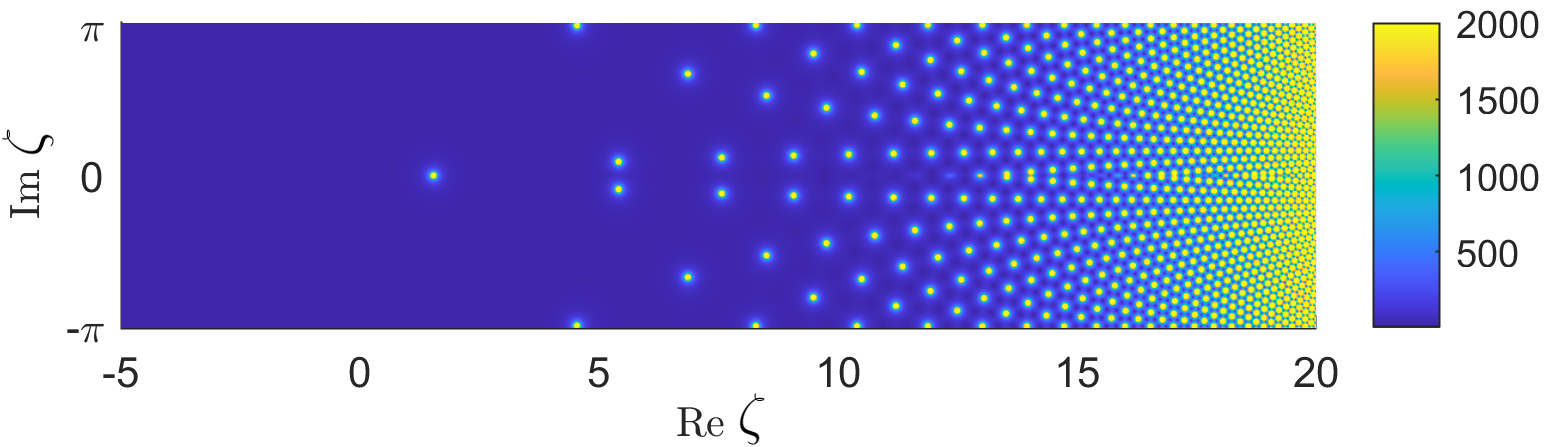}\\
\centering    \includegraphics[width = 0.45\textwidth]{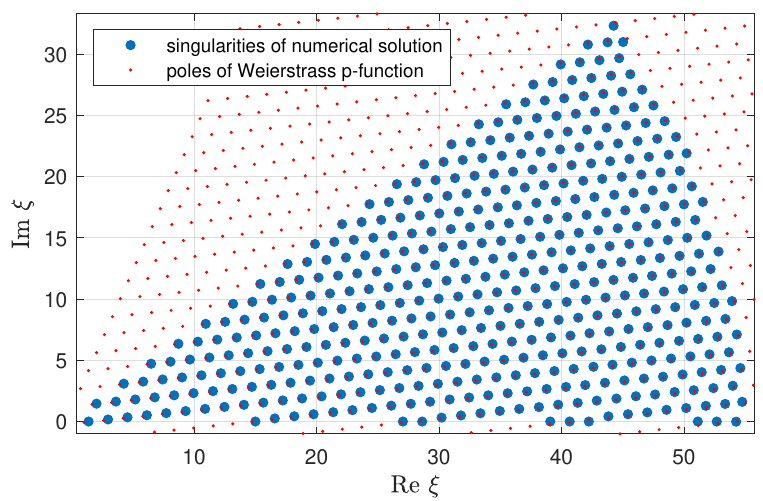}
    \caption{Top: Modulus of the solution to (\ref{eq:nonlinODEprob2}) with initial conditions given by (\ref{eq:nonlinODE2ics}), as computed by the procedure described in \ref{sec:appbcscompute}.
    Bottom: Comparison of the pole locations (\ref{eq:weierpoles}) of the Weierstrass function (with $\alpha = r\exp(i\theta)$, $r \approx 0.28910$, $\theta \approx 0.1066$ and $\zeta_0 \approx -0.603 - 0.2574i$)  and the singularity locations on the upper half-plane in the top frame of the figure after
    mapping to the $\xi$-plane ($\xi = e^{\zeta/5}$). 
    }\label{fig:ODE_exp_mod}
\end{figure}

The left frame of Figure~\ref{fig:1stturnaround} compares the numerically determined singularity location (as described in section~\ref{sec:fourier}) and the approximation (\ref{eq:sigmatimescale1m}), but with $\zeta_*(t)$ neglected. The right frame of the figure shows the difference between these two quantities, which gives an estimate for $\zeta_*(t)$. As predicted by the asymptotic analysis, $\zeta_*(t)$ increases from zero and for large $t$ and $\epsilon \to 0$, $\zeta_* \to \widetilde{\zeta}_*-\log 2$, where $\widetilde{\zeta}_* \approx 1.53767$, see the top frame of  Figure~\ref{fig:ODE_exp_mod}.
The latter figure shows that the far-field condition in (\ref{eq:nonlinODEprob2}) leads to multiple singularities on the principal Riemann sheet of the ODE solution, in contrast to the solution  in Figure~\ref{fig:ODE_complex_phase_modulus} 
corresponding to the  condition (\ref{eq:V0mbc}). 
 This suggests that the NLH solution has multiple singularities that in the small-time limit live off its principal Riemann sheet but, for small amplitude and large time, move onto its principal Riemann sheet, a possibility we alluded to above.

 \ref{sec:appivp2} shows that the far-field solution of (\ref{eq:nonlinODEprob2}) can also be expressed to leading order in terms of the (equianharmonic) Weierstrass function in the variable $\xi = e^{\zeta/5}$. Hence the far-field singularities of the solution in the top frame of Figure~\ref{fig:ODE_exp_mod} lie approximately on the same lattice as the Weierstrass function's singularities in $\xi$-plane, as shown in the bottom frame of the figure. 

\subsection{$t = \mathcal{O}(\epsilon^{-2})$}\label{sec:timescale2} On the second timescale, the solution on the real axis is
\begin{equation}
    u \sim \frac{1}{t_c - t} + \frac{16\, {\rm e}^{-t}}{\epsilon^3(t_c- t)^2}\cos x  +  \left( \frac{128\, e^{-4t}}{\epsilon^6(t_c - t)^2}\int_{-\infty}^{t}\frac{e^{2s}}{(t_c - s)^2}ds\right)\cos 2x
    \label{eq:utimescale2real}
\end{equation}
where the leading order estimate of the blow-up time is $t_c \sim 4/\epsilon^2$ (see \ref{sec:appsmallamp}). 
 For the solution in Figure~\ref{fig:small_amp_point_bu} with $\epsilon = 1/2$, this gives the estimate $t_c \approx 16$, while the numerically computed blow-up time is $t_c \approx 15.53$.

\begin{figure}
	\centering
	\includegraphics[width = 0.5\textwidth]{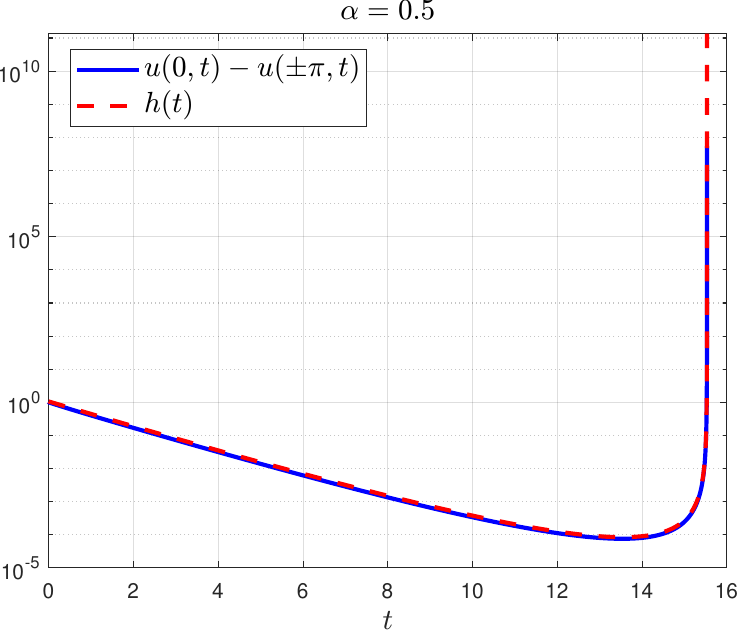}
	\caption{The peak-to-trough height of the solution in Figure~\ref{fig:small_amp_point_bu} (solid curve) compared to the approximation (\ref{eq:uheightapprox}) (dashed curve).  
 As expected, the approximation (\ref{eq:uheightapprox}) breaks down close to the blow-up time.  } \label{fig:small_amp_sol_height}
\end{figure}

Solutions with initial data (\ref{eq:ic}) have a maximum at $x = 0$ (see Figure~\ref{fig:blowup}) and, as shown in Appendix~B of~\cite{fasondini2023blow}, the peak-to-trough height of the maximum satisfies
\begin{equation*}
    u(0,t) - u(\pm \pi, t) \approx 4c_1(t) := h(t)
\end{equation*}
for even solutions, provided $t$ is not close to the blow-up time.  Here $c_1(t)$ is a Fourier coefficient of the solution; see (\ref{eq:uFourier}). From (\ref{eq:utimescale2real}) it follows that 
\begin{equation}
h(t) = \frac{32\, e^{-t}}{\epsilon^3(t_c - t)^2}, \quad \Rightarrow \quad \min h(t) = h(t_c-2) = \frac{8\,e^2}{\epsilon^3}e^{-t_c} \sim \frac{8\,e^2}{\epsilon^3}e^{-4/{\epsilon^2}},   \label{eq:uheightapprox}
\end{equation} 
which implies
the minimum peak-to-trough height of the solution decreases exponentially with $\epsilon$.
With $\epsilon = 0.3$ it follows that $\min h(t) \approx 10^{-16}$, which means that  the solution is completely flat in double-precision arithmetic and a uniform blow up is computed. This is incorrect since, as we shall see, point blow up eventually occurs for any $\epsilon > 0$.  Consequently, higher precision or rescaling methods would be required to compute the solution for $\epsilon < 0.3$.  Figure~\ref{fig:small_amp_sol_height} compares the asymptotic estimate (\ref{eq:uheightapprox}) to the numerically computed peak-to-trough height and confirms the validity of the estimate away from the blow-up time.

In the neighbourhood of the singularity, i.e.~for $\zeta = \mathcal{O}(1)$, the NLH solution is given by
\begin{equation*}
  u(z,t) \sim \phi(\zeta),\qquad  z = i\Big( t + 2\log\left( t_c - t \right) + 3\log\left(\epsilon/2  \right)   + \zeta\Big),
\end{equation*}
where $\phi$ is again defined as the solution to (\ref{eq:nonlinODEprob2}).
Hence an estimate of the singularity location is 
\begin{equation}
 \sigma =   t + 2\log\left( t_c - t \right) + 3\log\left(\epsilon/2  \right)   + \zeta_*,  \label{eq:sigmatimescale2m}
\end{equation}
where $\zeta_*$ is the first singularity of 
(\ref{eq:nonlinODEprob2})
on the real axis, which was found to be $\zeta_*\approx 1.53767$ in section~\ref{sec:smallamptscale1main}.  This estimate for $\sigma$, as well as the estimates in the upcoming sections, will be compared to numerical results in section~\ref{sec:smallampnumcomp}.

\subsection{$t_c - t = \mathcal{O}(1)$}  On this time scale, the solution on the real axis is given again by (\ref{eq:utimescale2real}).  The singularity location evolves on the imaginary axis according to (cf. (\ref{eq:sigmatimescale2m}))
\begin{equation}
 \sigma(t) \sim   t + 2\log\left( t_c - t \right) + 3\log\left(\epsilon/2  \right)   + \zeta_*(t),  \label{eq:sigmatimescale3m}
\end{equation}
where $\zeta_*(t)$ is the position of the first singularity of another nonlinear backward diffusion PDE initial-value problem ((\ref{eq:PDEsmallt3})--(\ref{eq:t3UYlarge})) whose   solution is not known explicitly.  However, the limiting behaviour of $\zeta_*(t)$ is shown in \ref{sec:apptimescale3} to be as follows: for $1 \ll t_c-t \ll 1/\epsilon^2$ with $\epsilon \to 0$,
$\zeta_*(t) \to \zeta_*$, where $\zeta_*$ is the first singularity of (\ref{eq:nonlinODEprob2}) on the real axis, i.e.~$\zeta_*\approx 1.53767$ and therefore (\ref{eq:sigmatimescale3m}) tends to (\ref{eq:sigmatimescale2m}).  For $t \to t_c$, to leading order $\zeta_*(t) \sim  -\log(t_c-t)$ and  thus
\begin{equation}
    \sigma(t) \sim   t + \log\left( t_c - t \right) + 3\log\left(\epsilon/2  \right), \qquad t \to t_c.  \label{eq:sigmatimescale3m2}
\end{equation}

\subsection{Fourth time scale}\label{sec:timescale4} 
The fourth time scale, which is valid exponentially close to the blow-up time, is defined via
\begin{equation*}
    t_c - t = \frac{s}{\epsilon^3 e^{4/\epsilon^2}}
\end{equation*}
with $s = {\mathcal O}(1)$, and to leading order
\begin{equation}
    u(x,t) \sim \frac{\epsilon^3 e^{4/\epsilon^2}}{s - 16 \cos x}.   \label{eq:urealtimescale4}
\end{equation}
Therefore blow up occurs at $s \sim 16$ and the modification to the algebraic expansion for $t_c(\epsilon)$ resulting from the previous timescales is the exponentially small, and hence in practice irrelevant quantity
\begin{equation*}
    -16 e^{-4/\epsilon^2}/\epsilon^3.
\end{equation*}
An approximation for the locations of the closest singularities (which will be used in section~\ref{sec:smallampnumcomp}) is $z = \pm i \sigma(t)$ with
\begin{equation}
\sigma \sim \cosh^{-1}(s/16) =  \log(s + (s^2 - 256)^{1/2}) - 4\log 2.  \label{eq:sigmatimescale4m}
\end{equation}

\subsection{Fifth time scale}\label{sec:timescale5}

As blow up is approached, $s \to 16$, $x \to 0$ apply in (\ref{eq:urealtimescale4}) and therefore
\begin{equation*}
    u \sim \frac{\epsilon^3 e^{4/\epsilon^2}}{s - 16 + 8x^2}.
\end{equation*}
This suggests the self-similar form
\begin{equation*}
    u = \frac{\epsilon^3 e^{4/\epsilon^2}}{s - 16}f\left( \frac{x}{(s - 16)^{1/2}} \right),
\end{equation*}
cf.~(\ref{eq:largeampblowup}).
Again, the appropriate ODE solution fails to exist, which necessitates the introduction of a final (doubly exponentially short) time variable to capture the logarithmic corrections.   The approach to blow up is analysed in more detail in \ref{sect:appsectblowup}.

\subsection{Comparisons to numerical results}\label{sec:smallampnumcomp}  
\begin{figure}
	\centering
	\includegraphics[width = 0.425\textwidth]{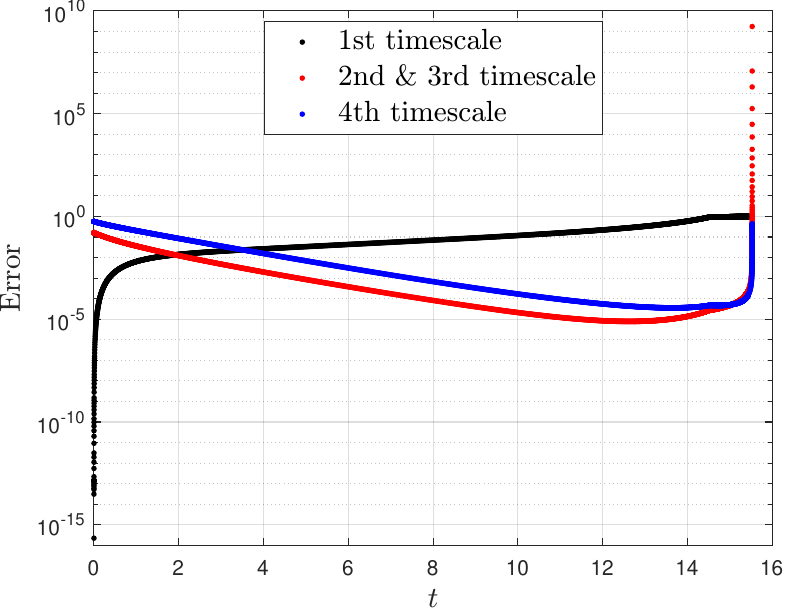}
\includegraphics[width = 0.45\textwidth]{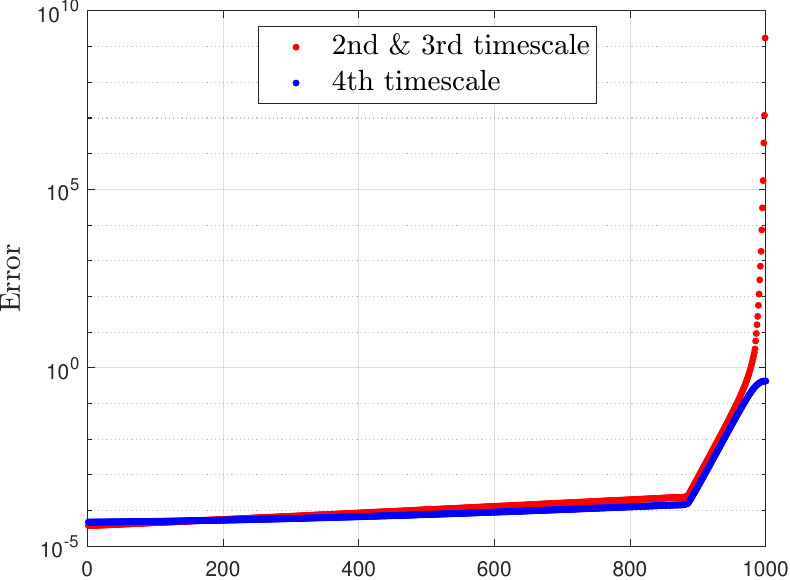}
	\caption{The error of the approximations (\ref{eq:utimescale1real}) (first timescale), (\ref{eq:utimescale2real}) (second and third timescales) and (\ref{eq:urealtimescale4}) (fourth timescale) compared to the numerical solution of the NLH with $u(x,0) = \epsilon \cos x$, $\epsilon = 0.5$.  The error is calculated at every time step as the minimum of the absolute and relative errors on $x \in [-\pi, \pi]$.  In the right frame, we `zoom in' on the error close to the blow-up time by showing the error for the last 1000 steps of the numerical time integrator, which corresponds to the time interval $t \in [t_1, t_c]$, with  $t_1 = 14.87\ldots$ and  $t_c = 15.53\ldots$. } \label{fig:asymptreal}
\end{figure}

Having derived a large number of estimates for the small amplitude case, we now compare some of them to numerical approximations. 
Figure~\ref{fig:asymptreal} shows the accuracy of the asymptotic approximations to the solution on the real axis given in (\ref{eq:utimescale1real}), (\ref{eq:utimescale2real}) and (\ref{eq:urealtimescale4}) for the 
 case $\alpha = \epsilon = 0.5$.
 We note that the need for a comparatively large value of $\epsilon$ reflects the asymptotic structure whereby the small quantity $\exp(-4/\epsilon^2)$ is prominent: the extent to which $u$ becomes near uniform is a striking feature of the analysis.
 In (\ref{eq:utimescale2real}), we choose $t_c$ to be the numerical blow-up time and in (\ref{eq:urealtimescale4}) we choose $t_c$ such that the blow-up time of (\ref{eq:urealtimescale4}), namely $s = 16$ (i.e., at $t = t_c - 16\epsilon^{-3}e^{-4/\epsilon^2}$), coincides with the numerical blow-up time.

\begin{figure}
   \centering
   \includegraphics[width = 0.5\textwidth]{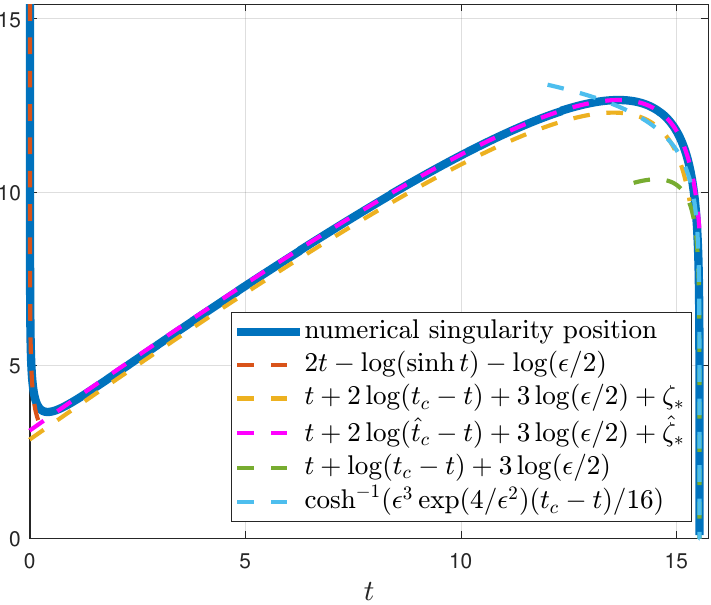}
    \caption{The position of the closest singularity on the positive real axis of the small-amplitude NLH solution with  $u(x,0) = \epsilon\cos x$ with $\epsilon = 0.5$ compared to the asymptotic approximations.
    }  
    \label{fig:smallamp_turnaround2}
\end{figure}

Figure~\ref{fig:smallamp_turnaround2} compares the asymptotic estimates of the closest singularities at $z = \pm i \sigma(t)$ (see (\ref{eq:sigmatimescale1m}), (\ref{eq:sigmatimescale2m}), (\ref{eq:sigmatimescale3m}), (\ref{eq:sigmatimescale3m2}) and (\ref{eq:sigmatimescale4m})) to the numerically computed singularity position.  For all these asymptotic approximations, we let $t_c$ be the numerically computed blow-up time.   For the estimate (\ref{eq:sigmatimescale2m}), we use the value of $\zeta_* = 1.53767$, while for (\ref{eq:sigmatimescale3m}), since the function $\zeta_*(t)$ is not known explicitly, we replace $t_c$ and $\zeta_*(t)$ with constants\footnote{In Figure~\ref{fig:smallamp_turnaround2}, these constants are $\hat{t}_c = 15.65$ and $\hat{\zeta}_* = 1.78$.} $\hat{t}_c$ and $\hat{\zeta}_*$ such that the local maximum of the resulting estimate matches that of the numerical singularity position.

\section{Blow-up limit}\label{sec:blowupmain}  
\begin{figure}[h!]
   \centering
\mbox{\includegraphics[width = 0.495\textwidth]{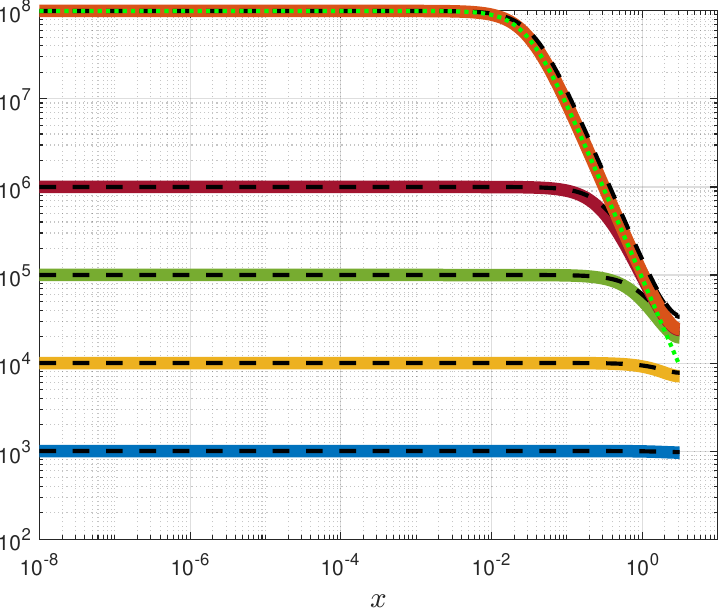}
   \includegraphics[width = 0.495\textwidth]{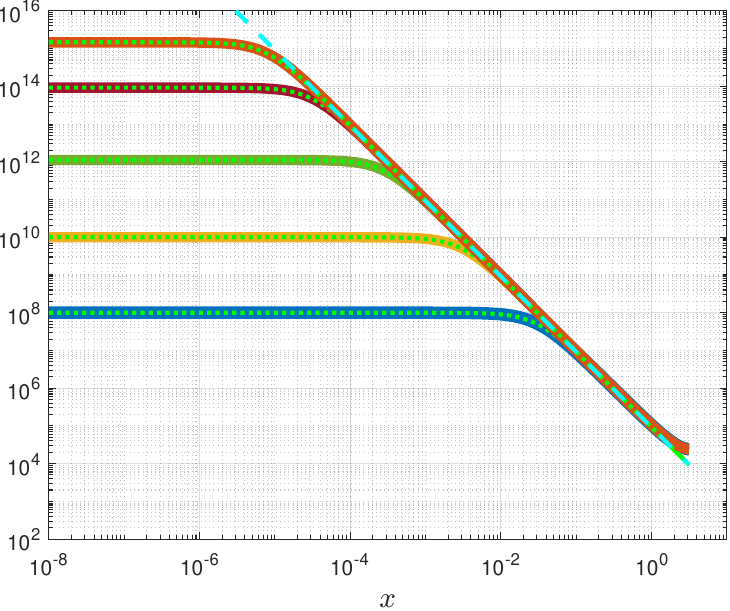}}
    \caption{Left: The five solid curves are NLH solutions with initial data $u(x,0) = \alpha \cos x$, $\alpha = 0.5$ corresponding to (from bottom to top) $t_c - t = 10^{-3}, 10^{-4}, 10^{-5}, 10^{-6}, 10^{-8}$.  The dashed curves show the estimate (\ref{eq:urealtimescale4}) and the single dotted curve is the estimate (\ref{eq:ublowupzeta}). 
    Note the logarithmic scale on the $x$-axis.  Right: The solid curves show the NLH solution (with the same initial data as in the left frame) with $t_c - t = 10^{-8}, 10^{-10}, 10^{-12}, 10^{-14}, 10^{-15}$.
    The dotted curves show the estimate (\ref{eq:urealtimescale4}) and the single dashed curve is the blow-up profile (\ref{eq:blowupprofile}).  For the asymptotic estimates, we use the numerically determined value of $t_c = 15.530458826185942$ and the values $C = 92000$ and $\beta_1 = -3/32$, which were chosen to fit the numerical solution.   }  \label{fig:blowupprofile}
\end{figure}

Figure~\ref{fig:blowupprofile} shows 
 a small-amplitude NLH solution at times approaching the blow-up time.  Since the solution is even, it is shown only for $x > 0$ with $x \in [10^{-8}, \pi]$. The figure illustrates the well-known fact that, as the blow-up time is approached, the solution is flat for $\eta = \mathcal{O}(1)$, where $\eta = x/\sqrt{t_c - t}$; see \ref{sec:appbuscale1}.  That is, for a fixed $t$ with $0< t_c -t \ll 1$, the solution is flat with $u \sim (t_c - t)^{-1}$ for $x$ sufficiently small.  More precisely, as shown in \ref{sec:appbuscale2}, for $x = \mathcal{O}((t_c - t)\log(t_c-t)^{1/2})$,
 \begin{equation}
 u(x,t) \sim \left[ t_c - t + \frac{x^2}{C - 8 \log(t_c - t)} \right]^{-1}, \qquad t \to t_c, \label{eq:ublowupzeta} 
 \end{equation}
 where $C$ is a constant that depends on the initial data.
 This approximation is shown as dotted curves in Figure~\ref{fig:blowupprofile} and matches the numerical solution well for sufficiently small $x$.  


 As is clear from Figure~\ref{fig:blowupprofile},  the solution is asymptotically flat on an interval whose width shrinks to zero as $t \to t_c^{-}$ 
 and, as shown in \ref{sec:appblowup3rdscale}, the solution acquires the  blow-up profile,
 \begin{equation}
    u(x,t_c) \sim \frac{8}{x^2}\left( 2\log(1/x) + \log(\log(1/x)) + 4\log 2 + C/8 + 8\beta_1 \right), \quad x \to 0^{+},  \label{eq:blowupprofile}
\end{equation}
where the constant $\beta_1$ is defined in~\ref{sect:appsectblowup}.
This profile is shown as a dashed curve in the right frame of Figure~\ref{fig:blowupprofile}, which matches the non-flat part of the numerical solution well for sufficiently small $x$.  We note that, before the blow-up time, the closest singularities are second-order poles to leading order with a logarithmic term at fourth order (see (\ref{eq:expansion})); however, in the blow-up limit, the leading-order behaviour for $x \to 0$, namely $u \sim 16\log(1/\vert x \vert)/x^2$, acquires a logarithmic contribution.



\begin{figure}[h!]
   \centering
\mbox{\includegraphics[width = 0.495\textwidth]{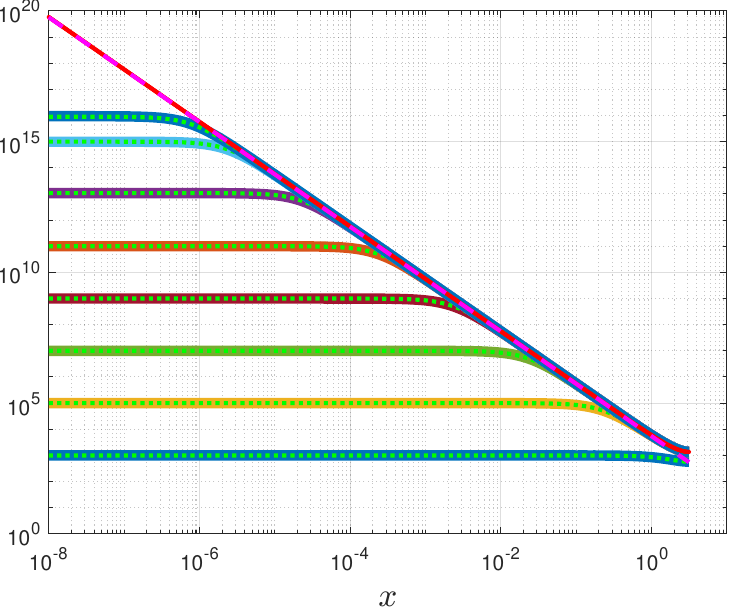}\includegraphics[width = 0.495\textwidth]{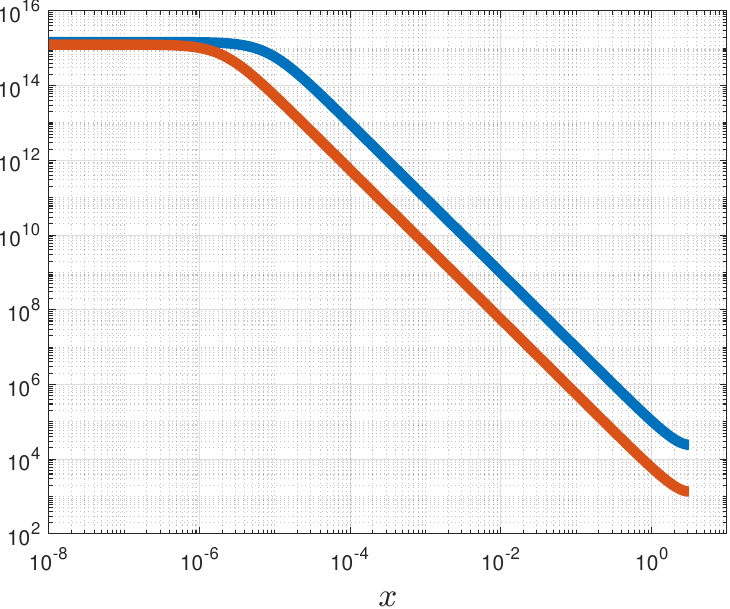}}
    \caption{Left: The solid curves show NLH solutions subject to (\ref{eq:icflat}) with $\alpha = 1$, $\epsilon = 0.001$ and $t_c-t = 10^{-3}, 10^{-5}, \ldots, 10^{-13}, 10^{-15}, 10^{-16}$; the dotted curves are defined by (\ref{eq:vtimescale2}) and the red and pink dashed curves are given by, respectively, (\ref{eq:vblowupprofsimp}) and (\ref{eq:BUexpsmallx}).  Right: the solution on the left with $t_c-t=10^{-15}$ (red curve) and  the small-amplitude solution in the right frame of Figure~\ref{fig:blowupprofile} with $t_c-t=10^{-15}$ (blue curve).}  \label{fig:blowupprofile2}
\end{figure}

Unlike the asymptotic estimates (\ref{eq:ublowupzeta}) and (\ref{eq:blowupprofile}), the asymptotics of the blow-up profile in~\cite{fasondini2023blow} are valid on the entire interval $[-\pi, \pi]$ (they are $2\pi$-periodic) and the constants are expressed explicitly in terms of the initial data considered in that paper, namely (\ref{eq:icflat}). For comparison purposes, we restate the analogue of (\ref{eq:ublowupzeta}) from~\cite{fasondini2023blow}: as $t \to t_c$
\begin{eqnarray}
& \hspace{-1.5cm} u(x,t) \sim \biggl[ t_c-t + 2\epsilon\, e^{-\alpha}\sin^2(x/2)  +2\epsilon^2 \log \epsilon\, e^{-2\alpha}\sin^2 x +  \epsilon (t - t_c)e^{-\alpha}\cos x \nonumber \\
&  \qquad  + 2\epsilon^2\sin^2 x\bigg(e^{-2\alpha} \log\left( \frac{t_c-t}{\epsilon} + 2e^{-\alpha}\sin^2(x/2)   \right) + C_1 + C_2\bigg)\biggr]^{-1}, 
\label{eq:vtimescale2}
\end{eqnarray}
where $\alpha$ and $\epsilon$ are parameters in the initial data (\ref{eq:icflat}) and
\begin{equation*}
C_1 = e^{-2\alpha}\log\alpha, \qquad 
C_2 = e^{-4\alpha}\int_0^{\alpha}\frac{e^{2t}-e^{2\alpha}}{\alpha - t}dt. 
\end{equation*}
Setting $t = t_c$ in (\ref{eq:vtimescale2}), we obtain
\begin{equation}
\hspace{-2cm}u(x,t_c) \sim  \left[2 \epsilon\, e^{-\alpha} \sin^2 (x/2) +   
    2 \epsilon^2 \sin^2 x  \Big(e^{-2\alpha}   
   \log\left(2 \epsilon\, e^{-\alpha}\sin^2(x/2)\right)    + C_1 + C_2 \Big)\right]^{-1}.
\label{eq:vblowupprofsimp}
\end{equation}
As shown in~\cite{fasondini2023blow}, at the blow-up time and for $x$ exponentially small with respect to $\epsilon$, the following analogue of (\ref{eq:blowupprofile}) holds,
\begin{equation}
u(x,t_c) \sim 
\frac{8}{x^2}\left(2\log(1/x) + (4\epsilon\, e^{-\alpha})^{-1}   \right),
\qquad x \to 0.  \label{eq:BUexpsmallx}
\end{equation}
The left frame of Figure~\ref{fig:blowupprofile2} shows the accuracy of the asymptotic approximations (\ref{eq:vtimescale2})--(\ref{eq:BUexpsmallx}) and the right frame compares the small-amplitude NLH solution with $t_c-t = 10^{-15}$ in Figure~\ref{fig:blowupprofile} with another NLH solution with $t_c-t = 10^{-15}$ but subject to the initial data (\ref{eq:icflat}).

\begin{figure}[h!]
    \centering
  \mbox{\includegraphics[width = 0.495\textwidth]{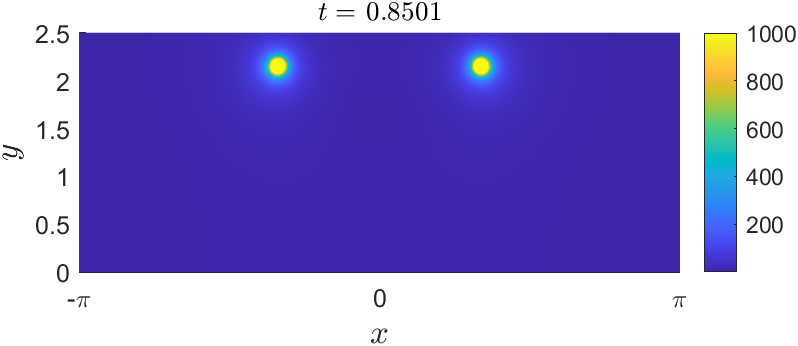}\includegraphics[width = 0.495\textwidth]{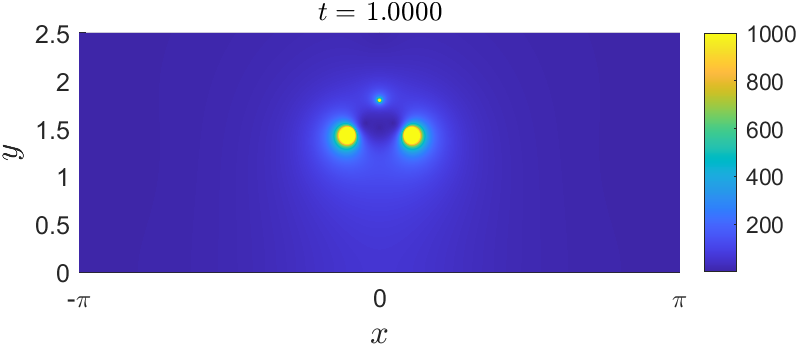}} 
  \mbox{\includegraphics[width = 0.495\textwidth]{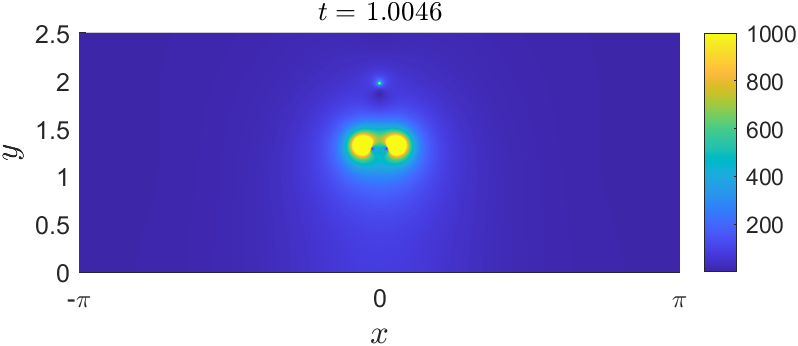}\includegraphics[width = 0.495\textwidth]{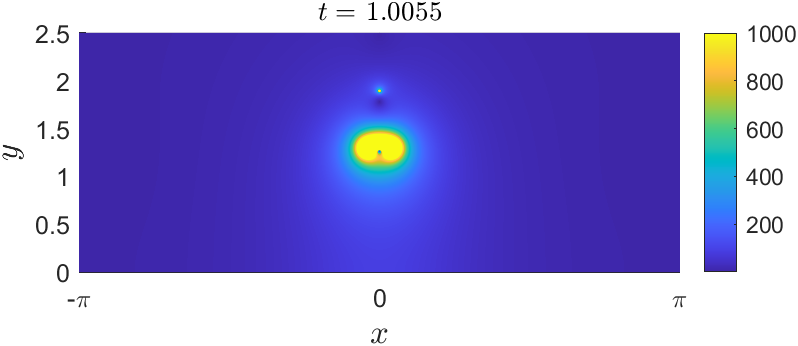}}
  \mbox{\includegraphics[width = 0.495\textwidth]{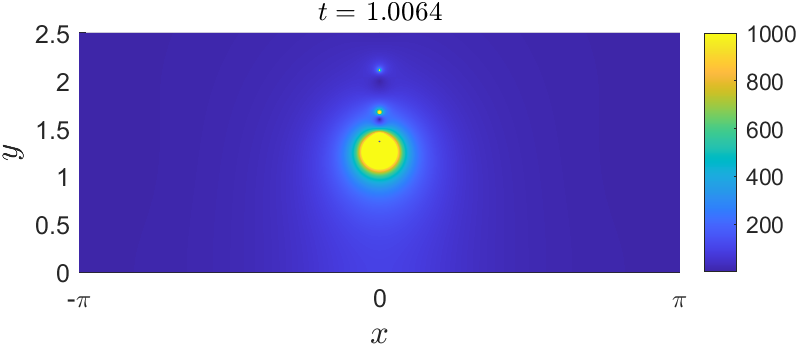}\includegraphics[width = 0.495\textwidth]{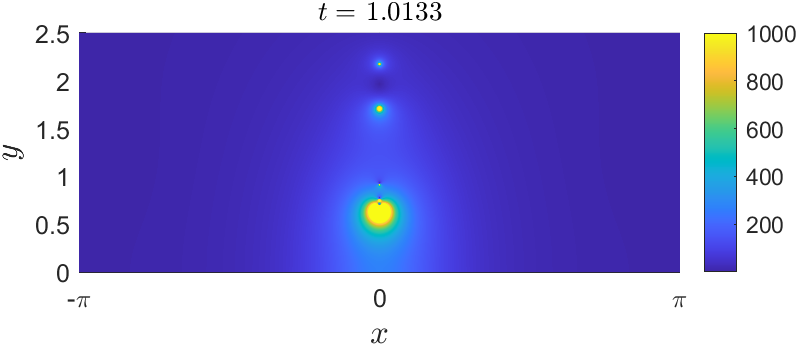}}
  \caption{Modulus plot of  the solution in the right frame of Figure~\ref{fig:2peaksreal} in the upper half-plane.  In the latter figure, the maxima (shown as red dots) coalesce at $t \approx 0.855$, whereas the singularities above collide $t \approx 1.006$.
  }\label{fig:singscollide}
   \end{figure}

\begin{figure}[h!]
    \centering
  \mbox{\includegraphics[width = 0.495\textwidth]{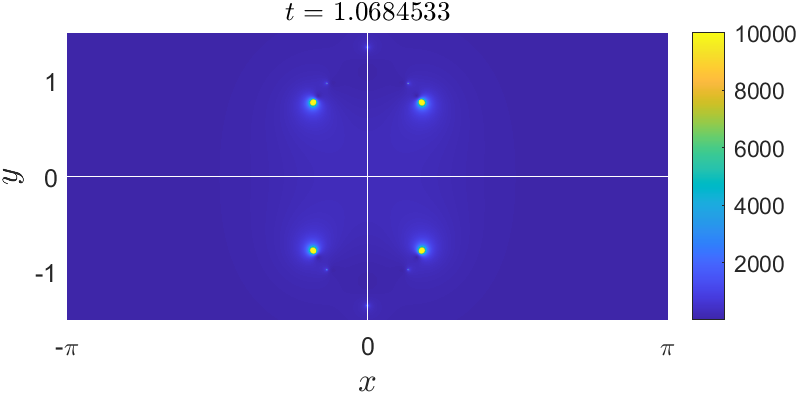}\includegraphics[width = 0.495\textwidth]{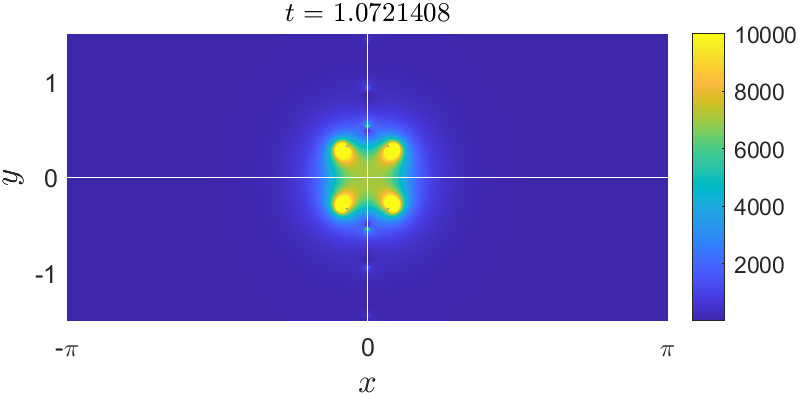}} 
  \mbox{\includegraphics[width = 0.495\textwidth]{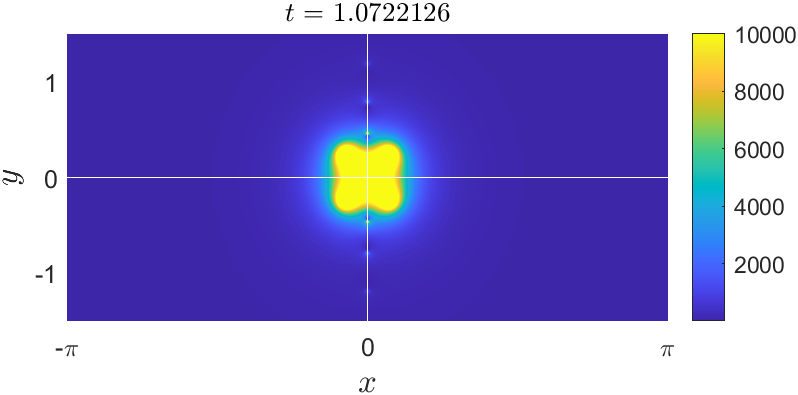}\includegraphics[width = 0.495\textwidth]{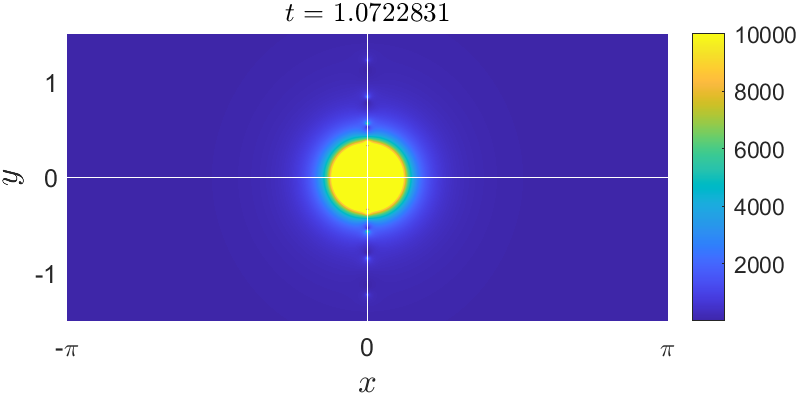}}
  \caption{Modulus plot of an NLH solution exhibiting non-generic blow up due to four singularities colliding at $x = 0$ at the blow-up time.  This solution has initial data
 $u(x,0) = \alpha \exp(\mu \cos(x + \delta x)-\mu) + \alpha \exp(\mu \cos(x + \delta x)+\mu)$ with $\alpha = 6$, $\mu = 50$ and $\delta = 0.4363\pi$.  Note that there is a factor of 10 difference between the colour maps (indicating the modulus of the solution) in this figure and in Figure~\ref{fig:singscollide}.   }\label{fig:4singscollide}
   \end{figure}

\begin{figure}[h!]
    \centering
  \includegraphics[width = 0.49\textwidth]{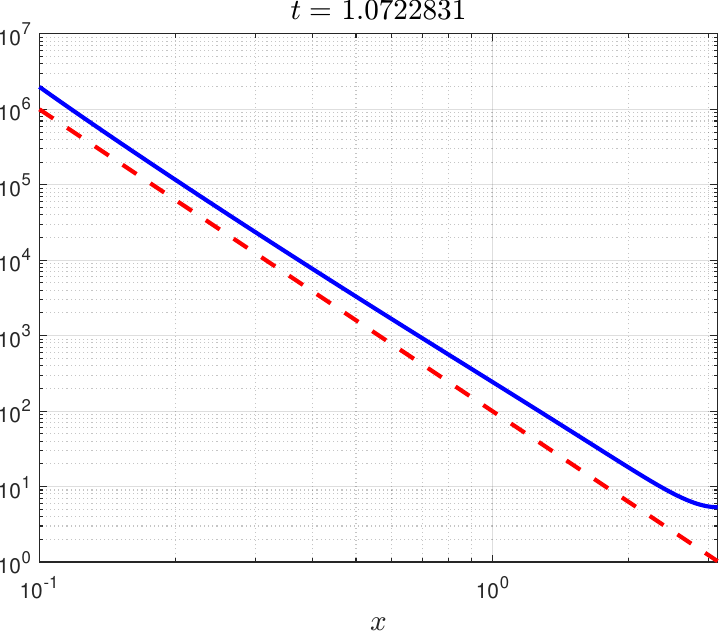}
  \caption{The (non-generic) blow-up profile of the solution in the bottom-right frame of Figure~\ref{fig:4singscollide} (solid line) and  a curve that grows as $\mathcal{O}(x^{-4})$, $x \to 0$ (dotted line). 
} \label{fig:4thorderblowup}
   \end{figure}

As discussed in the introduction, for an even initial condition with two local maxima, blow up of the type we have discussed here (generic blow-up) can occur simultaneously at two points (see the left frame of Figure~\ref{fig:2peaksreal}) or at a single point (right frame of Figure~\ref{fig:2peaksreal}). In the former case, two singularities collide on the real axis at both blow-up points and in the latter case (see Figure~\ref{fig:singscollide}), a pair of singularities coalesces in each of the upper and lower half-planes before the resulting singularities collide on the real axis at blow-up.  The borderline case, in which two singularities from the upper half-plane and two singularities from the lower half-plane collide on the real axis at the same point at the blow-up time (see Figure~\ref{fig:4singscollide}) is a type of non-generic blow up for which the leading-order behaviour  at blow up  is $u \sim C/x^4$, $x \to 0$, where $C$ depends on the initial data. (See \ref{sec:non-generic} and Figure~\ref{fig:4thorderblowup}.)  

\section{Conclusion}  For the NLH we have given asymptotic descriptions (along with numerical confirmation) of: its branch-point-type singularities, the solution in the neighbourhood of the closest singularities in the complex plane, the complex-plane dynamics of the closest singularities and the solution in the small-time, large-amplitude, small-amplitude and blow-up limits.


 


Numerous generalisations suggest themselves; perhaps the most immediate is the quasilinear power-law case
\begin{equation*}
    \frac{\partial u}{\partial t} = \frac{\partial}{\partial x}   \Big(u^m  \frac{\partial u}{\partial x}   \Big) + u^p,  
\end{equation*}
with $m \neq 0$, $p>1$, for which the blow-up behaviour is well-known to differ from that for the semilinear case $m=0$ with which we have been concerned here; indeed, the latter can be viewed as the borderline case between two distinct classes of behaviour, the logarithmic terms that are prevalent in the above being associated with this borderline status.

We related the NLH solution to nonlinear ODE solutions (which have interesting properties in their own right) in certain limits and found the singularity locations of the ODE solutions numerically (via the pole field solver) and asymptotically.
A possibility for future work would be to develop an analogue of the pole field solver for PDEs.  Just as the pole field solver can accurately compute multivalued ODE solutions by using the ODE and adaptive Pad\'e approximation to  continue analytically the solution onto multiple Riemann sheets,
so a pole field solver for PDEs might be able to compute, for a fixed time $t$,
the NLH solution on multiple Riemann sheets in the complex $x$ plane.  The numerical analytic continuation method for the NLH that we used in this paper (Pad\'e and quadratic Pad\'e approximation using the Fourier expansion of the solution) is accurate in a neighbourhood of the closest singularities of the NLH.  However, it rapidly loses accuracy as one moves further away from the real axis.  A pole field solver for PDEs would presumably maintain accuracy much further away from the real axis and also onto neighbouring Riemann sheets.


\section{Acknowledgements} The authors would like to thank the Isaac Newton Institute for Mathematical Sciences for support and hospitality during the programme \textit{Complex analysis: techniques, applications and computations} when work on this paper was undertaken. This programme was supported by EPSRC grant number EP/R014604/1. The work of the first author was also supported by the Leverhulme Trust Research Project Grant RPG-2019-144.  The second author gratefully acknowledges a Royal Society Leverhulme Trust Senior Fellowship. The third author
acknowledges a grant from the
H.B.~Thom Foundation of Stellenbosch University that enabled participation in the above mentioned programme.

\section*{References}
\bibliographystyle{elsarticle-num}

\bibliography{main_Nonlinearity}
\appendix
\section{Analysis of the small-time limit}\label{sec:appsmalltime} The Taylor expansion of the NLH in $t$ with initial condition (\ref{eq:ic}) yields
\be
u = \alpha \cos x + t \Big(-\alpha \cos x + \frac12 \alpha^2 (1+\cos 2x)\Big) + \mathcal{O}(t^2),  \label{eq:smalltexp1}
\ee
which holds sufficiently close to the real $x$ axis and for small values of $t$.
Setting $x = i y$ with $y \to \infty$, the terms
$\alpha \cos x \sim \frac12 \alpha e^y$ and $\frac12 t \alpha^2 \cos(2x) \sim 
\frac14 t \alpha^2 e^{2y}$ come into balance, leading to a non-uniformity in (\ref{eq:smalltexp1}) when $Y = \mathcal{O}(1)$, where 
\begin{equation}
y = - \log (\alpha t/2) + Y. \label{eq:smallts}
\end{equation}

Setting $u = U/t$, (\ref{eq:smalltexp1}) implies the matching condition
\begin{equation}
U \sim e^{Y} + e^{2Y} - t \, e^{Y}, \qquad \mathrm{as\: }  t\to 0, Y\to -\infty. \label{eq:Uexp1}
\end{equation}
In these variables, the NLH becomes
\begin{equation*}
{\frac {\partial U}{
		\partial Y}} -U- U^{2} +
t\left( {\frac {\partial U}{\partial t}}  + {\frac {
		\partial ^{2} U}{\partial {Y}^{2}}}   \right) 
  =0.
\end{equation*}
Letting $U \sim U_0(Y) + tU_1(Y)$, as $t \to 0$ with $Y = \mathcal{O}(1)$ we find that $U_0$ and $U_1$ must satisfy
\be
\frac{dU_0}{dY} =U_{{0}}+U_{{0}}^{2}, \qquad 
{\frac {d U_{{1}}}{dY}}- 2\,U_{{0}}  U_{{1}}  = -{\frac {d^{2} U_{{0}}}{
		d{Y}^{2}}}.
  \label{eq:U0U1eqs}
\ee
Given (\ref{eq:Uexp1}), the required solutions are
\begin{equation}
U_0 = \frac{e^{Y}}{1 - e^{Y}}, 
\qquad 
U_1 =   {\frac {{{\rm e}^{Y}}+2\,\log  \left( 1-{{\rm e}^{Y}} \right)}{ \left( 1-{{\rm e}^{Y}} \right) ^{2}}}. \label{eq:riccati}
\end{equation}
Note that $U_0$ has a simple pole at $Y=0$, which does not match the type of singularity of the NLH found in section~\ref{sec:local}, while $U_1$ has a double pole at $Y=0$, which implies a further non-uniformity for small $Y$.
Since
\begin{equation}
U_1 \sim   {\frac {  1+2\,\log  \left( -Y \right)   }{{Y}^{2}}}, \qquad \mathrm{as\: } Y \to 0^{-}, \label{eq:Using}
\end{equation}
the resulting inner rescalings are
\begin{equation}
Y = 2t\log(1/t) + \zeta t, \qquad U = \frac{V}{t},  \label{eq:smalltimeYform}
\end{equation}
and, in these variables, the NLH becomes 
\begin{equation*}
{\frac {\partial ^{2}V}{\partial {\zeta}
		^{2}}}
+{\frac {\partial V}
	{\partial \zeta}}    - V^{2}
	+
t\Big( (2\log t + 2 - \zeta){\frac {\partial V}{\partial \zeta}} -2\,V  \Big)
+{t}^{2}{\frac {\partial V}{\partial t}} 
         =0.
\end{equation*}
Setting 
$V(\zeta,t) \sim V_0(\zeta)$
as $t\to 0$ with 
$\zeta = \mathcal{O}(1)$,
$V_0$ satisfies
\begin{equation}
 {\frac {d^{2}V_{0}}{d\zeta^{2}}}      +  {\frac{d V_{0} }
	{d\zeta}}   = V_{0}^{2}.
 \label{eq:V0}
\end{equation}
In order to match the behaviour (\ref{eq:Using}) for 
$\zeta \to -\infty$,
we require that $V_0$ satisfies the far-field condition
\begin{equation}
V_0 \sim -\frac{1}{\zeta} + \frac{2\log(-\zeta)}{\zeta^2} + \frac{1}{\zeta^2}, \qquad \mathrm{as\: } \zeta \to -\infty;  \label{eq:V0bc}
\end{equation}
prescribing the coefficient of the 
$1/\zeta^2$
term (which reflects the translation invariance of (\ref{eq:V0bc})) specifies $V_0$ uniquely (as can be confirmed by a phase-plane analysis, see Figure~\ref{fig:ODEphaseplot}). 
Making the change of variables $\zeta \mapsto -\zeta - 1$ and relabelling $V_0 \mapsto \phi$, (\ref{eq:V0})--(\ref{eq:V0bc}) is equivalent to (\ref{eq:V0m})--(\ref{eq:V0mbc}).  Unlike the first equation in (\ref{eq:U0U1eqs}), the singularities of (\ref{eq:V0}) have the same form as that of the NLH equation, see (\ref{eq:philocexpm}) and (\ref{eq:expansion}). 
Hence, no further rescalings are required.   Finally, we note that combining $u = U/t$, (\ref{eq:smallts}), (\ref{eq:smalltimeYform}) and $\zeta \mapsto -\zeta - 1$, $V_0 \mapsto \phi$, one obtains (\ref{eq:NLHcomplexsmallt}).

\section{A nonlinear ordinary differential equation}\label{sect:appode} The ODE
\begin{equation}
    \frac{d^2 \phi}{dx^2} - \frac{d \phi}{dx} = \phi^2 \label{eq:apnonlinODE}
\end{equation}
arises in the analysis of the small-time limit (see (\ref{eq:V0m})), the large-amplitude limit and the small-amplitude limit (see (\ref{eq:Uode}) and (\ref{eq:U0sol})).  Here we elaborate on some of its properties. 

\subsection{Phase plane} 
\begin{figure}[h]
    \centering
    \mbox{\includegraphics[width = 0.67\textwidth]{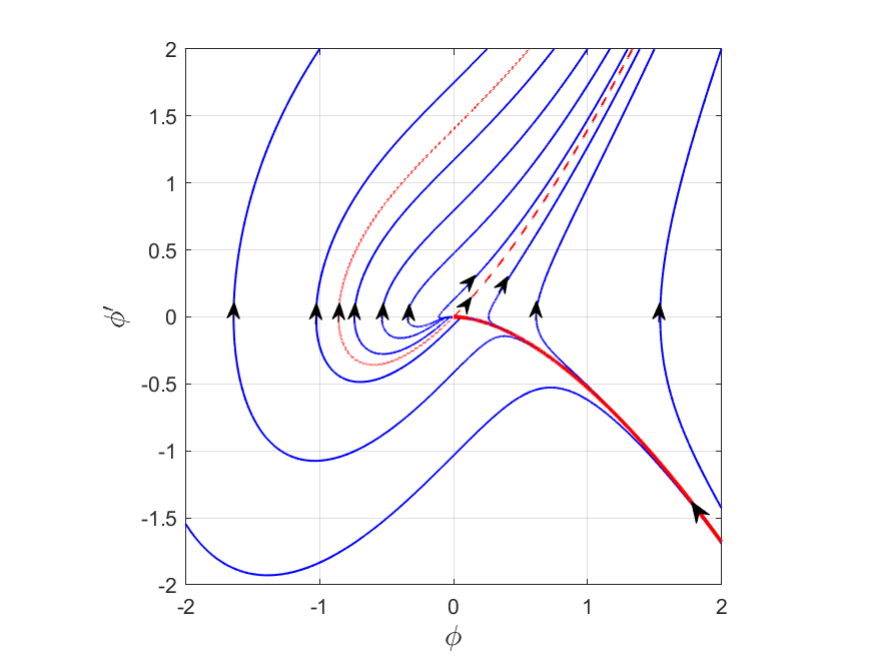}
    \includegraphics[width = 0.31\textwidth]{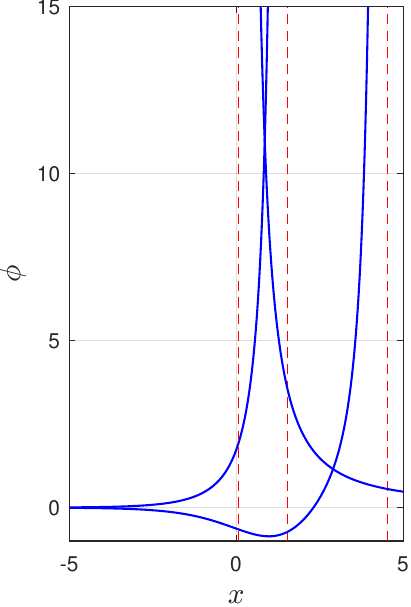}}
    \caption{  
    Left: phase plot of (\ref{eq:apnonlinODE}) for real-valued solutions on the real axis. The exceptional trajectories are indicated in red by a dashed line (\ref{eq:except1}), a dotted line (\ref{eq:except2}) and a thick solid line (\ref{eq:except3}). Right: the solutions corresponding to the exceptional trajectories with $a=1$ and $x_0 = 0$. The singularity locations are indicated by the dashed lines and occur at $x = 1.53767\ldots$ (for (\ref{eq:except1})), $x = 4.53879\ldots$ (for (\ref{eq:except2})) and $x = 0.05695\ldots$ (for (\ref{eq:except3})).
    }  \label{fig:ODEphaseplot}
\end{figure}

The phase plane of real solutions on the real line is immediately instructive (see the left frame of Figure~\ref{fig:ODEphaseplot}). 

There are precisely three exceptional trajectories,  each of which is relevant to our analysis, with the following behaviour:
\begin{equation}
    \phi \sim a e^{x}, \qquad x \to -\infty, \label{eq:except1}
\end{equation}
($a$ is a positive constant),
\begin{equation}
    \phi \sim -a e^{x}, \qquad x \to -\infty, \label{eq:except2}
\end{equation}
and
\begin{equation}
    \phi \sim \frac{1}{x}  + \frac{2\log x}{x^2} + \frac{x_0}{x^2} + \mathcal{O}\left( {x^{-3}} \right), \qquad
    x \to \infty \label{eq:except3}.
\end{equation}
For the latter case, $x_0$ is arbitrary and corresponds to a translation of $x$.  Notice that (\ref{eq:except1}) and (\ref{eq:except3}) characterise the solutions that arise in the analysis above (cf.~(\ref{eq:V0mbc}), 
(\ref{eq:Uodeic}) and (\ref{eq:U0sol})). The solutions specified by (\ref{eq:except1})--(\ref{eq:except3}) blow up on the real axis (see the right frane in Figure~\ref{fig:ODEphaseplot}) and in particular the singularities of (\ref{eq:except1}) and (\ref{eq:except3}) on 
the real line 
are of interest since these are used in the asymptotic estimates of the NLH singularity locations in (\ref{eq:sigmasmalltime}), (\ref{eq:sigmalargeamp}), (\ref{eq:sigmatimescale1m}) and (\ref{eq:sigmatimescale2m}). 

As noted in (\ref{eq:philocexpm}), the ODE (\ref{eq:apnonlinODE}) has movable logarithmic branch point singularities (which have the behaviour of second-order poles to leading order, however) and thus it fails the Painlev\'e test.


\subsection{Computing asymptotic boundary conditions}\label{sec:appbcscompute}  
\begin{figure}[h]
    \centering
    \mbox{
    \includegraphics[width = 0.5\textwidth]{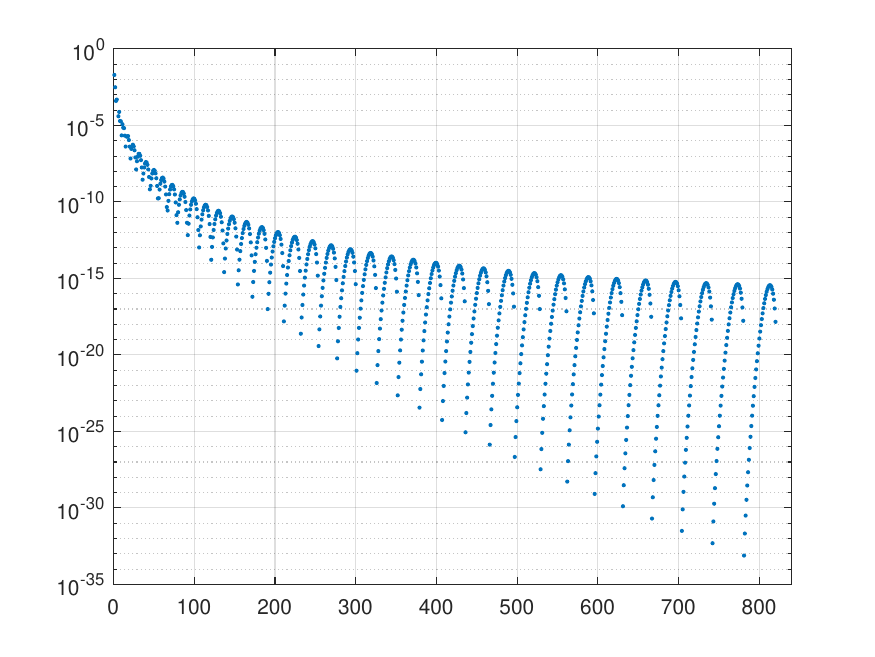}}
    \caption{  
    The magnitude of the terms $c_{k,j}(\log x)^{j-1-k}(-x)^{-j}$ in the asymptotic expansion (\ref{eq:odelogasexp}) for $x = 50$. The terms are plotted in the same order in which they appear in the series (\ref{eq:odelogasexp}).
    }  \label{fig:asexpterms}
\end{figure}

To  compute initial conditions for the solutions (\ref{eq:except1}) and (\ref{eq:except3}) accurately,  we generate higher-order terms in their asymptotic boundary conditions for fixed $x \ll -1$ or $x \gg 1$, respectively, and truncate the expansion when the first omitted term is below the level of machine precision. For (\ref{eq:except1}), the asymptotic expansion takes the form 
\begin{equation*}
\phi \sim \sum_{k = 1}^{\infty}c_ke^{kx}, \qquad x \to -\infty,
\end{equation*}
where $c_1$ is arbitrary (for the solution we require, $c_1 =1$ (see (\ref{eq:nonlinODEprob2}))
and 
\begin{equation}
c_{k+1} =\frac{1}{k(k+1)} \sum_{j=1}^{k}c_jc_{k+1-j},\qquad k \geq 1.  \label{eq:expexp}
\end{equation}
For (\ref{eq:except3}), the expansion takes the form
\begin{equation}
\phi \sim \sum_{j = 1}^{\infty}\left(\sum_{k = 0}^{j-1}c_{k,j}(\log x)^{j - 1 - k}\right)(-x)^{-j},\qquad  x \to \infty. \label{eq:odelogasexp}
\end{equation}
Comparing this expansion to the first few terms in (\ref{eq:except3}), one gets
\begin{equation*}
    c_{0,1} = -1, \qquad c_{0,2} = 2, \qquad c_{1,2} = x_0,
\end{equation*}
where $x_0$ is arbitrary (for the solution we require, $x_0 = 0$, see (\ref{eq:V0mbc})) and for $j = 4, 5, \ldots$, $k=0, \ldots, j-2$,
\begin{eqnarray*}
  \hspace{-1cm}  c_{k,j-1} =&\frac{1}{3-j}\left[ b_{k,j} - \left( j-1 \right)  \left( j-2 \right) c_{{k-1,j-2}}-
 \left( j-1-k \right) c_{{k-1,j-1}}+ \right.\\
& \hspace{1cm} \left.\left( 2\,j-3 \right)  \left( j-1-k \right) c_
{{k-2,j-2}}
- \left( j-1-k \right)  \left( j-k
 \right) c_{{k-3,j-2}}
\right],
\end{eqnarray*}
where
\begin{equation*}
    b_{k,j} = \sum _{r=2}^{j-2}  \sum _{p=\max \left( 0,k-j+r+1 \right) }^{
\min \left( k,r-1 \right) }c_{{p,r}}c_{{k-p,j-r}}.
\end{equation*}
Figure~\ref{fig:asexpterms} gives an example of the number of terms of (\ref{eq:odelogasexp}) that is required to compute a truncated expansion to an accuracy of roughly $10^{-15}$.

\subsection{Intial-value problem I}\label{sec:appivp1} The first initial-value problem for (\ref{eq:apnonlinODE}) with which we are concerned has initial condition (\ref{eq:except1}) with $a = 1$ (we impose this as $\mathrm{Re \,}x \to -\infty$ for all $\mathrm{Im\, }x$). On rays $\mathrm{Im\, }x = 2N\pi$ for integer $N$ the solution is given by (\ref{eq:except1});  $\phi$ is again real on $\mathrm{Im\, }(2N+1)\pi$ and (\ref{eq:except2}) applies on these rays. Indeed, $\phi$ is
$2\pi$-periodic in $x$.
Note from Figure~\ref{fig:ODEphaseplot} that if (\ref{eq:except1}) applies, the first singularity is encountered at $x_0 = 1.53767\ldots$, whereas for (\ref{eq:except2}), $x_0 = 4.53879\ldots$ and hence the same applies on the rays $\mathrm{Im\, }x = 2N\pi$ and $\mathrm{Im\, }x = (2N+1)\pi$. 

Figure~\ref{fig:ODE_exp_mod} shows the solution satisfying (\ref{eq:except1}) with $a=1$ in the complex plane, which was computed with the pole field 
solver~\cite{fasondini2017,fornberg2011}. The initial conditions for the solution, computed in the manner described above using (\ref{eq:expexp}), are
\begin{equation}
    \phi(-5) = 6.760698048065327\cdot10^{-3}, \qquad \phi'(-5) = 6.783500281706220\cdot10^{-3}. \label{eq:nonlinODE2ics} 
\end{equation}

To understand better the singularity structure shown in Figure~\ref{fig:ODE_exp_mod} we now analyze the behaviour as $\mathrm{Re\, }x \to +\infty$, which is expedited by the change of variables
\begin{equation}
    \xi = e^{x/5}, \qquad \phi = \frac{1}{25}e^{2x/5}v, \label{eq:apptrans1}
\end{equation}
to give
\begin{equation}
    \frac{d^2v}{d\xi^2} - \frac{6}{\xi^2}v = v^2, \qquad v \sim \xi^3 \quad \mathrm{as} \quad \vert \xi \vert \to 0. \label{eq:apxieq}
\end{equation}
Two remarks about this change of variables are in order. Firstly, it was identified on application of Kuzmak's method to determine the far-field asymptotic behaviour of (\ref{eq:apnonlinODE}); once identified, it seems more transparent to apply it at the onset. Secondly, disregarding the second term in (\ref{eq:apxieq}) and reversing the transformation gives
\begin{equation}
      \frac{d^2\phi}{dx^2} - \frac{d\phi}{dx} + \frac{6}{25}\phi = \phi^2  \label{eq:apphieq}
\end{equation}
so $x = -5z/\sqrt{6}$, $\phi = 6u/25$ yields
\begin{equation*}
   \frac{d^2u}{dz^2} + \frac{5}{\sqrt{6}}\frac{du}{dz} + u(1-u) = 0 
\end{equation*}
i.e. the transformation (\ref{eq:apptrans1}) corresponds to the inverse of that identified by~\cite{ablowitz1979} in constructing explicit solutions to Fisher's equation; the resemblance of (\ref{eq:apnonlinODE}) and (\ref{eq:apphieq}) for large $\phi$ is clear. 

The solution to (\ref{eq:apxieq}) satisfies
\begin{equation}
 v \sim   \frac{6^{1/3}}{\alpha^2}\wp\left( \frac{\xi + \xi_0}{6^{1/3}\alpha};0, 1 \right) - \frac{3}{\xi^2} + \mathcal{O}\left( \frac{1}{\xi^4}\right), \qquad \mathrm{as\: } \vert \xi \vert \to \infty, \label{eq:appvsol}
\end{equation}
where $\wp(z;0,1)$ denotes the equianharmonic case of the Weierstrass elliptic function~\cite[Ch.~23]{NIST:DLMF}. 
In (\ref{eq:appvsol}) the constants $\xi_0$ and $\alpha$ require numerical calculation 
and provide the requisite two degrees of freedom in the far-field behaviour. 
Thus the far-field singularities are given by those of the doubly periodic function $\wp$ which lie on the lattice
\begin{equation}
\hspace{-2cm}    \xi = -\xi_0 + 6^{1/3}\alpha\left[  2N\omega_1 + 2M\omega_3  \right], \qquad N, M \in \mathbb{Z}, \qquad \omega_1 = \frac{\left[\Gamma(1/3) \right]^{3}}{4\pi}, \omega_3 = e^{\pi i/3}\omega_1.  \label{eq:weierpoles}
\end{equation}
The only singularities of $\wp$ are second-order poles, hence (\ref{eq:philocexpm}) shows that the singularities of $\wp$ match those of $\phi$ (which are branch points) only to the leading order. Nevertheless, consistent with (\ref{eq:appvsol}), the bottom frame of Figure~\ref{fig:ODE_exp_mod} shows that the singularity locations of $\wp$ approximate those of $\phi$ shown in Figure~\ref{fig:ODE_exp_mod} .



\subsection{Initial-value problem II}\label{sec:appivp2} We now consider solutions satisfying (\ref{eq:except3}). Figure~\ref{fig:ODE_complex_phase_modulus} depicts the solution in the complex $x$-plane with $x_0 = 0$. The initial conditions for this solution are
\begin{equation}
\phi(50) = 2.359876891765835\cdot10^{-2}, \qquad
\phi'(50) = -5.332816483593814\cdot 10^{-4},  \label{eq:ODE2ics}
\end{equation}
and were computed using more than 800 terms from the expansion (\ref{eq:odelogasexp}), see Figure~\ref{fig:asexpterms}.
The singularity at $x_*=0.05695$ is the only singularity on the initial Riemann sheet. Integrating clockwise around the branch point onto the next sheet of the Riemann surface of the solution, however, we find more singularities, as shown in Figure~\ref{fig:ODE_complex_mod_sheet2_rotate}.

To clarify the singularity structure shown in Figures~\ref{fig:ODE_complex_phase_modulus} and~\ref{fig:ODE_complex_mod_sheet2_rotate}, we follow the trajectory (\ref{eq:except3}) through decreasing $x$. The Liouville--Green (or JWKB) method implies that for $x$ real, again up to translations in $x$ (translating out $x_0$ for brevity),
\begin{equation}
\phi \sim \frac{1}{x} + \frac{2\log x}{x^2} + \cdots + a x^2e^{x}, \qquad \mathrm{as\: } x \to -\infty,  \label{eq:apalgstokes}
\end{equation}
wherein $a$ is an arbitrary constant and the ellipsis denotes a divergent series in the usual way (the associated Stokes phenomenon will play an important role below).

Now despite the finite $x$ blow up at $x_*$, (\ref{eq:except3}) holds as $\vert x \vert \to \infty$ for all $\mathrm{arg}(x)$ on the initial Riemann sheet\footnote{Indeed, the exponential 
term in~(\ref{eq:apalgstokes})
is present only in the Stokes lines about $\arg(x) = \pm \pi$ wherein it is exponentially subdominant.}
and, as noted above, the numerical results in Figure~\ref{fig:ODE_complex_phase_modulus} suggest that $x=x_*$ is the only singularity on that Riemann sheet on which we let the branch cut be $(-\infty, x_*)$. The Stokes line on which the exponential in (\ref{eq:apalgstokes}) (wherein $a$ requires numerical calculation and will be imaginary) is turned on thus coincides with the branch cut, 
with the jump across the cut being balanced by the Stokes phenomenon. 
However, one can use (\ref{eq:apalgstokes}) in transitioning to the far field on the next Riemann sheet; hence we now impose (\ref{eq:apalgstokes})  with $a$ treated as given and rotate in the far field of that Riemann sheet to the anti-Stokes line, $\vert x \vert \to \infty$, $\mathrm{arg}(x) = \frac{\pi}{2}$; thus we perturb about the location $x = ir$ with $r\in \mathbb{R}$ in the form 
\begin{equation}
    x = ir - \log(-ar^3) - i(r\mathrm{\, mod}(2\pi)) + \widehat{X} \label{eq:apstokesdchange}
\end{equation}
to give the matching condition
\begin{equation}
    \phi \sim \frac{1}{r}\left(-i + e^{\widehat{X}}   \right), \qquad \mathrm{as\: } \widehat{X} \to -\infty, \qquad r \to \infty.  \label{eq:apphiasympt}
\end{equation}
The intermediate terms in (\ref{eq:apstokesdchange}) are included to attain this balance with the exponential having a coefficient of unity. The leading-order behaviour for $\widehat{X} = \mathcal{O}\left(1\right)$ is straightforward: since $\phi = \mathcal{O}(1/r)$ the right-hand side of (\ref{eq:apnonlinODE}) is negligible, so the leading-order solution simply reproduces (\ref{eq:apphiasympt}); the important scaling thus involves a further translation of $x$,
\begin{equation*}
    \widehat{X} = \log r + X
\end{equation*}
so the first terms in (\ref{eq:apalgstokes}) and (\ref{eq:apphiasympt}) become negligible and we are left with
\begin{equation*}
    \frac{d^2\phi}{dX^2} - \frac{d\phi}{dX} = \phi^2, \qquad \mathrm{as\: } X \to -\infty, \phi \sim e^{X}
\end{equation*}
at leading order and the above analysis associated with (\ref{eq:except1}) can be reused, albeit only at leading order. These results are consistent with the numerical solutions on the second Riemann sheet, see Figure~\ref{fig:ODE_complex_mod_sheet2_rotate}, which shows that the singularities are distributed according to those of the (equianharmonic) Weierstrass function in the variable $\xi$.

\section{Analysis of solutions with large-amplitude initial data}\label{sec:applargeamp}
Setting $\alpha = 1/\epsilon$ with $0 < \epsilon \ll 1$ in (\ref{eq:ic}), the rescaling
$t = \epsilon T$, $u = \epsilon^{-1} U$ transforms~(\ref{eq:pde}) 
and~(\ref{eq:ic}) into
\begin{equation}
\frac{\partial U}{\partial T} =  \epsilon \frac{\partial^2 U}{\partial x^2} + U^2, \label{eq:largeamp}
\end{equation}
and
\[
U(x,0) = \cos x.
\]
The zeroth-order approximation is therefore
\be
U_0 = \frac{\cos x}{1 - T \cos x}, \label{eq:U0largeamp}
\ee
with leading-order singularity locations $x = \pm i y = \pm i \sigma(T)$, where $\sigma(T) \sim \cosh^{-1}(1/T)$, $\epsilon \to 0$. 
 Moreover,
 \begin{equation*}
     \frac{\partial U_1}{\partial T} = \frac{\partial^2 U_0}{\partial x^2} + 2U_0U_1
 \end{equation*}
 leads to
 \begin{equation}
 U_1  = - \frac{\cos x (1 + 2\tan^2 x)T}{(1 - T \cos x)^2} - \frac{2\tan^2 x \log (1 - T\cos x)}{(1 - T\cos x)^2}.  \label{eq:largeampU1sol}
 \end{equation}

The singularities of (\ref{eq:U0largeamp}) are simple poles after a single rescaling and the estimate therefore ceases to be valid close to the singularities. Since (\ref{eq:U0largeamp}) and (\ref{eq:largeampU1sol}) imply
\begin{equation*}
\hspace{-1cm}U_0 \sim \frac{1}{T\sqrt{1 - T^2}(-Y)}, \qquad U_1 \sim \frac{2\log(-Y)}{(-Y)^2} + \frac{1-2T^2 + (1-T^2)\log(1 - T^2)}{(1-T^2)(-Y)^2},
\end{equation*}
as $Y \to 0$, where $Y = y - \cosh^{-1}(1/T)$, it follows that the inner scalings required to reproduce (\ref{eq:V0m})--(\ref{eq:V0mbc}) at leading order read
\begin{eqnarray*}
&Y = \epsilon T \sqrt{1-T^2}\left( 2\log(1/\epsilon) - \frac{1 - 2T^2}{1 - T^2} - 2\log T - 2\log(1 - T^2) - \zeta  \right) \\
& U \sim \frac{1}{\epsilon T (1 - T^2)}\phi(\zeta).
\end{eqnarray*}
Converting these expressions to the original variables, we obtain (\ref{eq:uzlargeamp})--(\ref{eq:largeampzzeta}).

\section{Analysis of solutions with small-amplitude initial data}\label{sec:appsmallamp}

Now setting $\alpha = \epsilon$ with $0 < \epsilon \ll 1$, the rescaling
$u = \epsilon w$ gives
\[
\frac{\partial w}{\partial t} =  \frac{\partial^2 w}{\partial x^2} + \epsilon w^2,
\]
with
\[
w(x,0) = \cos x.
\]
There are five time scales that need to be considered.

\subsection{$t = \mathcal{O}(1)$}

Setting $w = w_0 + \epsilon w_1 + \epsilon^2 w_2 + \cdots$, where $w_0(x,0) = \cos x$, $w_k(x,0) = 0$, $k \geq 0$, one obtains a sequence of linear PDEs:
\begin{equation*}
    {\frac {\partial w_{{0}}}{\partial t}}  = {\frac {
\partial ^{2}w_{{0}}}{\partial {x}^{2}}},  \qquad
 {\frac {\partial w_{{1}}}{
\partial t}} =  {\frac {\partial ^{2}w_{{1}}}{
\partial {x}^{2}}}   +   w_{{0}}^{2},  \qquad
      {\frac {
\partial w_{{2}}}{\partial t}} = {\frac {\partial ^{
2}w_{{2}}}{\partial {x}^{2}}} + 2\,w_{{0}}  w_{{1}}.
\end{equation*}
Solving these, we find 
\begin{eqnarray}
  & \hspace{-2cm} w  \sim  \mathrm{e}^{-t}\cos x + \frac{\epsilon}{4}\left[ 1-\,{{\rm e}^{-2\,t}}+\, \left( {{\rm e}^{-2\,t}}-{{\rm e}^{-4
\,t}} \right) \cos  2\,x  \right] \nonumber \\
 &\hspace{-1.5cm}  + \frac{\epsilon^2}{48}\left[ \left( 24 t+6\,{{\rm e}^{-2\,t}}+3\,{{\rm e}^{-4\,t}}-9
 \right) {{\rm e}^{-t}}\cos  x  + \left( 2 -3\,{{\rm e}^{-2\,t}}+{{\rm e}^{-6\,t}}\right) {{\rm e}^{-3\,t}}\cos
  3\,x   \right].
 \label{eq:smallampexp}
\end{eqnarray}
Setting $x = i y$ and moving up the imaginary axis, a non-uniformity occurs at
\begin{equation*}
    \frac{\epsilon\left( {{\rm e}^{-2\,t}}-{{\rm e}^{-4
\,t}} \right) \cos  2\,x }{4 e^{-t}\cos x} \sim  \frac{\epsilon}{2}\,{{\rm e}^{y-2\,t}}\sinh t  = \mathcal{O}(1).  
\end{equation*}
This motivates the change of variable  
\begin{equation}
    y = 2t - \log(\sinh t) - \log(\epsilon/2) + \zeta. 
    \label{eq:Yform}
\end{equation}
In this variable, the expansion (\ref{eq:smallampexp}) yields
\begin{eqnarray}
u & 
\sim  {\frac {{
{\rm e}^{\zeta+t}}}{\sinh  t  }}+   \left( {\frac {{{\rm e}^{2\,t}}-1}{ 2 \sinh^{2} t}}
 \right)  {{\rm e}^{2\zeta}} +\left({\frac {  2\,{{\rm e}^{3\,t}}-3\,{{\rm e}^{\,t}}+ {{\rm e}^{-3\,t}}}{  12\sinh^{
3}t}}
 \right) {{\rm e}^{3\zeta}} \nonumber \\
 & + \frac{\epsilon^2}{4} \left[ \left(1 -   {{\rm e}^{-2t}}\right) + \left( {\frac {  8\,t{{\rm e}^{\,t}}-3\,{{\rm e}^{\,t}}+2\,{
{\rm e}^{-\,t}}+{{\rm e}^{-3\,t}}}{4 \sinh t  }}
 \right) {{\rm e}^{\zeta}}  \,\right],
\label{eq:uY}
\end{eqnarray}
which will provide matching conditions as $\zeta \to -\infty$.

Setting $u\sim U(\zeta,t)$ in the NLH gives
\begin{equation}
{\frac {\partial U}{\partial t}}
  -   \left( 2- \coth t
 \right) {\frac {\partial U}{\partial \zeta}} =  -{\frac {\partial ^{2} U}{\partial {\zeta}^{2}}}+  
U^{2}.   \label{eq:backdiff1}
\end{equation}
From (\ref{eq:uY}), the initial and boundary conditions are 
\begin{equation}
\hspace{-2cm} U \sim \frac{{\rm e}^{\zeta}}{t}, \qquad \mathrm{as\: }t \to 0^{+}, \zeta = \mathcal{O}(1), \qquad    U \sim \frac{{\rm e}^{t}}{\sinh t}{\rm e}^{\zeta}, \qquad \mathrm{as\:} \zeta \to -\infty,  t = \mathcal{O}(1). \label{eq:backdiffbcs}
\end{equation}

A leading-order approximation for the location of the singularity is 
\begin{equation}
    \sigma(t) \sim 2t - \log(\sinh t) - \log(\epsilon/2) + \zeta_*(t), \qquad \epsilon \to 0, \label{eq:smallampsigmaf1}
\end{equation}
where extracting $\zeta_*(t)$, the first singularity of (\ref{eq:backdiff1})--(\ref{eq:backdiffbcs}) on the real axis, would require the numerical solution to this backward--diffusion problem, which, in some respects, is more involved than numerically solving the NLH equation. We note that the approximation (\ref{eq:smallampsigmaf1}) is consistent with (\ref{eq:sigmasmalltime}) (i.e., $\sigma \sim \log(2/(\alpha t))$) in the small-$t$ limit. 

The limiting behaviour of $\zeta_*(t)$ can be determined as follows. As $t \to 0^+$, 
\[
U \sim \frac{e^{\zeta}}{t(1 - \mathrm{e}^{\zeta})},
\]
for $\zeta = \mathcal{O}(1)$, in accord with the analysis of \ref{sec:appsmalltime} 
(cf.~(\ref{eq:riccati})), so that $\zeta_*(0)=0$. As $t \to \infty$ according to $t = \mathcal{O}(1/\epsilon)$, $U$, and hence $\zeta_*$, becomes independent of $t$, so that
\be
\frac{\partial^2 U}{\partial \zeta^2} \sim \frac{\partial U}{\partial \zeta} + U^2, \label{eq:Uode}
\ee
and $\zeta_*$ is determined from the solution to this ODE subject to the initial data
\be
U \sim 2\mathrm{e}^{\zeta}, \qquad \mathrm{as\:} \zeta \to -\infty.   \label{eq:Uodeic}
\ee
Changing variables according to $\zeta \mapsto \zeta - \log 2$ and $U \mapsto \phi$ in (\ref{eq:Uodeic}) and (\ref{eq:smallampsigmaf1}), one obtains (\ref{eq:nonlinODEprob2}) and (\ref{eq:NLHcomplextimescale1}).

\subsection{$t = \mathcal{O}\left( \epsilon^{-2} \right)$}\label{sec:smallamptimescale2}

The secular term in (\ref{eq:smallampexp}) suggests the introduction of the long timescale $t = T/\epsilon^{2}$. From (\ref{eq:smallampexp}) it also implies the matching condition $u \sim \epsilon^2/4$. Hence, we 
set $t = T/\epsilon^{2}$, $u = \epsilon^2 v$ to give
\begin{equation}
    \epsilon^2 \frac{\partial v}{\partial T} = \frac{\partial^2 v}{\partial x^2} + \epsilon^2 v^2.  \qquad \label{eq:nlht2}
\end{equation}
Since matching back into the expansion (\ref{eq:smallampexp}) implies that the solution is uniform in $x$ up to terms exponentially small in $\epsilon$ (making tracking non-uniformities numerically challenging) we have
\be
v \sim \frac{1}{T_c - T},
\ee
as $T \to T_c$ with $T_c \sim 4$, so a leading order approximation to the blow-up time is $t \sim 4/\epsilon^2$.
Next including the first Fourier mode via $v \sim 1/(T_c   - T) + a_1(T)\cos x$ ($\cos x$ being associated with the slowest decaying exponential term in (\ref{eq:smallampexp})) implies
\begin{equation}
    \epsilon^2 \dot{a_1} \sim -a_1 + \frac{2\epsilon^2a_1}{T_c - T}, \qquad \mbox{so \ that} \qquad a_1 \sim \frac{16}{\epsilon}\frac{{\rm e}^{-T/\epsilon^2}}{(T_c - T)^2},  \label{eq:a1s}
\end{equation}
in keeping with the above asserted exponential smallness, the arbitrary  constant in the general solution of $a_1$ being determined at leading order by matching back to the $\cos x$ term on the previous timescale in (\ref{eq:smallampexp}).

Proceeding further to the next Fourier mode in the solution expansion, $a_2(T)\cos 2 x$, one finds 
\begin{equation*}
    \epsilon^2\, \dot{a_{{2}}}\sim 
-4\,a_{{2}}  +\frac{{\epsilon}^{2}}{2}  a_{{1}}^{2}+\frac{2{\epsilon}^{2} a_{{2}}}{T_c-T},
\end{equation*}
with solution on matching back as $T \to 0^{+}$,
\begin{equation*}
    a_2 \sim \frac{128}{\epsilon^2}\frac{{\rm e}^{-4T/\epsilon^2}}{(T_c-T)^2}\int_{-\infty}^{T} \frac{{\rm e}^{2\widetilde{T}/\epsilon^2}}{(T_c-\widetilde{T})^2} \mathrm{d}\widetilde{T}, 
\end{equation*}
which implies via integration by parts that 
\begin{equation}
    a_2 \sim \frac{64}{(T_c - T)^4}{\rm e}^{-2T/\epsilon^2}, \qquad \label{eq:a2s}
\end{equation}
for $T = \mathcal{O}(1)$, matching back correctly to (\ref{eq:smallampexp}).
Hence, the solution expansion takes the form
\begin{equation}
    v \sim \frac{1}{T_c - T} + \left(\frac{16\, {\rm e}^{-T/\epsilon^2}}{\epsilon(T_c - T)^2} + \cdots \right)\cos x +  \left( \frac{64\, {\rm e}^{-2T/\epsilon^2}}{(T_c - T)^4}+ \cdots \right)\cos 2x.  \label{eq:phiexp}
\end{equation}
Setting $x = iy$ and moving up the imaginary axis, this expansion becomes non-uniform where $a_1 \cos x$ and $a_2 \cos 2x$ are of the same order (which is why we have included the $\cos 2x$ terms in the above calculation), implying the rescalings,
%
\begin{equation}
    y = \frac{T}{\epsilon^2} + \log\left( \frac{(T_c - T)^2}{8\epsilon}  \right)  + Y, \qquad v \sim \frac{U(Y) }{\epsilon^2}. \label{eq:yform}
\end{equation}
In this variable, the matching condition 
\begin{equation}
    U \sim {\rm e}^{Y} + \frac{1}{2}{\rm e}^{2Y}, \qquad \mathrm{as\: }Y\to -\infty,
    \label{eq:ut2rescaleexp}
\end{equation}
pertains, the  $a_1 \cos x$ and $a_2 \cos 2x$ terms dominating $a_0 = 1/(T_c - T)$. By (\ref{eq:nlht2}) and (\ref{eq:ut2rescaleexp}),
$U$ satisfies
\begin{equation}
    U^{\prime \prime} = U^{\prime} + U^2, \qquad \mathrm{with\: } U \sim {\rm e}^{Y} \ \mathrm{as\: }Y \to -\infty.  \label{eq:U0sol} 
\end{equation}
This is the same ODE as in~(\ref{eq:Uode}), but with a different boundary condition, cf.~(\ref{eq:Uodeic}).

\subsection{$T = T_c + \epsilon^2 \tau$}\label{sec:apptimescale3}

For $\tau = \mathcal{O}(1)$,
the expansion (\ref{eq:phiexp}) implies that the terms we need to take account of here take the form
\begin{equation}
    u = \epsilon^2\phi \sim \frac{1}{(-\tau)} + \frac{16\, {\rm e}^{-T_c/\epsilon^2}{\rm e}^{-\tau}}{\epsilon^3 (-\tau)^2}\cos x +   A(\tau) \frac{{\rm e}^{-2T_c/\epsilon^2}}{\epsilon^6}\cos2x. \label{eq:uexpt3}
\end{equation}
Substituting this into the NLH, the modulation factor $A(\tau)$ satisfies
\begin{equation*}
    \dot{A} = -4A + \frac{2}{(-\tau)}A + \frac{128}{(-\tau)^4} {\rm e}^{-2\tau},
\end{equation*}
and the solution that matches back to the expansion (\ref{eq:phiexp}) is
\begin{equation*}
 A =    128\,{\frac {{{\rm e}^{-4\,\tau}}}{{(-\tau)}^{2}}\int_{-\infty }^{\tau}\!{\frac {{
{\rm e}^{2\,\widetilde{\tau}}}}{{(\widetilde{\tau})}^{2}}}\,{\rm d}\widetilde{\tau}}.
\end{equation*}
The large- and small-$\tau$ behaviour of $A$ read 
\begin{equation*}
    A \sim \frac{64\, {\rm e}^{-2\tau}}{(-\tau)^4}, \qquad \mathrm{as\: } \tau \to -\infty, \qquad A \sim \frac{128}{(-\tau)^3}, \qquad \mathrm{as\: } \tau \to 0^{-}.
\end{equation*}
The rescaling needed to capture the nearest singularity in $y>0$ is thus 
\begin{equation}
    y = \frac{T_c}{\epsilon^2} + \tau - \log\left( \frac{8}{\epsilon^3 (-\tau)^2} \right) + Y,  \label{eq:t3yform}
\end{equation}
because (\ref{eq:uexpt3}) then becomes 
\begin{equation}
    u \sim \frac{1}{(-\tau)} + {\rm e}^{Y} + \frac{A(\tau)}{128}(-\tau)^4{\rm e}^{2Y + 2\tau}, \label{eq:ut3exp}
\end{equation}
and (\ref{eq:ut2rescaleexp}) is recovered in the limit $\tau \to -\infty$. What makes this $\tau$ time scale a distinguished limit is the fact that all three terms in (\ref{eq:ut3exp})
are of the same order,
whereas on the previous time scale only two of the terms were.
With the change of variables (\ref{eq:t3yform}), the NLH becomes
\begin{equation}
 {\frac 
{\partial U}{\partial \tau}} + \left( \frac{2}{(-\tau)} - 1\right)\,\frac {\partial U }{\partial Y} =  -{\frac {\partial ^
{2} U}{\partial {Y}^{2}}}  + U^{2},  \label{eq:PDEsmallt3}
\end{equation}
so that a full balance in the PDE again occurs for $\tau = \mathcal{O}(1)$
and we require that the solution satisfy
\begin{equation}
    U \sim \frac{1}{(-\tau)} + {\rm e}^{Y}, \qquad \mathrm{as\: } Y \to -\infty.   \label{eq:t3UYlarge}
\end{equation}
As $\tau \to -\infty$, steady state behaviour is required to match the previous timescale, i.e., (\ref{eq:U0sol}) again applies. 

In the limit $\tau \to 0^{-}$, the change of variables 
\begin{equation}
    Y = -\log(-\tau) + z, \qquad U = \frac{V}{(-\tau)}, \label{eq:t3Yrescale2}
\end{equation}
brings the two terms in (\ref{eq:t3UYlarge}) into balance, so that 
\begin{equation}
    \frac{d V_0}{dz} + V_0 = V_0^2, \qquad V_0 = \frac{1}{1 - {\rm e}^{z}}.  \label{eq:smalltsc3v0eq}
\end{equation}
We note that (\ref{eq:smalltsc3v0eq}) implies the presence of a finer scale to capture the double-pole nature of the singularity. Indeed, the rescalings
\begin{equation*}
z = (-\tau)(-Z), \qquad V = W/(-\tau)
\end{equation*}
imply the leading-order balance
\begin{equation*}
\frac{\partial^2 W}{\partial Z^2} - \frac{\partial W}{\partial Z} \sim W^2, \qquad W \sim 1/Z, \qquad \mathrm{as\: } Z \to \infty,
\end{equation*}
equivalent to (\ref{eq:V0m})--(\ref{eq:V0mbc}) but requiring the calculation of higher-order terms to fix the translation invariance with respect to $Z$, which we shall not pursue (but which will, as usual, involve a logarithmic contribution).

Combining (\ref{eq:t3yform}) and (\ref{eq:t3Yrescale2}) with $z = 0$ (the singularity location of $V_0$ in (\ref{eq:smalltsc3v0eq})) and converting to the original variables, we obtain the estimate (\ref{eq:sigmatimescale3m2}).

\subsection{$\tau = \frac{{\rm e}^{-T_c/\epsilon^2}}{\epsilon^3}\widehat{\tau}$}

In the limit $\tau \to 0^{-}$ as the blow-up time is approached, the expansion (\ref{eq:uexpt3}) becomes
\begin{equation}
    u \sim \frac{1}{(-\tau)} + \frac{16\, {\rm e}^{-T_c/\epsilon^2}}{\epsilon^3 (-\tau)^2}\cos x  + \frac{128\, {\rm e}^{-2T_c/\epsilon^2}}{\epsilon^6 (-\tau)^3}\cos 2x.  \label{eq:uexpt4}
\end{equation}
Under the rescalings
\begin{equation}
    \tau = \frac{{\rm e}^{-T_c/\epsilon^2}}{\epsilon^3}\widehat{\tau}, \qquad 
    u = \epsilon^3 {\rm e}^{T_c/\epsilon^2}\widehat{u},  \label{eq:t4vars}
\end{equation}
the expansion (\ref{eq:uexpt4}) implies
\begin{equation}
 \widehat{u} \sim \frac {1 }{(-\widehat{\tau})} + \frac{16}{(-\widehat{\tau})^2}\cos x  + \frac{128}{(-\widehat{\tau})^3}\cos 2x.
 \label{eq:utausmall}
\end{equation}
while the NLH becomes
\begin{equation*}
    \frac{\partial \widehat{u}}{\partial \widehat{\tau}} = \frac{{\rm e}^{-T_c/\epsilon^2}}{\epsilon^3} \frac{\partial^2 \widehat{u}}{\partial x^2} + \widehat{u}^2.
\end{equation*}
To leading order,
\begin{equation*}
    \frac{\partial \widehat{u}_0}{\partial \widehat{\tau}} =   \widehat{u}^2_0, 
\end{equation*}
so that
\begin{equation}
     \widehat{u}_0 = \frac{1}{-\widehat{\tau} - 16\cos x}, \label{eq:uhat0}
\end{equation}
in view of (\ref{eq:utausmall}), which represents a matching condition in the limit $\widehat{\tau} \to -\infty$. 

The fifth time scale is discussed in section~\ref{sec:timescale5}.

\section{Generic and non-generic blow-up behaviour}\label{sect:appsectblowup}




In order to keep our analysis self-contained we here re-derive well-known results (see the review articles cited in section~\ref{sec:intro}), though some aspects of what follows are novel (including the extension into the complex plane) and the perspective may at times differ.

We start from the well-established observation 
that blow up is generically close, in a sense that becomes explicit in what follows, to the scale-invariant self-similar reduction of the NLH; 
i.e.~choosing the origins of space and time such that blow up occurs at $(x, t) = (0, 0)$, we set
\begin{equation}
u = \frac{f(\eta,\tau)}{-t}, \quad \eta = \frac{x}{(-t)^{1/2}}, \quad \tau = \log(1/(-t)),  \label{eq:buvars}
\end{equation}
to give the $\tau$-translation-invariant PDE
\begin{equation*}
    \frac{\partial f}{\partial \tau} + f + \frac{1}{2}\eta \frac{\partial f}{\partial \eta} = \frac{\partial^2 f}{\partial \eta^2} + f^2
\end{equation*}
that we wish to analyse in the limit $\tau \to +\infty$. This will require us to consider three different scales, namely $\eta = \mathcal{O}(1)$, $\eta = \mathcal{O}(\tau^{1/2})$ and $\log \eta = \mathcal{O}(\tau)$. 

\subsection{$\eta = \mathcal{O}(1)$}\label{sec:appbuscale1} On the first scale we exploit the known property  (that we further evidence below) that to leading order the blow-up profile for $\eta = \mathcal{O}(1)$ is flat, $f \sim 1$, by writing
\begin{equation}
    f = 1 + F \label{apeq:fF}
\end{equation}
with $F \to 0$ as $\tau \to \infty$ and 
\begin{equation}
    \frac{\partial F}{\partial \tau} - F + \frac{1}{2}\eta \frac{\partial F}{\partial \eta} = \frac{\partial^2 F}{\partial \eta^2} + F^2. \label{apeq:F}
\end{equation}
Since $F \to 0$ it is natural first to discard the nonlinear term in (\ref{apeq:F}) and consider the role of the Hermite polynomial solutions to the heat equation, the first two of which (the zeroth and first eigenmodes),
\begin{equation*}
 e^{\tau}, \qquad e^{\tau/2}\eta,   
\end{equation*}
are inadmissible since their growth in $\tau$ would invalidate the assumption implicit in (\ref{apeq:fF}) (they in any case simply represent, respectively, shifts in the time and location of blow up). The next three eigenmodes are
\begin{equation}
    1-\eta^2/2, \quad e^{-\tau/2}\left(\eta - \eta^3/6 \right), \quad e^{-\tau}\left(1 - \eta^2 + \eta^4/12 \right). \label{apeq:eigenmodes}
\end{equation}
The first of (\ref{apeq:eigenmodes}) (i.e., the second eigenmode) is borderline (`neutrally stable') with respect to (\ref{apeq:fF}) and it is this borderline property that is responsible for the algebraic dependence upon $\tau$ that arises in and complicates what follows. The first in the hierarchy of non-generic blow up scenarios occurs when the second mode is absent, in which case the third must also be since $\eta = \mathcal{O}(1)$ is required to contain the maximum of $u$, and a multiple of the third expression in (\ref{apeq:eigenmodes}) then dominates $F$, with the second term in (\ref{apeq:fF}) exponentially smaller than the first as $\tau \to +\infty$; the need for such a non-generic blow up scenario can be clarified by considering two-peaked initial data symmetric about $x = 0$: for sufficiently widely spaced peaks generic blow up will occur simultaneously at two locations (see the left frame of Figure~\ref{fig:2peaksreal}), while for closely spaced ones diffusion will result in generic blow up at $x = 0$ only  (as in the right frame of Figure~\ref{fig:2peaksreal}). Non-generic blow up  occurs as the borderline between these two scenarios and we return to such matters below.

Focusing now on the generic case, in light of the above insights we seek a solution to (\ref{apeq:F}) of the form
\begin{equation}
    F \sim a_0(\tau)\left( 1 - \eta^2/2   \right) + a_1(\tau)F_1(\eta), \qquad \tau \to +\infty,  \label{apeq:Fexp}
\end{equation}
with $a_1 \ll a_0 \ll 1$ in that limit. Hence
\begin{equation}
    a_1\left( \frac{d^2 F_1}{d\eta^2} -  \frac{1}{2} \eta \frac{dF_1}{d\eta} + F_1    \right) = \frac{d a_0}{d\tau}\left( 1 - \eta^2/2  \right) - a_0^2\left( 1 - \eta^2/2  \right)^2. \label{apeq:F1}
\end{equation}
The requirement that $F_1$ not grow exponentially as $\eta \to \pm \infty$ leads to a solvability condition on $a_0(\tau)$, namely 
\begin{equation}
    \frac{d a_0}{d\tau} = - 4a_0^2,  \label{apeq:a0solve}
\end{equation}
so that
\begin{equation*}
a_0 = \frac{1}{4(\tau + \tau_0)} \sim \frac{1}{4\tau}, \qquad \tau \to +\infty,
\end{equation*}
for some constant $\tau_0$ that reflects the scale invariance of the NLH and is dependent on the initial data; the solution to (\ref{apeq:F1}) then takes the form
\begin{equation}
    a_1(\tau)F_1(\eta) = -\frac{5}{32(\tau +\tau_0)^2}\eta^2 + \frac{1}{64(\tau +\tau_0)^2}\eta^4 + \widehat{a}_1(\tau)\left(1 - \frac{1}{2}\eta^2  \right)  \label{apeq:Fa0a1}
\end{equation}
where the calculation of $\widehat{a}_1(\tau)$ requires a similar solvability condition at the next order of the expansion and, without pursuing those details, it can readily be shown that
\begin{equation}
    \widehat{a}_1 = \frac{\alpha_1}{(\tau + \tau_0)^2},  \label{apeq:a1hat}
\end{equation}
with $\alpha_1$ constant.

\subsection{$\eta = \mathcal{O}(\tau^{1/2})$}\label{sec:appbuscale2} We turn next to the outer region $\zeta = \mathcal{O}(1)$, where\footnote{We remark that as $\tau \to +\infty$ the constant $\tau_0$ contribution can be disregarded in such expressions, but since $\tau$ appears throughout only in the combination $\tau + \tau_0$ there is value in consistently retaining it.}
\begin{equation*}
    \zeta = \frac{\eta}{(\tau + \tau_0)^{1/2}}
\end{equation*}
Hence
\begin{equation}
    \frac{\partial f}{\partial \tau} - \frac{\zeta}{2(\tau + \tau_0)}\frac{\partial f}{\partial \zeta} + f + \frac{1}{2}\zeta \frac{\partial f}{\partial \zeta} = \frac{1}{\tau+\tau_0}\frac{\partial^2 f}{\partial \zeta^2} + f^2  \label{apeq:fzeta}
\end{equation}
and the expansion
\begin{equation*}
    f \sim f_0(\zeta) + \frac{1}{\tau+\tau_0}f_1(\zeta)
\end{equation*}
pertains. Therefore
\begin{equation*}
    f_0 + \frac{1}{2}\zeta \frac{d f_0}{d \zeta} = f_0^2
\end{equation*}
and
\begin{equation}
    f_0 = \frac{1}{1 + \alpha_0\zeta^2/2 }, \qquad \alpha_0 = \frac{1}{4},  \label{apeq:f0sol}
\end{equation}
follows on matching into (\ref{apeq:a0solve}) and (\ref{apeq:Fexp}). Combining (\ref{apeq:f0sol}) and (\ref{eq:buvars}) and converting to the original variables with blow up occurring at $t = t_c$, one obtains (\ref{eq:ublowupzeta}).  

We proceed to next order, in part because the solvability condition that determines $\alpha_0$ is sometimes derived on this outer scale instead, so it seems helpful to clarify its presence.  At next order (\ref{apeq:fzeta}) implies
\begin{equation*}
    \frac{1}{2}\zeta \frac{d f_1}{d\zeta} + f_1 - 2f_0f_1 = \frac{d^2 f_0}{d \zeta^2} + \frac{\zeta}{2}\frac{d f_0}{d \zeta},
\end{equation*}
with solution
\begin{equation}
    f_1 = \frac{(4\alpha_0^2 - \alpha_0)\zeta^2\log \zeta + \alpha_0 - 2 \alpha_0^2\zeta\log(1 + \alpha_0\zeta^2/2) + \beta_1\zeta^2 }{(1 + \alpha_0\zeta^2/2)^2}  \label{apeq:f1sol}
\end{equation}
for constant $\beta_1$. Appealing to analyticity again implies $\alpha_0 = 1/4$ (indeed, (\ref{apeq:a0solve}) can be derived in this fashion without pre-specification of the dependence on $\tau$), but such an appeal is questionable given that the non-analyticity could in principle be smoothed over the inner scale: the issue can conveniently be clarified by differentiating (\ref{apeq:Fa0a1}) twice and introducing $G_1 = \frac{1}{a_1}\frac{d^2 F_1}{d\eta^2}$ to give
\begin{equation*}
    \frac{d^2 G_1}{d\eta^2} - \frac{1}{2}\eta \frac{d G_1}{d \eta} = -\frac{da_0}{d\tau} + (2 - 3\eta^2)a_0^2.
\end{equation*}
Setting
\begin{equation}
    G_1 = 3a_0^2\eta^2 - \left(\frac{da_0}{d\tau} + 4a_0^2   \right)H  \label{apeq:G1sub}
\end{equation}
implies
\begin{equation*}
    \frac{d^2 H}{d\eta^2} - \frac{1}{2}\eta\frac{dH}{d\eta} = 1
\end{equation*}
even solutions to which satisfy
\begin{equation*}
    \frac{dH}{d\eta} = \sqrt{\pi}e^{\eta^2/4}\mathrm{erf}(\eta/2) \sim \sqrt{\pi}e^{\eta^2/4} - 2/\eta, \quad \eta \to +\infty,
\end{equation*}
so that $H$ contains contributions as $\eta \to +\infty$ of both the forms
\begin{equation}
    \frac{2\sqrt{\pi}}{\eta}e^{\eta^2/4}, \qquad -2\log \eta; \label{apeq:Hforms}
\end{equation}
eliminating the former by requiring $a_0$ in (\ref{apeq:G1sub}) to satisfy (\ref{apeq:a0solve}) necessarily simultaneously removes the latter. Accordingly, analyticity in (\ref{apeq:f1sol}) is subsidiary to, and predicated upon, solvability in the inner region, the latter eliminating the highly problematic exponentially growing term in (\ref{apeq:Hforms}), to which the logarithmic term is inextricably tied. Such comments have implications for the various concise derivations of the form of the blow up in~\cite{Keller1993} that adopts $1/u$ as the dependent variable and derives a low-dimensional representation via a Taylor expansion (cf.~(\ref{apeq:Fexp}))\footnote{We remark in passing that truncating the Taylor expansion at higher order gives a certain amount of insight into non-generic blow up of the type we discuss below. }: this approach implicitly involves the outer scale only and its success is again indirectly a (fortuitous) consequence of the fact that suppressing the second term in (\ref{apeq:Hforms}) necessarily suppresses the first also. Finally, matching (\ref{apeq:f1sol}) with (\ref{apeq:a1hat}) and (\ref{apeq:Fa0a1}) implies
\begin{equation}
    \beta_1 = -\frac{3}{32} - \frac{\alpha_1}{2},  \label{eq:beta1v}
\end{equation}
the $\eta^4$ term in (\ref{apeq:Fa0a1}) already being captured by expanding (\ref{apeq:f0sol}) in its inner limit. Hence evaluation of $\beta_1$ requires the higher-order calculation necessary to determine $\alpha_1$.

\subsection{$\log \eta = \mathcal{O}(\tau)$}\label{sec:appblowup3rdscale}The third and final real-line scaling we need to pursue sets
\begin{equation*}
    f(\eta,\tau) = \frac{\phi(\sigma, \tau)}{\eta^2}, \qquad \sigma = \log \eta
\end{equation*}
to give
\begin{equation}
    \frac{\partial \phi}{\partial \tau} + \frac{1}{2}\frac{\partial \phi}{\partial \sigma} = e^{-2\sigma}\left(\frac{\partial^2 \phi}{\partial \sigma^2} - 5\frac{\partial \phi}{\partial \sigma} + 6 \phi + \phi^2  \right)  \label{apeq:phieq}
\end{equation}
the relevant scaling is $\sigma = \mathcal{O}(\tau)$, $\sigma > 0$, so the right-hand side of (\ref{apeq:phieq}) is exponentially small, corresponding to 
\begin{equation*}
    u \sim u(x,0)
\end{equation*}
applying as $t \to 0^{-}$ away from the blow-up point $x=0$, and we now obtain the local behaviour of $u(x,0)$. Matching into  (\ref{apeq:f0sol}) and (\ref{apeq:f1sol})   (with $\alpha_0 = 1/4$) implies the matching condition
\begin{equation*}
    \phi \sim 8(\tau + \tau_0) - 16 \sigma + 8 \log(8(\tau+\tau_0)) + 64\beta_1,
\end{equation*}
so disregarding the right-hand side of (\ref{apeq:phieq}) gives
\begin{equation*}
    \phi \sim 8(\tau + \tau_0) - 16 \sigma + 8 \log(8(\tau+\tau_0) - 16\sigma) + 64\beta_1;
\end{equation*}
since $u = \phi/x^2$, $\sigma = \log x + \tau/2$ this implies that
\begin{equation}
    u(x,0) \sim \frac{8}{x^2}\left( 2\log(1/x) + \log(\log(1/x)) + 4\log 2 + \tau_0 + 8\beta_1 \right), \quad x \to 0^{+}.  \label{apeq:uexp}
\end{equation}
The $\tau_0$ dependence of (\ref{apeq:uexp}) captures the leading-order dependence on the initial data (except of course for the time and location of blow up).

\subsection{Complex-plane behaviour} One further region is needed to complete the complex-plane analysis in the current context. The relevant terms from (\ref{apeq:f0sol}) and (\ref{apeq:f1sol}) for the purposes of matching give, as $\zeta^2 \to -8$ that 
\begin{equation*}
    f \sim \frac{1}{1 + \zeta^2/8} + \frac{1}{\tau+\tau_0}\frac{\log(1+\zeta^2/8) + 1/4 - 8\beta_1}{\left(1 + \zeta^2/8 \right)^2}.
\end{equation*}
The inner scalings 
\begin{equation*}
    \zeta = i\left[  2\sqrt{2} + \frac{\sqrt{2}\log(\tau + \tau_0) + \frac{1}{\sqrt{2}}\log 2 -\frac{1}{4\sqrt{2}} + 4\sqrt{2}\beta_1 + \xi }{\tau + \tau_0} \right], \quad f = (\tau + \tau_0)\rho,
\end{equation*}
then furnish the leading-order initial-value problem
\begin{equation*}
\sqrt{2}\frac{d\rho_0}{d \xi} = -\frac{d^2\rho_0}{d\xi^2} + \rho_0^2,   
\end{equation*}
\begin{equation*}
    \rho_0 \sim \frac{\sqrt{2}}{-\xi} + \frac{2\log-\xi}{(-\xi)^2} + \mathcal{O}\left(\frac{1}{(-\xi)^2} \right), \qquad \xi \to -\infty,
\end{equation*}
equivalent up to a rescaling to that discussed in \ref{sect:appode}, with the complex-singularity structure uncovered in Figures~\ref{fig:ODE_complex_phase_modulus} and~\ref{fig:ODE_complex_mod_sheet2_rotate} being present for $\xi = \mathcal{O}(1)$.

\subsection{Higher modes and non-generic blow up}\label{sec:non-generic} Two final issues are worth briefly noting here. Firstly, the generic blow-up behaviour should by definition be stable up to $t$ and $x$ translations, which respectively provide the zeroth and first modes in the associated linear stability problem. The constant $\tau_0$ embodies the second mode, leading to an algebraically decaying contribution of relative size $\tau_0/\tau$ as $\tau \to +\infty$. The remaining modes (providing a complete set) are exponentially decaying but are also (in this generic case) algebraically modulated for the same reason.  Specifically, on the inner scale one obtains within the linearisation
\begin{equation*}
    A_Me^{-(M-2)\tau/2}H_M(\eta/2)/\tau^M, \qquad M = 3, 4, \ldots, 
\end{equation*}
for constant $A_M$, as the leading order representation of the $M$th mode, the algebraic factor being determined from a solvability condition akin to that arising from (\ref{apeq:Fa0a1}) above, a calculation facilitated by the fact that only the two highest powers in the Hermite polynomials $H_M$ contribute to each solvability condition.

Finally, we revisit the non-generic form of blow up associated with the third expression in (\ref{apeq:eigenmodes}), whereby on the inner scale
\begin{equation}
    f \sim 1 - Ae^{-\tau}\left(1 - \eta^2 + \eta^4/12   \right)  \label{apeq:fnongen}
\end{equation}
for some constant $A>0$ that depends on the initial data. The outer variable is then $\zeta = \eta/e^{\tau/4}$ so that
\begin{equation*}
    f_0 + \frac{1}{4}\zeta \frac{df_0}{d\zeta} = f_0^2
\end{equation*}
giving 
\begin{equation}
    f_0 = \frac{1}{1 + A\zeta^4/12}, \label{apeq:nongen}
\end{equation}
and
\begin{equation}
    u(x,0) \sim \frac{12}{Ax^4}, \qquad x \to 0.  \label{eq:4thorderblowup}
\end{equation}
From (\ref{apeq:nongen}) it follows that four complex singularities collide at the blow-up time, giving further (and specifically complex-plane) insight into the non-generic nature of this scenario. We have given above an interpretation of its status as a borderline case, upon which we can now expand: given a solution symmetric about $x=0$ with $x=0$ being a local minimum at $t=0$, we can expect there initially to be four singularities nearest to the real axis, one in each quadrant. In the case of two-point blow-up the pairs in the left and right half planes will simultaneously collide with the real axis.
In the generic single-point blow-up  case the pairs in the upper and lower half planes  separately and simultaneously collide  on the imaginary axis prior to blow up, with singularities subsequently propagating down that axis to collide on the real axis at blow up; see Figure~\ref{fig:singscollide} and the right frame of Figure~\ref{fig:2peaksreal}.  In these figures, the local maximum of the real-line solution at $x = 0$ is converted to a maximum before the singularities collide and coalesce in the complex plane. 
The case (\ref{apeq:fnongen}) is exceptional in that the loss of the local minimum coincides with blow up.  Figures~\ref{fig:4singscollide} and \ref{fig:4thorderblowup} illustrate this case and confirm the fourth-order blow up in (\ref{eq:4thorderblowup}). 

\section{More general initial data} \label{sect:appicab}
\subsection{Preliminaries}

We now concisely generalise to the initial data
\begin{equation}
u(x,0) = \alpha \cos x + \beta,  \label{eq:abt0data}
\end{equation}
our purpose being twofold: firstly, to illustrate the modifications needed for more general initial data (indeed, some of our comments will concern arbitrary even and $2\pi$-periodic\footnote{Both of these restrictions are important in what follows.} initial data, namely
\begin{equation}
u(x,0) = U(x)  \label{eq:abgent0data}
\end{equation}
when such further generality does not meaningfully complicate the discussion) and, secondly and more specifically, to highlight a transition that can occur in the ultimate behaviour for $\beta < 0$.  Defining $\langle \psi \rangle$ to be the mean of a function $\psi(x)$ on $[-\pi, \pi]$, we have
\begin{equation*}
\frac{d}{dt}\langle u \rangle = \langle u^2 \rangle > 0
\end{equation*}
and if $\langle u \rangle \geq 0$ at $t=0$, then blow up necessarily ensues.   However, for $\langle u \rangle < 0$ two generic outcomes are possible, namely finite time blow up (as above) and
\begin{equation}
    u \sim -\frac{1}{t}, \qquad t \to \infty,
    \label{eq:heatdeath}
\end{equation}
for all $x$; these are divided by the separable solution
\begin{equation}
    u \sim A e^{-t}\cos x  \label{eq:heatsolint}
\end{equation}
of the heat equation, wherein $A$ is an arbitrary constant (this solution is singly unstable; as usual, there is a subsequent hierarchy of less and less stable such solutions, the next 
\begin{equation*}
u \sim Be^{-4t}\cos(2x)
\end{equation*}
being doubly unstable and separating the cases $A>0$ and $A<0$ in (\ref{eq:heatsolint}).

\begin{figure}[h!]
    \centering
   \includegraphics[width = 0.4\textwidth]{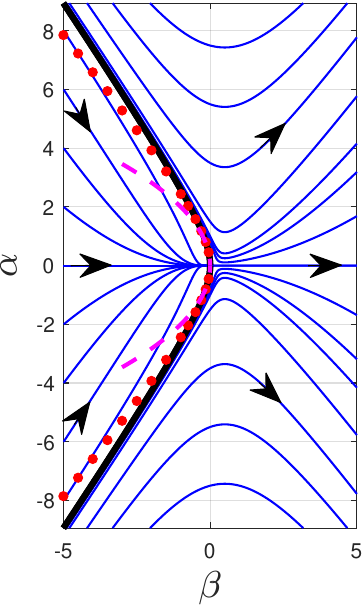}
   \caption{Trajectories of solutions of the two-mode approximation~(\ref{eq:NLH2modesystem}) to the NLH.  The trajectories are symmetric about the $\alpha$-axis because $\alpha \mapsto -\alpha$ corresponds to a translation of the solution ($x \mapsto x \pm \pi$). 
   The trajectories to the right of the exceptional trajectories (the thick black curves) are those of blow-up solutions, while the trajectories to the left asymptote to $0$ according to~(\ref{eq:heatdeath}).  The magenta-coloured dashed curve is defined by $\beta = -\alpha^2/4$, which, for $\vert \alpha \vert, \vert \beta \vert \ll 1$, is not only the borderline case between solutions that blow up and those that do not for the two-mode system (\ref{eq:NLH2modesystem}), but also for the NLH     (see section~\ref{sec:absmallamp}).  The red dots indicate numerically calculated values of $\alpha$ and $\beta$ in the initial data (\ref{eq:abt0data}) for the NLH that are at the threshold between solutions that blow up and solutions that tend to $0$ (i.e.~those that `extinguish' resulting in `heat death') according to (\ref{eq:heatdeath}).  For example, the red dot in the top-left corner corresponds to the values of $\alpha$ and $\beta$ used for the initial data of the `heat-death' solution in the left frame of Figure~\ref{fig:reddotfigs}.  A perturbation of these initial data leads to a blow-up solution, as in the right frame of Figure~\ref{fig:reddotfigs}.   As expected, the distance between the red dots and the exceptional trajectories of the two-mode system decreases as $\beta \to 0$.  
   }   
    \label{fig:nlh2modephaseplane}
\end{figure}

It is instructive to consider a two-mode approximation to the solution.  That is, if one sets
\begin{equation*}
    u = \beta(t) + \alpha(t)\cos x,
\end{equation*}
in the NLH and ignore higher-order modes (which gives accurate approximations provided $\vert \alpha \vert, \vert \beta \vert \ll 1$, a limit we shall consider below), then one obtains
\begin{eqnarray}
 \dot{\beta} &= \beta^2 + \frac{\alpha^2}{2}, \nonumber \\
 \dot{\alpha} &= (2\beta - 1)\alpha. \label{eq:NLH2modesystem}
\end{eqnarray}
Figure~\ref{fig:nlh2modephaseplane} shows the trajectories of this system. 


\subsection{Small-time limit}
In the case of (\ref{eq:abt0data}),
\begin{equation*}
u \sim \alpha\cos x + \beta + t\left( \beta^2 + \frac{1}{2}\alpha^2 + (2\beta-1)\alpha \cos x + \frac{1}{2}\alpha^2 \cos 2x  \right)
\end{equation*}
holds as $t \to 0$ sufficiently close to the real axis, this expansion disordering for $x = i\log t + \mathcal{O}(1)$, with the $\beta$ terms then playing no further leading-order role for small $t$, eliminating the need for any further discussion of this limit.

\subsection{$\vert \alpha \vert, \vert\beta\vert \ll 1$}\label{sec:absmallamp}
In capturing this case we set $\alpha = \epsilon a$, $\beta = \epsilon b$, and 
\begin{equation*}
    U(x) = \epsilon V(x)
\end{equation*}
in (\ref{eq:abgent0data}) with $0 < \epsilon \ll 1$.  Setting
\begin{equation*}
    u \sim \epsilon v(x,t), \qquad v \sim v_0(x,t) + \epsilon v_1(x,t)
\end{equation*}
yields
\begin{equation*}
    \frac{\partial v_0}{\partial t} = \frac{\partial^2 v_0}{\partial x^2}
\end{equation*}
so that
\begin{equation*}
    v_0 \sim b + a e^{-t}\cos x, \qquad t \to \infty,
\end{equation*}
where 
\begin{equation*}
    b = \langle V \rangle, \qquad a = 2\langle V \cos x \rangle
\end{equation*}
in the general case. Moreover,
\begin{equation*}
 \frac{\partial v_1}{\partial t} = \frac{\partial^2 v_1}{\partial x^2} + v_0^2   
\end{equation*}
so that
\begin{equation*}
    v_1 = b^2t + \frac{1}{4}a^2(1 - e^{-2t}) + 2abt e^{-t}\cos x + \frac{1}{4}a^2(e^{-2t} - e^{-4t} )\cos 2x
\end{equation*}
in the case of (\ref{eq:abt0data}), to which we now revert (similar results apply more generally as 
 $t \to \infty$ -- most notably, the deviation from spatially uniform is exponentially small in $t$ -- but we shall not pursue this further here).  Unless $b=0$ (the special case analysed in detail above, highlighting its important status) the above becomes non-uniform for $t = T/\epsilon$, for which
\begin{equation*}
v \sim \frac{k(\epsilon)}{1 - k(\epsilon)T} + \mathrm{ exponentially \: small\: terms }
\end{equation*}
where 
\begin{equation*}
k \sim b + \frac{1}{4}\epsilon a^2, \qquad \epsilon \to 0.
\end{equation*} 
Hence if $k >0$, finite-time blow up occurs, the previously exponentially small non-uniform terms ultimately coming into play in a manner akin to that described in section~\ref{sec:smallamp}, while (\ref{eq:heatdeath}) (heat death) holds for $k<0$.  The borderline case (\ref{eq:heatsolint}) applies for 
\begin{equation*}
b \sim -\frac{1}{4}\epsilon a^2,
\end{equation*}
and, in general, in the current limit the sign of $\langle V(x) \rangle$ is the key determinant of the ultimate outcome.

\subsection{$\vert \alpha \vert, \vert\beta\vert \gg 1$}

Here it suffices to consider the general case with
\begin{equation*}
U(x) = \frac{1}{\epsilon}V_0(x) + \frac{1}{\epsilon^{1/2}}V_1(x),
\end{equation*}
the reasons for the inclusion of the $V_1$ term being apparent shortly; we take $U(x)$ to have a single global maximum, at $x=0$.  Setting $t = \epsilon T$, $u = v/\epsilon$ and 
\begin{equation*}
v \sim v_0(x,T) + \epsilon^{1/2}v_1(x,T) 
\end{equation*}
gives
\begin{equation*}
v_0 = \frac{V_0(x)}{1 - TV_0(x)}, \qquad v_1 = \frac{V_1(x)}{(1 - TV_0(x))^2}
\end{equation*}
Hence (\ref{eq:heatdeath}) ensues if $V_0(0)<0$, while finite-time blow up (as analysed in section~\ref{sec:largeamp}) occurs for $V_0(0) > 0$, a point worth emphasis being that, while the borderline criterion in section~\ref{sec:absmallamp} is dominated by the average of $V(x)$, here it is the maximum of $V(x)$ that is critical.  The borderline is thus associated with $V_0(0) = 0$; we then set
\begin{equation*}
V_0(x) \sim \nu_0 x^2, \qquad x \to 0, \qquad V_1(0) = \mu_1,
\end{equation*} 
$\nu_0 > 0$, $\mu_1 = \mathcal{O}(1)$ then capturing the transition via the rescalings 
\begin{equation}
T = \epsilon^{-1/2}\tau, \qquad x = \epsilon^{-1/4}\xi, \qquad v = \epsilon^{1/2}\phi \label{eq:n6}
\end{equation}
whereby
\begin{eqnarray}
& \frac{\partial \phi_0}{\partial \tau} = \frac{\partial^2 \phi_0}{\partial \xi^2} + \phi_0^2, \nonumber \\
& \mathrm{at\: } \tau = 0, \qquad \phi_0 = \mu_1 - \nu_0 \xi^2, \qquad \mathrm{for\: }  \xi = \mathcal{O}(1),
\label{eq:n7}
\end{eqnarray}
this leading-order formulation applying on the whole real line (and amenable to continuation into the complex plane, of course) because of the spatial rescaling in (\ref{eq:n6}).  By rescaling, one can set $\nu_0=1$ in (\ref{eq:n7}), leaving a single parameter $\mu_1$, and we conjecture that the initial-value problem (\ref{eq:n7}) then has the following properties.
\begin{itemize}
    \item[(I)] $\phi_0 \sim -\frac{1}{\tau}$ as $\vert \xi \vert \to \infty$ with
\begin{equation*}
\phi_0 \sim - \frac{\xi^2}{1 + \xi^2 \tau} \qquad \mathrm{as\: } \tau \to 0^{+} \mathrm{\: with\: } \xi = \mathcal{O}\left(\tau^{-1/2}\right).
\end{equation*} 
\item[(II)] For $\mu_1<\mu_*$ for some $\mu_*<0$, $\phi_0 \sim -\frac{1}{\tau}$ as $\tau \to + \infty$ for all $\xi$.
\item[(III)] For $\mu_1 > \mu_*$ finite-time blow up occurs.
\item[(IV)] For $\mu_1 = \mu_*$,
\end{itemize}
\begin{equation}
    \phi_0 \sim \frac{1}{\tau}f(\eta), \qquad \eta = \frac{\xi}{\tau^{1/2}} \qquad \mathrm{ as\: } \tau \to \infty \qquad \mathrm{with\: } \qquad \eta = \mathcal{O}(1),  \label{eq:n8}
\end{equation}
\begin{itemize}
\item[]
where $f(\eta)$ satisfies the boundary value problem
\begin{equation*}
-f -\frac{1}{2}\eta \frac{d f}{d \eta} = \frac{d^2 f}{d\eta^2} + f^2,\qquad
f'(0) = 0, \qquad
f \to -1 \mathrm{\: as\: } \vert \eta \vert \to \infty,
\end{equation*}
this similarity reduction being singly unstable.
\end{itemize}
We leave the calculation of the critical value $\mu_*$ as a worthwhile open problem.  For $\tau = \mathcal{O}(\epsilon^{-1/2})$ (i.e. $t = \mathcal{O}(1)$) the formulation (\ref{eq:n8}) will break down, spatial periodicity coming into play and necessitating a further (higher-order) refinement of the borderline.

\subsection{$\vert \alpha \vert \ll 1$, $\beta = \mathcal{O}(1)$}

This is the final regime that seems to warrant separate treatment here.  We set $\alpha = \epsilon$ so that
\begin{equation*}
u_0 = \frac{\beta}{1 - \beta t}, \qquad u_1 = \frac{1}{(1 - \beta t)^2}e^{-t}\cos x
\end{equation*}
for $t < 1/\beta$; incidentally, the algebraic convenience of such calculations provides a motivation for our specific choice of initial conditions.  The subsequent timescale sets 
\begin{equation*}
    t = \frac{1}{\beta} + \epsilon T, \qquad u = \frac{1}{\epsilon}v
\end{equation*}
with
\begin{equation*}
    v_0 = \frac{1}{-T - \frac{1}{\beta^2}e^{-\beta}\cos x}
\end{equation*}
so blow up occurs at
\begin{equation*}
    T \sim -\frac{1}{\beta^2}e^{-1/\beta}
\end{equation*}
having non-monotonic dependence upon $\beta$.  The final approach to blow up  parallels that of section~\ref{sec:smallamp}.

The main observation arising from this regime relates to the exponential  dependence on $\beta$, which points towards the role of exponentially small terms in the preceding analysis.

\subsection{Complex-plane singularities}

The behaviour of the singularities closest to the real axis has been explored in detail above for $\beta = 0$ and no qualitatively new effects arise in the finite-time blow up case for $\beta \neq 0$.  When (\ref{eq:heatdeath})  or (\ref{eq:heatsolint}) apply instead, the singularities unsurprisingly move away from the real axis, at least for sufficiently large $t$.  For (\ref{eq:heatsolint}) the behaviour
\begin{equation}
    \sigma(t) =  t + \mathcal{O}(1), \qquad t \to \infty,  \label{eq:n9}
\end{equation}
follows immediately from the criterion that $u$ be of $\mathcal{O}(1)$; for (\ref{eq:heatdeath}) we need the correction term
\begin{equation*}
    u \sim -\frac{1}{t} + \frac{A}{t^2}e^{-t}\cos x, \qquad t \to \infty,
\end{equation*}
in inferring that
\begin{equation}
    \sigma(t) = t + 2\log t + \mathcal{O}(1), \qquad t \to \infty,  \label{eq:n11}
\end{equation}
which is consistent with the linear increase of the blue curve in Figure~\ref{fig:heatdeathsing}.
It is noteworthy that (\ref{eq:n9}) and (\ref{eq:n11})  coincide at leading order.








\end{document}